\documentclass[reqno]{amsart}
\usepackage{amsmath,amsrefs,amssymb,hyperref, xcolor}

\numberwithin{equation}{section}

\DeclareMathOperator{\divergence}{div}
\DeclareMathOperator{\Span}{span}
\DeclareMathOperator{\Vol}{Vol}
\DeclareMathOperator{\Scal}{Scal}

\DeclareMathOperator{\loc}{loc}
\DeclareMathOperator{\dv}{dv}
\DeclareMathOperator{\dy}{dy}
\DeclareMathOperator{\dz}{dz}
\DeclareMathOperator{\tr}{tr}
\DeclareMathOperator{\dsigma}{d\sigma}
\DeclareMathOperator{\dr}{dr}
\DeclareMathOperator{\dmu}{d\mu}
\DeclareMathOperator{\dist}{d}
\DeclareMathOperator{\B}{B}
\DeclareMathOperator{\bigO}{O}
\DeclareMathOperator{\smallo}{o}
\newcommand{\R}{\mathbb{R}}
\renewcommand{\S}{\mathbb{S}}
\newcommand{\BB}{\mathbb{B}}
\newcommand{\N}{\mathbb{N}}

\newcommand{\<}{\left<}
\renewcommand{\>}{\right>}
\renewcommand{\[}{\left[}
\renewcommand{\]}{\right]}
\renewcommand{\(}{\left(}
\renewcommand{\)}{\right)}

\newtheorem{theorem}{Theorem}[section]
\newtheorem{lemma}[theorem]{Lemma}

\begin{document}

\title[Nonexistence of minimizers for the second conformal eigenvalue]{Nonexistence of minimizers for the second conformal eigenvalue near the round sphere\\in low dimensions}

\author{Bruno Premoselli}

\address{Bruno Premoselli, Universit\'e Libre de Bruxelles, Service d'analyse, CP 218, Boulevard du Triomphe, B-1050 Bruxelles, Belgique}
\email{bruno.premoselli@ulb.be}

\author{J\'er\^ome V\'etois}

\address{J\'er\^ome V\'etois, Department of Mathematics and Statistics, McGill University, 805 Sherbrooke Street West, Montreal, Quebec H3A 0B9, Canada}
\email{jerome.vetois@mcgill.ca}


\thanks{To appear in {\it Journal of Differential Geometry}}

\thanks{The second author was supported by the NSERC Discovery Grant RGPIN-2022-04213.}

\date{September 20, 2026}

\begin{abstract}
We consider the problem of minimizing the second conformal eigenvalue of the conformal Laplacian in a conformal class of metrics with renormalized volume. We prove, in dimensions $n\in\left\{3,\dotsc,10\right\}$, that a minimizer for this problem does not exist for metrics sufficiently close to the round metric on the sphere. This is in striking contrast with the situation in dimensions $n \ge 11$, where Ammann and Humbert~\cite{AH} obtained the existence of minimizers for the second conformal eigenvalue on any smooth closed non-locally conformally flat manifold. As a byproduct of our techniques, we also obtain a lower bound on the energy of sign-changing solutions of the \nobreak Yamabe equation in dimensions 3, 4 and 5, which extends a result obtained by Weth~\cite{W} in the case of the round sphere. 
\end{abstract}

\maketitle

\section{Introduction and main result}

We let $\(M,g\)$ be a smooth closed Riemannian manifold of dimension $n\ge3$. We let $\[g\]$ be the conformal class of the metric $g$. For each metric $\hat{g}\in\[g\]$, we denote by $L_{\hat{g}}$ the conformal Laplacian of $\(M,\hat{g}\)$, i.e.
$$L_{\hat{g}}:=\Delta_{\hat{g}}+c_n\Scal_{\hat{g}},$$
where $\Delta_{\hat{g}} := - \divergence_{\hat{g}}\(\nabla \cdot\)$ is the Laplace--Beltrami operator of $\(M,\hat{g}\)$, $\Scal_{\hat{g}}$ is the scalar curvature of $\(M,\hat{g}\)$ and $c_n:=\frac{n-2}{4\(n-1\)}$. Since $M$ is closed, for each $\hat{g}\in\[g\]$ the eigenvalues of $L_{\hat{g}}$ form a nondecreasing sequence $\(\lambda_k\(L_{\hat{g}}\)\)_{k\in\N}$ such that
$$ \lambda_1\(L_{\hat{g}}\)< \lambda_2\(L_{\hat{g}}\) \le \cdots \le \lambda_k\(L_{\hat{g}}\) \le \cdots \to \infty . $$
For each $k\in\N$, the $k$-th conformal eigenvalue of $\(M,\[g\]\)$ is defined as 
$$\Lambda_k\(M,\[g\]\):=\inf_{\hat{g}\in\[g\]} \( \lambda_k\(L_{\hat{g}}\)\Vol\(M,\hat{g}\)^{\frac{2}{n}}\),$$
where $\lambda_k\(L_{\hat{g}}\)$ is the $k$-th eigenvalue of $L_{\hat{g}}$ and $\Vol\(M,\hat{g}\)$ is the volume of $\(M,\hat{g}\)$. This invariant was first introduced and studied by Ammann and Humbert~\cite{AH} and further studied by El Sayed~\cite{ElS}. It is not difficult to see that if $\Lambda_1\(M,\[g\]\)\ge0$, then $\Lambda_1\(M,\[g\]\)$ coincides with the classical Yamabe invariant, i.e.
$$\Lambda_1\(M,\[g\]\)=c_n\inf_{\hat{g}\in\[g\]}\(\Vol\(M,\hat{g}\)^{\frac{2-n}{n}}\int_M\Scal_{\hat{g}}\dv_{\hat{g}}\),$$
where $\dv_{\hat{g}}$ is the volume element of $\(M,\hat{g}\)$. In this paper, we consider the case of metrics with positive Yamabe invariant. In this case, it follows from the work of Trudinger~\cite{Tr}, Aubin~\cite{A} and Schoen~\cite{Sc1} that $\Lambda_1\(M,\[g\]\)$ is always attained by some smooth positive function, and, moreover,
$$\Lambda_1\(M,\[g\]\)\le\Lambda_1\(\S^n,\[g_0\]\)$$
with equality if and only if $\(M,g\)$ is conformally equivalent to the round $n$-sphere $\(\S^n,g_0\)$.  

\smallskip
In this paper, we focus on the case where $k=2$ and $\Lambda_1\(M,\[g\]\)>0$. In this case, Ammann and Humbert~\cite{AH} obtained
\begin{equation}\label{IntroEq1}
2\Lambda_1\(M,\[g\]\)^{\frac{n}{2}}\le\Lambda_2\(M,\[g\]\)^{\frac{n}{2}}\le\Lambda_1\(M,\[g\]\)^{\frac{n}{2}}+\Lambda_1\(\S^n,\[g_0\]\)^{\frac{n}{2}}.
\end{equation}
They also obtained that $\Lambda_2\(M,\[g\]\)$ is attained provided the second inequality in \eqref{IntroEq1} is strict, i.e.
\begin{equation} \label{IntroEq2}
 \Lambda_2\(M,\[g\]\)^{\frac{n}{2}} <\Lambda_1\(M,\[g\]\)^{\frac{n}{2}}+\Lambda_1\(\S^n,\[g_0\]\)^{\frac{n}{2}}
 \end{equation}
and that if \eqref{IntroEq2} is satisfied, then $\Lambda_2\(M,\[g\]\)$ is attained at a ``generalized'' metric $\hat{g}:=\left|u\right|^{2^*-2}g$ for some $u\in C^{3,\vartheta}\(M\)$ for some $\vartheta<2^*-2$ (see Section~\ref{Sec2} for more details). Ammann and Humbert also obtained in \cite{AH} that, by test-functions computations, if $n\ge 11$ and $\(M,g\)$ has non-vanishing Weyl tensor somewhere, then \eqref{IntroEq2} is satisfied and that, for each $n \ge 3$, $\Lambda_2\(\S^n, \[g_0\]\) = 2^{\frac{2}{n}}\Lambda_1(\S^n, \[g_0\])  $ is never attained. 

\smallskip
Except for the trivial case of the round sphere $\(\S^n, g_0\)$, the existence of minimizers for $\Lambda_2\(M,\[g\]\)$ in dimensions 3 to 10 is an open problem. Our main result provides a partial negative answer to this question: 

\begin{theorem}\label{Th1}
Assume that $3\le n\le10$. Let $m_0 \in \N$. Then there exist $\delta\in\(0,\infty\)$ and an integer $m \ge m_0$ such that, for every smooth metric $g$ on $\S^n$, if $\left\|g-g_0\right\|_{C^m\(\S^n\)}<\delta$, then
$$\Lambda_2\(\S^n,\[g\]\)^{\frac{n}{2}}=\Lambda_1\(\S^n,\[g\]\)^{\frac{n}{2}}+\Lambda_1\(\S^n,\[g_0\]\)^{\frac{n}{2}}$$
and $\Lambda_2\(\S^n,\[g\]\)$ is not attained by any generalized metric.
\end{theorem}

We recall that every smooth metric on $\S^n$ which is not conformally equivalent to $g_0$ is also not locally conformally flat. Theorem~\ref{Th1} establishes a striking dichotomy between the case where $3 \le n \le 10$ and the case where $n \ge 11$. First, when $3 \le n \le 10$, Theorem~\ref{Th1} has to be understood as a perturbative nonexistence result of extremals for $\Lambda_2\(\S^n,\[g\]\)$ for $g$ close to the round metric $g_0$, the analogue of which fails when $n \ge 11$ by the results of \cite{AH}. Second, Theorem~\ref{Th1} establishes that, when $3 \le n \le 10$, \eqref{IntroEq2} cannot be guaranteed by solely enforcing local geometric assumptions on $g$. This again strongly contrasts with the results of \cite{AH} in dimensions $n\ge11$, where a minimizer for $\Lambda_2\(M,\[g\]\)$ is proven to exist for any $g$ which is not locally conformally flat. Determining whether \eqref{IntroEq2} holds true and whether there exist extremals for $\Lambda_2\(M,\[g\]\)$ when $3 \le n \le 10$ therefore requires new ideas. Theorem 1.1 can be understood as a step forward in this direction : when $\(M,g\)=\(\S^n,g\)$ it reveals that a necessary condition for (1.2) to hold is that $g$ is sufficiently far from the round metric $g_0$ in a strong sense. Not being conformally diffeomorphic to $g_0$, in particular, is not enough, which is very surprising in view of the definition of $\Lambda_2$. How Theorem 1.1 may adapt on a general manifold is still unclear, but it seems to hint that global infomation on $\(M,g\)$ is needed to obtain (1.2).

\smallskip
Eigenvalue optimization problems in conformal classes have attracted a lot of attention in recent years. When $n \ge3$ and in the case of the conformal Laplacian $L_g$ which we consider here, the invariants $\Lambda_k\(M,\[g\]\)$ that we define are the only meaningful ones when $\Lambda_1\(M,\[g\]\) >0$. Indeed, it is for instance proven by Ammann and Jammes~\cite{AJ} that, if $\Lambda_k\(M,\[g\]\)>0$, then
$$\sup_{\hat{g}\in\[g\]}\(\lambda_k\(L_{\hat{g}}\)\Vol\(M,\hat{g}\)^{\frac{2}{n}}\)=\infty.$$
On the other hand, in the case where $\Lambda_k\(M,\[g\]\)<0$, it is natural to replace the infimum in the definition of $\Lambda_k$ by a supremum since otherwise $\Lambda_k\(M,\[g\]\)=-\infty$ (see Proposition~8.1 in \cite{AH}). This therefore leads to a maximization problem, which is very different in nature. We refer to the work of Gursky and P\'erez-Ayala~\cite{GP} where this problem is studied in the case where $k=2$. Most of previous work on eigenvalue optimization problems in conformal classes of closed manifolds of dimension larger than or equal to $2$ concern the Laplace--Beltrami operator $\Delta_g$, for which again only the maximization problem is interesting. In dimensions $n \ge 3$, the maximization of conformal eigenvalues of $\Delta_g$ was recently investigated by P\'etrides~\cite{Pe3}. In dimension $2$, this problem was investigated by many authors. In this case, we refer for instance to the work of Nadirashvili and Sire~\cite{NS}, Petrides~\cites{Pe1,Pe2}, Matthiesen and Siffert~\cite{MS} and Karpukhin and Stern~\cites{KS1,KS2}. To the best of our knowledge, another remarkable feature of Theorem~\ref{Th1} is that it is the first nonexistence result of extremals for conformal eigenvalues of any kind (for any dimension $n \ge 2$ and any of the operators $\Delta_g$ and $L_g$) for metrics that are not conformal to the round metric on $\S^n$. 

\smallskip
The structure of the paper is as follows. In Section~\ref{Sec2}, we discuss the connection between the second conformal eigenvalue and sign-changing solutions of the Yamabe equation of lowest energy. We also state a stronger result than Theorem~\ref{Th1}, in dimensions 3, 4 and 5, namely Theorem~\ref{Th2}. Section~\ref{Sec3} is devoted to a sharp bubbling analysis of sign-changing solutions of the Yamabe equation whose energies converge to $\Lambda_2\(\S^n,\[g_0\]\)^{\frac{n}{2}}$. We prove bubble-tree convergence results as well as sharp pointwise asymptotics. We then prove Theorems~\ref{Th1} and \ref{Th2} in Section~\ref{Sec4}.

\smallskip\noindent
\textbf{Acknowledgements:} The authors are truly indebted to the anonymous referees for their wise and important comments which greatly improved the quality and the presentation of the paper.

\section{The second conformal eigenvalue and sign-changing solutions of the Yamabe equation}\label{Sec2}

The second conformal eigenvalue of the conformal Laplacian has a strong connection with sign-changing solutions of the Yamabe equation
\begin{equation}\label{Sec2Eq}
L_gu=\left|u\right|^{2^*-2}u\quad\text{in }M,
\end{equation}
where $2^*:=\frac{2n}{n-2}$ is the critical Sobolev exponent. Indeed, Ammann and \nobreak Humbert~\cite{AH} proved that if $\Lambda_1\(M,\[g\]\)\ge0$ and $\Lambda_2\(M,\[g\]\)$ is attained, then there exists a sign-changing function $u\in L^{2^*}\(M\)$ such that the ``generalized'' metric $\hat{g}:=\left|u\right|^{2^*-2}g$ satisfies 
$$\lambda_2\(L_{\hat{g}}\) = \Lambda_2\(M,\[g\]\)\quad\text{and}\quad\Vol\(M,\hat{g}\)=1,\text{ i.e.}\int_M\left|u\right|^{2^*} \dv_g = 1$$
(see \cite{AH}*{Section 3.2} for the rigorous definition of $\lambda_2\(L_{\hat{g}}\)$ when $u\in L^{2^*}\(M\) \setminus \left\{0\right\}$). It is also shown in \cite{AH} that $u$ is a second ``generalized'' eigenvector associated to $\lambda_2\(L_{\hat{g}}\)$, which implies that $u\in C^{3,\vartheta}\(M\)$ for some $\vartheta <2^*-2$ and, up to a renormalization factor, $u$ can be made into a sign-changing solution of \eqref{Sec2Eq} with energy
$$\int_M\left|u\right|^{2^*}\dv_g=\Lambda_2\(M,\[g\]\)^{\frac{n}{2}}.$$
Furthermore, in this case $u$ is a sign-changing solution of \eqref{Sec2Eq} of least energy among all sign-changing solutions, and it has exactly two nodal domains. We again refer to \cite{AH}*{Section~3} for more details. As mentioned in the introduction, \eqref{IntroEq2} is sufficient to ensure that $\Lambda_2\(M,\[g\]\)$ is attained. Equation \eqref{Sec2Eq} has to be understood as the Euler-Lagrange equation for minimizers of $\Lambda_2\(M,\[g\]\)$. Its hidden meaning is that, for a generalized metric $\hat{g}:=\left|u\right|^{2^*-2}g$ attaining $\Lambda_2\(M,\[g\]\)$, $\lambda_2\(L_{\hat{g}}\)$ is simple. This unusual feature for an eigenvalue optimization problem is a direct consequence of the definition of $\Lambda_2\(M,\[g\]\)$ as an infimum (see for instance \cite{GP}*{Remark~6.1} for a detailed explanation). 

\smallskip
A consequence of Theorem~\ref{Th1} is that, for every smooth metric $g$ on $\S^n$ sufficiently close to $g_0$ in $C^m\(\S^n\)$ for some sufficiently large $m\in\N$, there does not exist any sign-changing solution $u$ of \eqref{Sec2Eq} such that 
$$\int_{\S^n}\left|u\right|^{2^*} \dv_g \le \Lambda_1\(\S^n,\[g\]\)^{\frac{n}{2}}+\Lambda_1\(\S^n,\[g_0\]\)^{\frac{n}{2}}.$$
Theorem~\ref{Th1} thus gives a lower bound on the energy of sign-changing solutions of the Yamabe equation \eqref{Sec2Eq} on the sphere of dimension lower than or equal to 10 when equipped with metrics sufficiently close to the round metric. In fact, in dimensions 3, 4 and 5, we obtain a stronger result: 

\begin{theorem}\label{Th2}
Assume that $n\in\left\{3,4,5\right\}$. Let $m_0 \in \N$. Then there exist $\delta,\varepsilon\in\(0,\infty\)$ and an integer $m \ge m_0$ such that, for every smooth metric $g$ on $\S^n$, if $\left\|g-g_0\right\|_{C^m\(\S^n\)}<\delta$, then the energy of every sign-changing solution of the Yamabe equation
$$L_gu=\left|u\right|^{2^*-2}u\quad\text{in }\S^n$$
is greater than $2\Lambda_1\(\S^n,\[g_0\]\)^{\frac{n}{2}}+\varepsilon$.
\end{theorem}

Theorem~\ref{Th2} extends a result obtained by Weth~\cite{W} in the exact case of the round sphere. Notice, however, that the result of Weth~\cite{W} holds for all dimensions $n \ge 3$. This is specific to the exact case $g = g_0$. Indeed, as shown by the results of Ammann and Humbert~\cite{AH}, at least in the case where $n \ge 11$, Theorem~\ref{Th2} is false when $g$ is not conformal to $g_0$.

\smallskip
We prove Theorems~\ref{Th1} and \ref{Th2} in the next two sections. By using \eqref{Sec2Eq} together with a contradiction argument, which is explained in details at the beginning of Section~\ref{Sec3}, the proof of Theorems~\ref{Th1} and \ref{Th2} amounts to ruling out the existence of sequences $\(u_k\)_{k \in \N}$ of sign-changing solutions of \eqref{Sec2Eq} with $g=g_k$, where $g_k$ converges to $g_0$ in $C^m\(\S^n\)$ as $k\to\infty$ for all $m\in\N$ such that the energies of $\(u_k\)_k$ converge to $2 \Lambda_1\(\S^n, \[g_0\]\)^{\frac{n}{2}}$. In dimensions 6 to 10, we assume in addition that, for each $k\in\N$, $\Lambda_2 \(\S^n, \[g_k\]\)$ is attained by the generalized metric $\left|u_k\right|^{2^*-2}g_k$. The proof of Weth~\cite{W} in the case where $g_k=g_0$ for all $k\in\N$ relies on the symmetries of $\(\S^n, g_0\)$ and uses the action of the conformal group of the sphere on the sign-changing solutions of the Yamabe equation. In our setting, however, the metrics $\(g_k\)_k$ do not have any symmetries in general, and we need to perform a much finer asymptotic analysis of $\(u_k\)_k$. As a first result, in Lemma~\ref{Lem1}, we prove that $\(u_k\)_{k}$ behaves like the difference between two positive solutions of the Yamabe equation on the sphere (see \eqref{Lem1Eq5}). This amounts to say that $\Lambda_2 \(\S^n, \[g_k\]\)$ is asymptotically attained by the disjoint union of two round spheres. The rest of Section~\ref{Sec3} is devoted to obtaining sharp pointwise estimates for the blow-up of $\(u_k\)_{k}$, which, in particular, captures the local geometry of $\(g_k\)_{k}$. This refined blow-up analysis is based on iterated estimates and makes crucial use of arguments previously developed in the context of positive solutions (see the work of Chen and Lin~\cite{CL}, Schoen~\cites{Sc2,Sc3}, Li and Zhu~\cite{LZhu}, Druet~\cite{Dr}, Marques~\cite{Mar}, Li and Zhang~\cites{LZha1,LZha2} and Khuri, Marques and Schoen~\cite{KMS}; see also the counterexamples in high dimensions of Brendle~\cite{B} and Brendle and Marques~\cite{BM1} and the survey article by Brendle and Marques~\cite{BM2}) and more recently in the context of sign-changing solutions (see Premoselli~\cite{Pr} and Premoselli and V\'etois~\cites{PV1,PV2,PV3}). We finally prove Theorems~\ref{Th1} and \ref{Th2} in Section~\ref{Sec4} by using the analysis developed in Section~\ref{Sec3}. 

\smallskip
Although similar at first glance, the settings of Theorems~\ref{Th1} and \ref{Th2} are quite different from that of the celebrated compactness result of Khuri, Marques and Schoen~\cite{KMS}. We point out two crucial differences: first, since the functions $\(u_k\)_{k}$ change sign, the concentration points are neither isolated nor simple; second, since $g_{k} \to g_0$ as $k\to\infty$ in $C^m\(\S^n\)$ for all $m\in\N$, the Riemannian masses of the metrics $\(g_k\)_{k}$ converge to $0$ at any point of $\S^n$, so that no local sign restriction argument is available to rule out blow-up. Therefore, and unlike in \cite{KMS}, our contradiction does not originate from a local sign restriction due to the Positive Mass Theorem, which we do not use in this paper. In dimensions 3, 4 and 5, instead, we obtain a contradiction by means of a Pohozaev-type identity in a region where we observe that a large virtual mass is created solely by the interaction between the two bubbles. In dimensions 6 to 10, the Pohozaev-type identity is not sufficient to conclude since additional lower-order terms appear which involve more of the geometry of the metrics $\(g_k\)_{k}$ at the concentration points. In this case, we still manage to obtain a sharp asymptotic estimate on $\Lambda_1\(\S^n,\[g_k\]\)$ (see \eqref{Lem9Eq2}). We then obtain a contradiction with this estimate by doing another estimation of $\Lambda_1\(\S^n,\[g_k\]\)$ based on a better family of test-functions, which construction again relies on the analysis of Section~\ref{Sec3}. The contradiction when $6 \le n \le 10$ thus really comes from the minimality of $\Lambda_2\(\S^n,\[g_k\]\)$. 

\smallskip
There is an abundant literature on sign-changing solutions of the Yamabe equation. In addition to the above-mentioned articles of Ammann and  Humbert~\cite{AH},  El Sayed~\cite{ElS} and  Gursky and  Perez-Ayala~\cite{GP}, we also refer on this topic to the historic work of Ding~\cite{Di} and the more recent work of Clapp~\cite{C}, Clapp, Pistoia and Weth~\cite{CPW}, del~Pino,  Musso, Pacard and  Pistoia~\cites{dPMPP1,dPMPP2}, Fernandez, Palmas and  Petean~\cite{FPP}, Fernandez and  Petean~\cite{FP}, Medina,  Musso and Wei~\cite{MMW}, Musso and  Medina~\cite{MM}, Musso and Wei~\cite{MW}, Premoselli and V\'etois~\cite{PV1} and Weth~\cite{W} in the case of the sphere, Clapp and Fernandez~\cite{CF} in the case of manifolds satisfying some symmetry assumptions and Clapp, Pistoia and  Tavares~\cite{CPT}, Premoselli and Robert~\cite{PR} and Premoselli and  V\'etois~\cite{PV3} in the case of more general manifolds. Other results in this spirit have been obtained for some classes of sign-changing solutions of equations with different potential functions than the Yamabe equation (see our previous articles~\cites{PV1,PV2,V}). Theorems~\ref{Th1} and \ref{Th2} can also be seen as a continuation of our study, initiated  in \cite{PV3}, of minimal energy sign-changing blowing-up solutions of the Yamabe equation.

\section{Asymptotic analysis and symmetry estimates}\label{Sec3}

The section and the next are devoted to the proofs of Theorems~\ref{Th1} and \ref{Th2}. We begin with recalling some well-known facts about solutions of the Yamabe equation in $\R^n$ and $\S^n$. By splitting the solutions into positive and negative parts, it is easy to see that there does not exist any sign-changing solution of the Yamabe equation
\begin{equation}\label{Sec3Eq1}
L_{g_0}u=\left|u\right|^{2^*-2}u\quad\text{in }\S^n
\end{equation}
with energy smaller than or equal to $2\Lambda_1\(\S^n,\[g_0\]\)^{\frac{n}{2}}$. Moreover, according to the classification result of Obata~\cite{O}, every nonzero nonnegative solution of \eqref{Sec3Eq1} has energy equal to $\Lambda_1\(\S^n,\[g_0\]\)^{\frac{n}{2}}$ and is either constant or of the form
$$u=\(\frac{2\sqrt{n\(n-2\)}\mu}{2\mu^2+\big(4-\mu^2\big)\big(1-\cos\(\dist_{g_0}\(\cdot,x\)\)\big)}\)^{\frac{n-2}{2}}$$
for some $x\in\S^n$ and $\mu\in\(0,2\)$. We also recall that, letting $\xi$ be the Euclidean metric on $\R^n$ and $\Delta_\xi:=-\sum_{i=1}^n\partial_{y_i}^2$, the stereographic projection gives a bijection between the solutions of \eqref{Sec3Eq1} and the solutions $u$ of the equation
\begin{equation}\label{Sec3Eq2}
\Delta_\xi u=\left|u\right|^{2^*-2}u\quad\text{in }\R^n,
\end{equation}
which belong to the energy space $D^{1,2}\(\R^n\)$ defined as the closure of $C^\infty_c\(\R^n\)$ with respect to the norm $\left\| \nabla \cdot \right\|_{L^2\(\R^n\)}^2$. By using this bijection, we obtain that every nonzero nonnegative solution of \eqref{Sec3Eq2} has energy equal to $\Lambda_1\(\S^n,\[g_0\]\)^{\frac{n}{2}}$ and is of the form
$$u=\(\frac{\sqrt{n\(n-2\)}\widetilde\mu}{\widetilde\mu^2+\left|\cdot-y\right|^2}\)^{\frac{n-2}{2}}$$
for some $y\in\R^n$ and $\widetilde\mu\in\(0,\infty\)$, and that there does not exist any sign-changing solution $u\in D^{1,2}\(\R^n\)$ of \eqref{Sec3Eq2} with energy smaller than or equal to $2\Lambda_1\(\S^n,\[g_0\]\)^{\frac{n}{2}}$. It is well-known that the positive solutions of \eqref{Sec3Eq2} are the extremals for the Sobolev inequality in $\R^n$, i.e.
$$ \Lambda_1\(\S^n,\[g_0\]\) = \inf_{v \in C^\infty_c(\R^n) \setminus \left\{0\right\}} \frac{\displaystyle\int_{\R^n}\left|\nabla v\right|^2\dy}{\displaystyle\(\int_{\R^n}\left|v\right|^{2^*} \dy\)^{\frac{n-2}{n}}} ,$$
where $\dy$ is the volume element of $\(\R^n,\xi\)$.

\smallskip
We prove Theorems~\ref{Th1} and \ref{Th2} by contradiction. From now until the end of the paper, we assume that $3\le n\le10$. By assuming the negation of the statement of Theorem~\ref{Th1} and applying a simple diagonal argument, we assume that there exists a sequence $\(g_k\)_{k\in\N}$ of smooth metrics on $\S^n$ such that, for each $m\in\N$, $g_k\to g_0$ in $C^m\(\S^n\)$ as $k\to\infty$ and for each $k\in\N$, there exists a sign-changing solution $u_k\in C^{3,\vartheta}\(\S^n\)$, $\vartheta<2^*-2$, of the Yamabe equation 
\begin{equation}\label{Sec3Eq3}
L_{g_k}u_k=\left|u_k\right|^{2^*-2}u_k\quad\text{in }\S^n.
\end{equation}
In the case where $3 \le n \le 5$, in view of Theorem~\ref{Th2}, we only assume that
\begin{equation}\label{Sec3Eq4}
\limsup_{k\to\infty}\int_{\S^n}\left|u_k\right|^{2^*}\dv_{g_k}\le2\Lambda_1\(\S^n,\[g_0\]\)^{\frac{n}{2}}.
\end{equation}
In the case where $6\le n\le10$, in view of Theorem~\ref{Th1}, we assume moreover that, for each $k\in\N$, $\Lambda_2\(\S^n,\[g_k\]\)$ is attained by the generalized metric $\left| u_k \right|^{2^*-2} g_k$ and
\begin{equation}\label{Sec3Eq5}
\int_{\S^n}\left|u_k\right|^{2^*}\dv_{g_k}=\Lambda_2\(\S^n,\[g_k\]\)^{\frac{n}{2}},
\end{equation}
which implies that $u_k$ is a sign-changing solution of \eqref{Sec3Eq3} of least-energy among all sign-changing solutions. We point out in passing that, since $u_k$ satisfies \eqref{Sec3Eq3}, the celebrated results of Hardt and Simon~\cite{HS} gives that $u_k$ vanishes on a set of measure zero and is therefore admissible in the definition of $\Lambda_2\(\S^n,\[g_k\]\)$ (see \cite{AH}*{Section~3}). By putting together \eqref{IntroEq1} and \eqref{Sec3Eq4}, we obtain 
\begin{equation}\label{Sec3Eq6}
\lim_{k\to\infty}\int_{\S^n}\left|u_k\right|^{2^*}\dv_{g_k}=2\Lambda_1\(\S^n,\[g_0\]\)^{\frac{n}{2}}.
\end{equation}
A first simple remark which follows from the previous discussion and from the convergence $g_k \to g_0$ in $C^m\(\S^n\)$ for all $m \in \N$ as $k \to \infty$ is that the sequence $\(u_k\)_k$ blows up as $k \to \infty$, i.e.
\begin{equation} \label{Sec3Eq7}
\left\|u_k \right\|_{L^\infty\(\S^n\)}\to\infty\quad\text{as }k\to\infty.
\end{equation}
The following result provides a first description of the blowing-up behavior of $\(u_k\)_k$:

\begin{lemma}\label{Lem1}
Let $\(g_k\)_k$ and $\(u_k\)_k$ be as defined above. Then, up to a subsequence and a change of sign, there exist $x_1,x_2\in\S^n$ and sequences $\(x_{1,k}\)_k$ and $\(x_{2,k}\)_k$ in $\S^n$ and $\(\mu_{1,k}\)_k$ and $\(\mu_{2,k}\)_k$ in $\(0,2\)$ such that 
\begin{enumerate}
\item[(i)]$\mu_{2,k}\le\mu_{1,k}$ for all $k\in\N$.\smallskip
\item[(ii)]For each $i\in\left\{1,2\right\}$, $x_{i,k}\to x_i$ and $\mu_{i,k}\to0$ as $k\to\infty$.\smallskip
\item[(iii)]$\displaystyle\frac{\mu_{1,k}}{\mu_{2,k}}+\frac{d_k^2}{\mu_{1,k}\mu_{2,k}}\to\infty$ as $k\to\infty$, where $d_k:=\dist_{g_0}\(x_{1,k},x_{2,k}\)$.\smallskip
\item[(iv)] For each $i\in\left\{1,2\right\}$ and $k\in\N$, define 
$$B_{i,k}:=\(\frac{2\sqrt{n\(n-2\)}\mu_{i,k}}{2\mu_{i,k}^2+\big(4-\mu_{i,k}^2\big)\big(1-\cos\(\dist_{g_0}\(\cdot,x_{i,k}\)\)\big)}\)^{\frac{n-2}{2}},$$
where $\dist_{g_0}$ is the distance function on $\(\S^n,g_0\)$. Then
\begin{equation}\label{Lem1Eq1}
\left\| \frac{u_k- \big(B_{1,k}-B_{2,k}\big)}{B_{1,k}+B_{2,k}} \right\|_{L^\infty\(\S^n\)} \to 0 \quad \text{as }k\to\infty.
\end{equation}
\item[(v)] For each $k\in\N$,
\begin{equation}\label{Lem1Eq2}
\displaystyle u_k\(x_{2,k}\)=\min_M u_k=-B_{2,k}\(x_{2,k}\)=-\(\frac{\sqrt{n\(n-2\)}}{\mu_{2,k}}\)^{\frac{n-2}{2}}.
\end{equation}
\item[(vi)] If, moreover,
\begin{equation}\label{Lem1Eq3}
\sqrt{\mu_{1,k}\mu_{2,k}}=\smallo\(d_k\)\quad\text{as }k\to\infty,
\end{equation}
then for each $k\in\N$,
\begin{equation}\label{Lem1Eq4}
u_k\(x_{1,k}\)=\max_M u_k=B_{1,k}\(x_{1,k}\)=\(\frac{\sqrt{n\(n-2\)}}{\mu_{1,k}}\)^{\frac{n-2}{2}}.
\end{equation}
\end{enumerate}
\end{lemma}

In what follows, for simplicity, we rewrite \eqref{Lem1Eq1} as 
\begin{equation}\label{Lem1Eq5}
u_k=B_{1,k}-B_{2,k}+\smallo\(B_{1,k}+B_{2,k}\)\quad\text{in }C^0\(\S^n\)\text{ as }k\to\infty.
\end{equation}

\smallskip
We point out that the assumption $n \le 10$ comes into play in this lemma. Indeed, as is explained in the proof below, it is crucial in order to obtain \eqref{Lem1Eq1}.

\proof[Proof of Lemma~\ref{Lem1}]
Since $\(u_k\)_{k}$ blows-up with finite energy as $k \to \infty$ by \eqref{Sec3Eq6} and \eqref{Sec3Eq7}, a celebrated result of Struwe~\cite{St} (see also the book of Druet, Hebey and Robert~\cite{DHR} and the article of Mazumdar~\cite{Maz} for versions in the Riemannian setting) shows that, up to a subsequence, there exist $m\in\left\{1,2\right\}$, $x_1,\dotsc,x_m\in\S^n$ and sequences $\(\widetilde{x}_{1,k}\)_k,\dotsc,\(\widetilde{x}_{m,k}\)_k$ in $\S^n$ and $\(\widetilde\mu_{1,k}\)_k,\dotsc,\(\widetilde\mu_{m,k}\)_k$ in $\(0,\infty\)$ such that 
\begin{equation}\label{Lem1Eq6}
\widetilde{x}_{i,k}\to x_i\in\S^n\quad\text{and}\quad\widetilde\mu_{i,k}\to0\quad\text{as }k\to\infty\quad\forall i\in\left\{1,\dotsc,m\right\}
\end{equation}
and
\begin{equation}\label{Lem1Eq7}
u_k=B_0+\sum_{i=1}^m\pm\widetilde{B}_{i,k}+\smallo\(1\)\quad\text{in }H^1\(\S^n\)\text{ as }k\to\infty,
\end{equation}
where $H^1\(\S^n\)=H^1\(\S^n,g_0\)$, $B_0$ is a constant-sign solution of \eqref{Sec3Eq1}, which may be equal to $0$, and $\widetilde{B}_{i,k}$ is given by 
$$\widetilde{B}_{i,k}:=\(\frac{2\sqrt{n\(n-2\)}\widetilde\mu_{i,k}}{2\widetilde\mu_{i,k}^2+\big(4-\widetilde\mu_{i,k}^2\big)\big(1-\cos\(\dist_{g_0}\(\cdot,\widetilde{x}_{i,k}\)\)\big)}\)^{\frac{n-2}{2}}.$$
Moreover, in the case where $m=2$, up to a subsequence, we may further assume that 
\begin{equation}\label{Lem1Eq8}
\widetilde\mu_{2,k}\le\widetilde\mu_{1,k}\quad\text{and}\quad\frac{\widetilde\mu_{1,k}}{\widetilde\mu_{2,k}}+\frac{\dist_{g_0}\(\widetilde{x}_{1,k},\widetilde{x}_{2,k}\)^2}{\widetilde\mu_{1,k}\widetilde\mu_{2,k}}\to\infty\quad\text{as }k\to\infty
\end{equation}
(see the remark at the end of Section~3.2 in \cite{DHR}). A consequence of \eqref{Lem1Eq7} is that 
\begin{equation}\label{Lem1Eq9}
\left\|u_k\right\|_{H^1\(\S^n\)}=\left\|B_0\right\|_{H^1\(\S^n\)}+m\Lambda_1\(\S^n,\[g_0\]\)^{\frac{n}{2}}+\smallo\(1\)\quad\text{as }k\to\infty.
\end{equation}
It follows from \eqref{Sec3Eq6} and \eqref{Lem1Eq9} that either [$m=2$ and $B_0 = 0$] or [$m=1$ and $B_0$ is a non-zero constant-sign solution of \eqref{Sec3Eq1}].

\smallskip
Assume first that $m=1$ and $B_0$ is a non-zero constant-sign solution of \eqref{Sec3Eq1}. Up to a change of sign, me may assume that $B_0$ is positive. Since the functions $\(u_k\)_{k}$ change sign, it then follows from \eqref{Lem1Eq7} that
\begin{equation} \label{Lem1Eq10}
u_k = B_0 - \widetilde{B}_{1,k} + \smallo(1) \quad\text{in }H^1\(\S^n\)\text{ as }k\to\infty.
\end{equation}
By using the pointwise blow-up theory for sign-changing solutions developed by Premoselli~\cite{Pr} (see also \cites{DHR,H} in the case of positive solutions), it follows from \eqref{Lem1Eq10} that
\begin{equation} \label{Lem1Eq11}
u_k = B_0 - \widetilde{B}_{1,k} + \smallo \big(\widetilde{B}_{1,k}\big) + \smallo \(1\) \quad\text{in }C^0\(\S^n\)\text{ as }k\to\infty.
\end{equation}
The proof of \cite{Pr} is written for a fixed metric $g$ but adapts straightforwardly to the case of a strongly converging sequence of metrics $\(g_k\)_{k \in \N}$ as is the case here. By using \eqref{Lem1Eq11}, since $n \le 10$ and the Weyl tensor of $\(\S^n,g_0\)$ vanishes everywhere, Theorem~1.2 of Premoselli and V\'etois~\cite{PV3} yields a contradiction with \eqref{Lem1Eq11}. The proof of \cite{PV3} is again stated for a fixed metric but its arguments adapt straightforwardly since they only rely on \eqref{Lem1Eq11} (see \cite{PV3}*{Section~5} for more details).

\smallskip
We have thus proven that $m=2$ and $B_0=0$ hold in \eqref{Lem1Eq7}. Up to a change of sign, since the functions $\(u_k\)_{k}$ change sign, we then obtain
\begin{equation} \label{Lem1Eq12}
u_k=\widetilde{B}_{1,k} - \widetilde{B}_{2,k}+\smallo \(1\) \quad\text{in }H^1\(\S^n\)\text{ as }k\to\infty.
\end{equation}
By using again the pointwise blow-up theory of \cite{Pr}, it follows from \eqref{Lem1Eq12} that
\begin{equation}\label{Lem1Eq13}
u_k=\widetilde{B}_{1,k}-\widetilde{B}_{2,k}+\smallo\big(\widetilde{B}_{1,k}+\widetilde{B}_{2,k}\big)\quad\text{in }C^0\(\S^n\)\text{ as }k\to\infty.
\end{equation}
By putting together \eqref{Lem1Eq6}, \eqref{Lem1Eq8} and \eqref{Lem1Eq13}, we obtain (i) to (iv) in Lemma~\ref{Lem1}.

\smallskip
We now prove that the centers and weights of $\(B_{1,k}\)_k$ and $\(B_{2,k}\)_k$ can be chosen so that (v) and (vi) are also satisfied. For each $k\in\N$, we let $x_{1,k}$, $x_{2,k}$, $\mu_{1,k}$ and $\mu_{2,k}$ be such that \eqref{Lem1Eq2} and \eqref{Lem1Eq4} hold true. For each $i\in\left\{1,2\right\}$, by using \eqref{Lem1Eq2}, \eqref{Lem1Eq4} and \eqref{Lem1Eq13}, we obtain 
\begin{align}
\mu_{i,k}^{\frac{2-n}{2}}&\le\(\frac{2\widetilde\mu_{i,k}\(1+\smallo\(1\)\)}{2\widetilde\mu_{i,k}^2+\big(4-\widetilde\mu_{i,k}^2\big)\big(1-\cos\(\dist_{g_0}\(x_{i,k},\widetilde{x}_{i,k}\)\)\big)}\)^{\frac{n-2}{2}}\nonumber\\
&\le\widetilde\mu_{i,k}^{\frac{2-n}{2}}\(1+\smallo\(1\)\)\quad\text{as }k\to\infty\label{Lem1Eq14}
\end{align}
and
\begin{align}
\mu_{i,k}^{\frac{2-n}{2}}&\ge\widetilde\mu_{i,k}^{\frac{2-n}{2}}\(1+\smallo\(1\)\)\nonumber\\
&\quad-\(\frac{2\widetilde\mu_{3-i,k}\(1+\smallo\(1\)\)}{2\widetilde\mu_{3-i,k}^2+\big(4-\widetilde\mu_{3-i,k}^2\big)\big(1-\cos\(\dist_{g_0}\(\widetilde{x}_{1,k},\widetilde{x}_{2,k}\)\)\big)}\)^{\frac{n-2}{2}}\quad\text{as }k\to\infty.\label{Lem1Eq15}
\end{align}
We now assume that
\begin{equation}\label{Lem1Eq16}
\sqrt{\widetilde\mu_{1,k}\widetilde\mu_{2,k}}=\smallo\(\dist_{g_0}\(\widetilde{x}_{1,k},\widetilde{x}_{2,k}\)\)\quad\text{as }k\to\infty.
\end{equation}
It follows from \eqref{Lem1Eq8} and \eqref{Lem1Eq16} that
$$\(\frac{2\widetilde\mu_{3-i,k}}{2\widetilde\mu_{3-i,k}^2+\big(4-\widetilde\mu_{3-i,k}^2\big)\big(1-\cos\(\dist_{g_0}\(\widetilde{x}_{1,k},\widetilde{x}_{2,k}\)\)\big)}\)^{\frac{n-2}{2}}=\smallo\(\widetilde{\mu}_{i,k}^{\frac{2-n}{2}}\)\quad\text{as }k\to\infty,$$
which together with \eqref{Lem1Eq14} and \eqref{Lem1Eq15} give
\begin{equation}\label{Lem1Eq17}
\mu_{i,k}\sim\widetilde\mu_{i,k}\quad\text{and}\quad\quad\dist_{g_0}\(x_{i,k},\widetilde{x}_{i,k}\)=\smallo\(\mu_{i,k}\)\quad\text{as }k\to\infty.
\end{equation}
By passing to a subsequence and exchanging $\(B_{1,k}\)_k$ and $\(B_{2,k}\)_k$ if necessary and using \eqref{Lem1Eq6}, \eqref{Lem1Eq8}, \eqref{Lem1Eq13} and \eqref{Lem1Eq17}, we now obtain that the sequences $\(x_{1,k}\)_k$, $\(x_{2,k}\)_k$, $\(\mu_{1,k}\)_k$ and $\(\mu_{2,k}\)_k$ simultaneously satisfy (i) to (vi) in Lemma~\ref{Lem1}.
\endproof

Lemma~\ref{Lem1} shows that the singular metric $\left|u_k\right|^{\frac{4}{n-2}}g_k$ decomposes asymptotically as the disjoint union of two round spheres centered at $x_{1,k}$ and $x_{2,k}$, respectively. The location of these two points is unknown. Lemma~\ref{Lem1} does not claim, in particular, that $d_k = \dist_{g_0}\(x_{1,k},x_{2,k}\)$ has a positive limit as $k \to \infty$. In order to prove Theorems~\ref{Th1} and \ref{Th2}, we need a more precise description of the blow-up behavior of $u_k$. As is often the case with the Yamabe equation, it is convenient to work with the conformal normal coordinate system introduced by Lee and Parker~\cite{LP}. We now explain how we do this in our setting. We first let $\varepsilon_0\in\(0,\infty\)$ be fixed and $\varphi_0$ be a smooth positive function on $\S^n\times\S^n$ such that 
\begin{equation}\label{Lem2Eq1}
\varphi_0\(x,y\):=\(\frac{2}{1+\cos\(\dist_{g_0}\(x,y\)\)}\)^{\frac{n-2}{2}}\quad\forall x\in\S^n,\,y\in \B_{g_0}\(x,r_0\),
\end{equation}
where
$$r_0:=2\tan^{-1}\(\varepsilon_0/2\).$$
For any $x \in \S^n$, we define 
$$g_{0,x}:=\varphi_0\(x,\cdot\)^{2^*-2}g_0\quad\text{and}\quad\hat{g}_{0,x}:={\exp_{0,x}}^*g_{0,x},$$ where $\exp_{0,x}$ is the exponential map at $x$ with respect to $g_{0,x}$ and where we identify $T_x \S^n$ with $\R^n$. A simple reformulation of the stereographic projection gives
\begin{equation}\label{Lem2Eq4}
\hat{g}_{0,x}=\xi\quad\text{in }\B_\xi\(0,\varepsilon_0\),
\end{equation}
where $\xi$ is the Euclidean metric on $\R^n$. The following result explains how a similar system of conformal normal coordinates can be obtained for the metrics $\(g_k\)_k$:

\begin{lemma}\label{Lem2}
Let $\(g_k\)_k$ be as in Lemma~\ref{Lem1}. Then there exists a sequence $\(\varphi_k\)_k$ of smooth positive functions on $\S^n\times\S^n$ such that the following hold:
\begin{enumerate}
\item[(i)]$\varphi_k\to\varphi_0$ in $C^m\(\S^n\times\S^n\)$ as $k\to\infty$ for all $m\in\N$.\smallskip
\item[(ii)] For each $k\in\N$ and $x\in\S^n$,
\begin{equation}\label{Lem2Eq2}
\varphi_k\(x,x\)=1\quad\text{and}\quad\nabla\varphi_k\(x,\cdot\)\(x\)=0.
\end{equation}
\item[(iii)] For each $k\in\N$ and $x\in\S^n$, let $g_{k,x}$ and $\hat{g}_{k,x}$ be the metrics on $\S^n$ and $\R^n$, respectively, defined as
$$g_{k,x}:=\varphi_k\(x,\cdot\)^{2^*-2}g_k\quad\text{and}\quad\hat{g}_{k,x}:={\exp_{k,x}}^*g_{k,x},$$
where $\exp_{k,x}$ is the exponential map at $x$ with respect to $g_{k,x}$. Then 
\begin{equation}\label{Lem2Eq3}
\dv_{\hat{g}_{k,x}}\(y\)=\big(1+\smallo\big(\left|y\right|^N\big)\big)\dy \quad\text{as }k\to\infty
\end{equation}
uniformly with respect to $x\in\S^n$ and $y\in\B_\xi\(0,\varepsilon_0\)$, where $\dv_{\hat{g}_{k,x}}$ and $\dy$ are the volume elements of $\(\R^n,\hat{g}_{k,x}\)$ and $\(\R^n,\xi\)$, respectively, and $N\in\N$ can be chosen arbitrarily large.
\end{enumerate}
\end{lemma}

\proof[Proof of Lemma~\ref{Lem2}]
We fix $N\in\N$ large enough ($N \ge n$ is sufficient in the rest of the paper). For any $x \in \S^n$ and $k \in \N$, we define $g_{k,x}^{\(1\)} = \varphi_0\(x, \cdot\)^{2^*-2} g_k$. Since $g_k\to g_0$ in $C^m\(\S^n\)$ for all $m\in\N$ as $k \to\infty$, we obtain $i_{g_k}\(\S^n\) \ge \pi + \smallo\(1\)$ as $k \to \infty$, where $i_{g_k}\(\S^n\)$ is the injectivity radius of $\(\S^n, g_k\)$, and  $g_{k,x}^{\(1\)}\to g_{0,x}$ in $C^m\(\S^n\)$ for all $m\in\N$ as $k \to \infty$. In particular, $g_{k,x}^{\(1\)}$ satisfies $ \det g_{k,x}^{\(1\)}\(y\) = 1 + \smallo \big( \left|y\right|^{2} \big)$ for $y \in \B_{\xi}\(0, 2 \varepsilon_0\)$ in its own exponential coordinates at $x$. We now follow the strategy of \cite{LP}*{Theorem~5.1} to inductively correct $g_{k,x}^{\(1\)}$ by a conformal factor in order to obtain \eqref{Lem2Eq3}. First, by \cite{LP}*{Theorem~5.2}, there exists a unique homogeneous polynomial of degree $2$ in exponential coordinates for $g_{k,x}^{\(1\)}$ at $x$, that we denote  $f_{k,x}^{\(2\)}$,  such that $g_{k,x}^{\(2\)}:= e^{ 2 f_{k,x}^{\(2\)}} g_{k,x}^{\(1\)}$ satisfies $ \det g_{k,x}^{\(2\)}\(y\) = 1 + \smallo \big( \left|y\right|^3 \big)$ for $y \in \B_{\xi}\(0, 2 \varepsilon_0\)$ in its own exponential coordinates at $x$ (that the error term is $\smallo\big(\left|y\right|^3\big)$ and not $\bigO\big(\left|y\right|^3\big)$ is again a consequence of the convergence of $g_k$ to $g_0$). The construction of $f_{k,x}^{\(2\)}$ is purely algebraic and its coefficients solely depend on covariant derivatives of order at most $2$ of $g_{k,x}^{\(1\)}$ (see for example the arguments in the proof of \cite{LP}*{~Theorem 5.2} and \cite{LP}*{Lemma 5.3}). As a consequence, $f_{k,x}^{\(2\)} \to 0$ in $C^m_{loc}\(\R^n\)$ for all $m\in\N$ as $k \to \infty$. Repeating the process, we construct by induction, for any $ \ell\in\N$, a homogeneous polynomial $f_{k,x}^{\(\ell+1\)}$ of degree $\ell+1$ in exponential coordinates for $g_{k,x}^{\(\ell\)}$ at $x$ such that $g_{k,x}^{\(\ell+1\)}:= e^{ 2 f_{k,x}^{\(\ell+1\)}} g_{k,x}^{\(\ell\)}$ satisfies $ \det g_{k,x}^{\(\ell+1\)}\(y\)= 1 + \smallo \big( \left|y\right|^{\ell+2} \big)$ for $y \in \B_{\xi}\(0, 2 \varepsilon_0\)$ in its own exponential coordinates at $x$ and such that $f_{k,x}^{\(\ell+1\)} \to 0$ in $C^m_{loc}\(\R^n\)$ for all $m \in\N$ as $k \to \infty$. We stop at $\ell= N-1$ and define 
$$ \varphi_{k,x}^{2^*-2}:= \varphi_0\(x,\cdot\)^{2^*-2}e^{2f_{1,k} + \dotsc+2f_{k,x}^{\(N-1\)}} \quad \text{ in } \B_{g_0}\(x, r_0/2\). $$
By construction, we obtain that $\varphi_{k,x}$ is smooth in $\B_{g_0}\(x,r_0/2\)$ and converges to $\varphi_0\(x, \cdot\)$ in $C^m\(\B_{g_0}\(x, r_0/2\)\)$ for all $m \in\N$ as $k \to \infty$. We then extend $\varphi_{k,x}$ to a smooth function on the whole of $\S^n$ by interpolating it to $\varphi_0\(x, \cdot\)$ in $\B_{g_0}\(x, r_0\) \setminus \B_{g_0}\(x, r_0/2\)$ so that the convergence to $\varphi_0\(x, \cdot\)$ holds true in $C^m\(\S^n\)$. Remark that by construction, we also obtain that $\varphi_{k,x}$ satisfies \eqref{Lem2Eq2} and \eqref{Lem2Eq3}. It remains to define $\varphi_k\(x,y\) = \varphi_{k,x}\(y\)$ for all $\(x,y\) \in \S^n \times \S^n$. The smoothness of $\varphi_k$ in the two variables and its convergence to $\varphi_0$ in $C^m\(\S^n\times\S^n\)$ as $k\to\infty$ for all $m\in\N$ follow from the smoothness of $\varphi_0$ in $\S^n \times \S^n$ and of the exponential map of $g_k$. This ends the proof of Lemma \ref{Lem2}.
\endproof

As observed by Khuri, Marques and Schoen~\cite{KMS}, it is convenient to express a given metric in exponential form when working with conformal normal coordinates. The following lemma explains how this applies to the metrics $g_k$ we investigate here:

\begin{lemma}\label{Lem3}
Let $\(\hat{g}_{k,x}\)_{k,x}$ and $\varepsilon_0$ be as in \eqref{Lem2Eq1}. Then, provided $k$ is large enough, the following hold
\begin{itemize}
\item[(i)] For each $k\in\N$ and $x\in\S^n$, there exists a smooth symmetric 2-covariant tensor $h_{k,x}$ in $\R^n$ such that
\begin{align}
&\hat{g}_{k,x}=\exp\(h_{k,x}\),\label{Lem3Eq1}\allowdisplaybreaks\\
&h_{k,x}\(y\)y=0\quad\forall y\in\R^n\label{Lem3Eq2}
\end{align}
and
\begin{equation}\label{Lem3Eq3}
\tr\(h_{k,x}\(y\)\)=\smallo\big(\left|y\right|^N\big)\quad\text{as }k\to\infty
\end{equation}
uniformly with respect to $x\in\S^n$ and $y\in \B_\xi\(0,\varepsilon_0\)$, where $\exp$ and $\tr$ are the matrix exponential and trace maps, respectively.\smallskip
\item[(ii)] The tensor $h_{k,x}$ satisfies
\begin{equation}\label{Lem3Eq4} 
h_{k,x}\(y\)=H_{k,x}\(y\)+\smallo\big(\left|y\right|^{\max\(n-3,2\)}\big)\quad\text{as }k\to\infty
\end{equation}
uniformly with respect to $x\in\S^n$ and $y\in \B_\xi\(0,\varepsilon_0\)$, where $H_{k,x}\(y\)= 0$ if $3 \le n \le 5$ and $H_{k,x}\(y\)$ is of the form
$$H_{k,x}\(y\)=\sum_{\left|\alpha\right|=2}^{n-4}h_{k,x,\alpha}y^\alpha$$
for some trace-free symmetric real matrices $h_{k,x,\alpha}$ which do not depend on $y$ if $n \ge 6$. Moreover, 
\begin{equation}
H_{k,x}\(y\)y=0\quad\text{and}\quad\tr\(H_{k,x}\(y\)\)=0\quad\forall y\in\R^n\label{Lem3Eq5}
\end{equation}
and \eqref{Lem3Eq4} can be differentiated.\smallskip
\item[(iii)] The scalar curvature of $\hat{g}_{k,x}$ satisfies
\begin{align}
\Scal_{\hat{g}_{k,x}}\(y\)&=\sum_{a,b=1}^n\partial_{y_a}\partial_{y_b}\(h_{k,x}\)_{ab}\(y\)-\sum_{a,b,c=1}^n\bigg(\partial_{y_b}\(\(H_{k,x}\)_{ab}\partial_{y_c}\(H_{k,x}\)_{ac}\)\nonumber\\
&\quad-\frac{1}{2}\partial_{y_b}\(H_{k,x}\)_{ab}\partial_{y_c}\(H_{k,x}\)_{ac}+\frac{1}{4}\(\partial_{y_c}\(H_{k,x}\)_{ab}\)^2\bigg)\(y\)\nonumber\\
&\quad+\bigO\(\sum_{\left|\alpha\right|=2}^{d_n}\left|h_{k,x,\alpha}\right|^2\left|y\right|^{2\left|\alpha\right|}\)+\smallo\big(\left|y\right|^{n-1}\big)\quad\text{as }k\to\infty\label{Lem3Eq6}
\end{align}
and
\begin{align}
\Scal_{\hat{g}_{k,x}}\(y\)&=\sum_{a,b=1}^n\partial_{y_a}\partial_{y_b}\(h_{k,x}\)_{ab}\(y\)+\bigO\(\sum_{\left|\alpha\right|=2}^{d_n}\left|h_{k,x,\alpha}\right|^2\left|y\right|^{2\left|\alpha\right|-2}\)\nonumber\\
&\quad+\smallo\big(\left|y\right|^{\max\(n-3,2\)}\big)\quad\text{as }k\to\infty\label{Lem3Eq7}
\end{align}
uniformly with respect to $x\in\S^n$ and $y\in \B_\xi\(0,\varepsilon_0\)$, where $\(h_{k,x}\)_{ab}$ and $\(H_{k,x}\)_{ab}$ are the coefficients of $h_{k,x}$ and $H_{k,x}$, respectively, and
$$d_n:=\[\frac{n-2}{2}\].$$
\end{itemize}
\end{lemma}

In \eqref{Lem3Eq6} and \eqref{Lem3Eq7}, $\left|h_{k,x,\alpha}\right|$ denotes the matrix norm of the square matrix $h_{k,x,\alpha}$. Remark that \eqref{Lem3Eq2} and the first identity in \eqref{Lem3Eq5} can also be written as
$$\sum_{b=1}^n h_{k,x}\(y\)_{ab} y_b = 0\quad\text{and}\quad\sum_{b=1}^n H_{k,x}\(y\)_{ab} y_b = 0\quad\forall a\in\left\{1,\dotsc,n\right\},\,y\in\R^n.$$
We point out that in the case where $3 \le n \le 5$, \eqref{Lem3Eq4} and \eqref{Lem3Eq7} simply give
$$ h_{k,x}\(y\)=\smallo\big(\left|y\right|^{2}\big)\quad \text{and} \quad \Scal_{\hat{g}_{k,x}}\(y\) = \smallo\big(\left|y\right|^{2}\big) \quad \text{as }k\to\infty$$
uniformly with respect to $x\in\S^n$ and $y\in \B_\xi\(0,\varepsilon_0\)$. When $n \ge 6$, we observe that 
$$\left|h_{k,x,\alpha}\right|\to0\quad\text{as }k \to \infty$$
uniformly with respect to $x\in\S^n$ and $\alpha$ such that $2\le\left|\alpha\right|\le n-4$, which simply follows from \eqref{Lem2Eq4} and Lemma~\ref{Lem2}~(i).

\proof[Proof of Lemma~\ref{Lem3}]
We let $x \in \S^n$ and $k \in \N$. First, since $\hat{g}_{k,x}$ is symmetric and positive definite, we obtain that there exists a unique symmetric 2-covariant tensor $h_{k,x}$ in $\R^n$ such that $\hat{g}_{k,x}=\exp\(h_{k,x}\)$ in $\R^n$. In addition, $h_{k,x}$ smoothly depends on $g_{k,x}$ as a matrix and is thus smooth in $\R^n$. This proves \eqref{Lem3Eq1}. Properties of the matrix exponential show that $\det \hat{g}_{k,x}\(y\) = e^{\tr\(h_{k,x}\(y\)\)}$ for all $y \in \R^n$, hence \eqref{Lem3Eq3} follows from \eqref{Lem2Eq3}. We now prove \eqref{Lem3Eq2}. Since $g_{k} \to g_0$ in $C^m\(\S^n\)$ for all $m \in\N$ as $k \to \infty$, Cartan's expansion shows that there exists a sequence $\(\eta_k\)_k$ of positive numbers such that $\eta_k \to 0$ as $k \to \infty$ and 
$$\left|h_{k,x}\(y\)\right| \le \eta_k \left|y\right|^2 \quad  \forall y\in \B_\xi\(0,\varepsilon_0\),$$ 
which, together with \eqref{Lem3Eq1}, implies that
$$ \left|\hat{g}_{k,x}\(y\)y -y - h_{k,x}\(y\)y \right| \le \eta_k\varepsilon_0^2e^{\eta_k\varepsilon_0^2}\left|h_{k,x}\(y\) y\right|  \quad  \forall y\in \B_\xi\(0,\varepsilon_0\).$$
Independently, Gauss's lemma states that
$$\hat{g}_{k,x}\(y\) y = y,\quad\text{i.e.}\quad\sum_{b=1}^n \hat{g}_{k,x}\(y\)_{ab} y_b = y_a\quad\forall a \in \left\{1, \dotsc, n\right\},\,y\in \B_\xi\(0,\varepsilon_0\).$$
This together with the previous inequality then shows that 
$$ \left| h_{k,x}\(y\)y \right| \le \eta_k\varepsilon_0^2e^{\eta_k\varepsilon_0^2} \left|h_{k,x}\(y\) y\right| \quad  \forall y\in \B_\xi\(0,\varepsilon_0\),$$
which implies that $h_{k,x}\(y\)y = 0$ provided $k$ is large enough. This proves \eqref{Lem3Eq2}. We remark that this argument also shows that \eqref{Lem3Eq2} holds true without passing to a subsequence in $k$ provided $\varepsilon_0$ is chosen small enough.

\smallskip
We now prove \eqref{Lem3Eq4} and \eqref{Lem3Eq5}. We begin with considering the case where $3 \le n \le 5$. It follows from \eqref{Lem2Eq4} and Lemma~\ref{Lem2}~(i) that $\hat{g}_{k,x}\to\xi$ in $C^m\(\B_\xi\(0,\varepsilon_0\)\)$ for all $m \in\N$ as $k\to\infty$ uniformly with respect to $x\in\S^n$. Cartan's expansion of $\hat{g}_{k,x}$ then shows that $\hat{g}_{k,x}\(y\) = \xi + \smallo\big(\left|y\right|^2\big)$ as $k \to \infty$ uniformly with respect to $x\in\S^n$ and $y\in\B_\xi\(0,\varepsilon_0\)$. Letting $H_{k,x} = 0$ for all $x \in \S^n$, this proves \eqref{Lem3Eq4} and \eqref{Lem3Eq5} since $\max\(n-3,2\) = 2$ in this case. We now assume that $n \ge 6$. In this case, we let $H_{k,x}$ be the Taylor polynomial of $h_{k,x}$ at $0$ of degree $n-4$. The fact that $H_{k,x}$ is of degree at least $2$ again follows from Cartan's expansion of $\hat{g}_{k,x}$ at $0$. Remarking that $\max\(n-3,2\) =n-3$ in this case, Taylor's formula now shows, since $h_{k,x}\(y\)\to0$ as $k \to \infty$  uniformly with respect to $x\in\S^n$ and $y\in\B_\xi\(0,\varepsilon_0\)$, that \eqref{Lem3Eq4} holds. Since $\tr\(H_{k,x}\(y\)\)$ is a polynomial in $y$ of degree at most $n-4$, by using \eqref{Lem3Eq3}, we obtain that the matrices $h_{k,x,\alpha}$ are trace-free, which gives $\tr\(H_{k,x}\(y\)\)=0$ for all $y\in\R^n$. Similarly, a straightforward expansion of \eqref{Lem3Eq2} in powers of $y$ using \eqref{Lem3Eq4} shows that 
$$\sum_{\left|\alpha\right|=d} y^{\alpha}h_{k,x,\alpha}\(y\)=0\quad\forall d\in\left\{2,\dotsc,n-4\right\},\,y \in \R^n,$$
which completes the proof of \eqref{Lem3Eq5}.

\smallskip
We finally prove \eqref{Lem3Eq6} and \eqref{Lem3Eq7}. For this, we use the precise expansion of scalar curvature for metrics of the form $g = \exp\(h\)$ which was proven in \cite{B}*{Proposition~26}. When applied to our metric $\hat{g}_{k,x} = \exp\(h_{k,x}\)$, this formula gives
\begin{align*}
\Scal_{\hat{g}_{k,x}}\(y\)&=\sum_{a,b=1}^n\partial_{y_a}\partial_{y_b}\(h_{k,x}\)_{ab}\(y\)-\sum_{a,b,c=1}^n\bigg(\partial_{y_b}\(\(h_{k,x}\)_{ab}\partial_{y_c}\(h_{k,x}\)_{ac}\)\\
&\quad-\frac{1}{2}\partial_{y_b}\(h_{k,x}\)_{ab}\partial_{y_c}\(h_{k,x}\)_{ac}+\frac{1}{4}\(\partial_{y_c}\(h_{k,x}\)_{ab}\)^2\bigg)\(y\)\\
&\quad+\bigO\( \left| h_{k,x}\(y\)\right|^2 \left| \nabla_\xi^2 h_{k,x}\(y\)\right| \) +\bigO\( \left| h_{k,x}\(y\)\right| \left| \nabla_\xi h_{k,x}\(y\)\right|^2 \), 
\end{align*}
where $\nabla_{\xi}$ denotes the covariant derivative with respect to the Euclidean metric. By using \eqref{Lem3Eq4}, we obtain, for any $a,b,c,d \in \left\{1, \dotsc, n \right\}$, 
$$\partial_{y_c}\(h_{k,x}\)_{ab}\(y\) =  \partial_{y_c}\(H_{k,x}\)_{ab}\(y\) + \smallo \big( \left|y\right|^{\max\(n-4, 1\)}\big)$$
and
$$\partial_{y_d} \partial_{y_c}\(h_{k,x}\)_{ab}\(y\) =   \partial_{y_d} \partial_{y_c}\(H_{k,x}\)_{ab}\(y\) + \smallo \big( \left|y\right|^{\max\(n-5, 0\)}\big).$$
We then obtain \eqref{Lem3Eq6} by plugging the latter formulae into the previous expansion of $\Scal_{\hat{g}_{k,x}}$. Finally, \eqref{Lem3Eq7} simply follows from \eqref{Lem3Eq6} by again using  \eqref{Lem3Eq4}.
\endproof

For each $k\in\N$, $x\in\S^n$ and $d\in\left\{2,\dotsc,n-4\right\}$, we define 
\begin{equation}\label{Lem3Eq8}
H_{k,x,d}\(y\) := \sum_{\left|\alpha\right| = d} h_{k,x,\alpha} y^{\alpha}\quad\forall y\in\R^n. 
\end{equation}
It is easy to see that $H_{k,x,d}$ is a homogeneous polynomial of degree $d$ and 
$$H_{k,x} = \sum_{d=2}^{n-4} H_{k,x,d}.$$
As a consequence of Lemma~\ref{Lem3}, we obtain
\begin{equation}\label{Lem3Eq9} 
\sum_{b=1}^n \(H_{k,x,d}\(y\)\)_{bb} = 0\quad\text{and}\quad \sum_{b=1}^n \(H_{k,x,d}\(y\)\)_{ab} y_b  = 0\quad\forall a\in\left\{1,\dotsc,n\right\},\,y\in\R^n. 
\end{equation}

\smallskip
For each $d \in \left\{0, \dotsc, n-4\right\}$, we define
\begin{equation} \label{Lem4Eq1}
w_{k,x,d}\(y\):= \sum_{\left|\alpha\right| = d} \sum_{a,b=1}^n  \(h_{k,x,\alpha}\)_{ab}\partial_{y_a}\partial_{y_b}\(y^\alpha\)\quad\forall y\in\R^n.
\end{equation}
It is easy to see that $w_{k,x,d}$ is a homogeneous polynomial of degree $d-2 \in \left\{0 , \dotsc, n-6 \right\}$, which by \eqref{Lem3Eq7} captures the main term in the Taylor expansion of $\Scal_{\hat{g}_{k,x}}\(y\)$ at 0. We claim that, for each $k\in\N$ and $x\in\S^n$ and for each dimension $n \ge 3$,
\begin{equation} \label{Lem4Eq2}
w_{k,x,d}\(y\)= 0 \quad \forall y\in\R^n,\, d\in\left\{0,1,2,3\right\}.
\end{equation} 
This is obvious in the case where $d\in\left\{0,1\right\}$. When $d\in\left\{2,3\right\}$, $w_{k,x,d}$ is a homogeneous polynomial of order $0$ or $1$, respectively. By using \eqref{Lem3Eq7} and remarking that 
$$\Scal_{\hat{g}_{k,x}}\(y\) = \bigO \big(\left|y\right|^2\big)$$
uniformly with respect to $y\in\B_{g_0}\(x,r_0\)$ and $k\in\N$, which follows from properties of the conformal normal coordinates (see \cite{LP}), we obtain \eqref{Lem4Eq2}. As a consequence of \eqref{Lem4Eq2}, we obtain that if $3\le n\le 7$, then
\begin{equation} \label{Lem4Eq3}
w_{k,x,d}\(y\)= 0 \quad \forall y\in\R^n,\, d\in\left\{0,\dotsc,n-4\right\},
\end{equation}
hence \eqref{Lem4Eq1} is trivial in this case. The error terms $w_{k,x,d}$ for $4 \le d \le n-4$ appearing in the expansion of $\Scal_{\hat{g}_{k,x}}\(y\)$ at 0 are therefore non-zero only in dimensions $n \ge 8$. The next result is crucial for taking into account these terms in order to refine the asymptotic analysis of Lemma~\ref{Lem1}:

\begin{lemma}\label{Lem4}
Assume that $n \ge 8$. Let $\(h_{k,x,\alpha}\)_{k,x,\alpha}$ be as in Lemma~\ref{Lem3}. For each $k\in\N$, $x\in\S^n$ and $d\in\left\{4,\dotsc,n-4\right\}$, let $w_{k,x,d}$ be defined by \eqref{Lem4Eq1}. Then there exists a unique family of real numbers $\(\gamma_{k,x,d,l,m}\)_{0\le m\le \[\(d-2\)/2\], 0\le l\le m+2}$ such that the function $v_{k,x,d}:\R^n\to\R$ defined by 
$$v_{k,x,d}\(y\):=\big(1+\left|y\right|^2\big)^{-\frac{n}{2}}\sum_{m=0}^{\[\(d-2\)/2\]} \sum_{l=0}^{m+2}\gamma_{k,x,d,l,m}\left|y\right|^{2l}\Delta_\xi^m w_{k,x,d}\(y\)\quad\forall y\in\R^n$$
solves the equation
\begin{equation}\label{Lem4Eq4}
\Delta_\xi v_{k,x,d}=\(2^*-1\)U_0^{2^*-2}v_{k,x,d}+U_0w_{k,x,d}\quad\text{in }\R^n,
\end{equation}
where
$$U_0\(y\):=\(\frac{\sqrt{n\(n-2\)}}{1+\left|y\right|^2}\)^{\frac{n-2}{2}}\quad\forall y\in\R^n.$$
Moreover, 
\begin{equation}\label{Lem4Eq5}
v_{k,x,d}\(0\)=\left|\nabla v_{k,x,d}\(0\)\right|=0
\end{equation}
and, for each $j\in\N$,
\begin{equation}\label{Lem4Eq6}
\left|\nabla^j v_{k,x,d}\(y\)\right| = \bigO \( \sum_{\left|\alpha\right|=d}\frac{\left|h_{k,x,\alpha}\right|}{\(1+\left|y\right|\)^{n-d-2+j}}\)
\end{equation}
uniformly with respect to $x \in \S^n$, $y \in \R^n$ and $k\in\N$.
\end{lemma}

\proof[Proof of Lemma~\ref{Lem4}]
We follow the arguments of Khuri, Marques and Schoen~\cite{KMS} (see also Li and Zhang~\cites{LZha1,LZha2}). We let $x \in \S^n$, $k \in \N$ and $d\in\left\{4,\dotsc,n-4\right\}$. First, it is easy to see that
\begin{equation} \label{Lem4Eq7}
w_{k,x,d}\(y\)= \sum_{a,b=1}^n \partial_{y_a} \partial_{y_b} \(H_{k,x,d}\)_{ab}\(y\)= \divergence_{\xi} \divergence_{\xi} H_{k,x,d}\(y\)\quad\forall y\in\R^n,
\end{equation}
so that, in particular, $w_{k,x,d}$ is a homogeneous polynomial of degree $d-2 \in \left\{2 , \dotsc, n-6 \right\}$. We recall that two homogeneous polynomials $p$ and $q$ in $\R^n$ are said to be orthogonal if 
$$\int_{\S^{n-1}} pq d\sigma= 0,$$ where $\dsigma$ is the volume element of the round metric on the $\(n-1\)$-sphere $\S^{n-1}$. We now verify that $w_{k,x,d}$ is orthogonal to $\Span\left\{1, y_1, \dotsc, y_n\right\}$. Since $w_{k,x,d}$ is homogeneous of degree $d-2$, we obtain
\begin{equation} \label{Lem4Eq8}
\int_{\S^{n-1}}w_{k,x,d}\dsigma = \(n+d-2\) \int_{\BB^n}   w_{k,x,d} \dy,
\end{equation}
where $\BB^n:=\B_\xi\(0,1\)$. On the other hand, by using \eqref{Lem3Eq9} and \eqref{Lem4Eq7} together with an integration by parts, we obtain
\begin{align}\label{Lem4Eq9}
\int_{\BB^n}   w_{k,x,d} \dy&=\sum_{a,b=1}^n \int_{\S^{n-1}}\partial_{y_b}\(H_{k,x,d}\)_{ab}\(y\)y_a\dsigma\(y\)\nonumber\\
&=\sum_{b=1}^n\int_{\S^{n-1}}\(\partial_{y_b}\(\sum_{a=1}^n\(H_{k,x,d}\)_{ab}\(y\)y_a\)-\(H_{k,x,d}\)_{bb}\(y\)\)\dsigma\(y\)\nonumber\\
&=0.
\end{align}
It follows from \eqref{Lem4Eq8} and \eqref{Lem4Eq9} that $w_{k,x,d}$ is orthogonal to 1. In a similar way, for each $i\in\left\{1,\dotsc,n\right\}$, we obtain
\begin{equation} \label{Lem4Eq10}
\int_{\S^{n-1}}w_{k,x,d}\(y\)y_i\dsigma\(y\) = \(n+d-1\) \int_{\BB^n}   w_{k,x,d}\(y\)y_i \dy
\end{equation}
and, by using again \eqref{Lem3Eq9},
\begin{align}\label{Lem4Eq11}
\int_{\BB^n}   w_{k,x,d}\(y\)y_i\dy&=\sum_{a=1}^n \int_{\BB^n}\partial_{y_a}\(\sum_{b=1}^n\partial_{y_b}\(H_{k,x,d}\)_{ab}\(y\)y_i-\(H_{k,x,d}\)_{ai}\(y\)\)\dy\nonumber\\
&=\sum_{a=1}^n \int_{\S^{n-1}}\(\sum_{b=1}^n\partial_{y_b}\(H_{k,x,d}\)_{ab}\(y\)y_i-\(H_{k,x,d}\)_{ai}\(y\)\)y_a\dsigma\(y\)\nonumber\allowdisplaybreaks\\
&=\sum_{b=1}^n\int_{\S^{n-1}}\Bigg(\partial_{y_b}\(\sum_{a=1}^n\(H_{k,x,d}\)_{ab}\(y\)y_a\)y_i-\(H_{k,x,d}\)_{bb}\(y\)y_i\nonumber\\
&\quad-\(H_{k,x,d}\)_{bi}\(y\)y_b\Bigg)\dsigma\(y\)\nonumber\\
&=0.
\end{align}

\smallskip
We now let $q$ be a homogeneous polynomial of degree $\ell\in\N$ in $\R^n$. We denote by $F\(q\)$ the set of all homogeneous polynomials in $\R^n$ which are obtained as linear combinations of homogeneous polynomials of the form $\left|y\right|^{2r} \Delta_{\xi}^s q$ for some nonnegative integers $r$ and $s$ such that $r \le s+2$. A simple observation is that any element in $F\(q\)$ has degree lower than equal to $\ell + 4$. It was shown in \cite{KMS}*{Proposition~4.1} that if $0 \le \ell < n-4$ and $q$ is orthogonal to $\Span\left\{1, y_1, \dotsc, y_n\right\}$, then there exists a unique polynomial $\psi \in F\(q\)$ such that the function $v:\R^n\to\R$ defined by
$$v\(y\):= \big(1+\left|y\right|^2\big)^{-\frac{n}{2}}\psi\(y\)\quad\forall y\in\R^n$$ 
solves the equation
$$\Delta_\xi v=\(2^*-1\)U_0^{2^*-2}v+U_0 q \quad\text{in }\R^n.$$
Since $w_{k,x,d}$ is orthogonal to $\Span\left\{1, y_1, \dotsc, y_n\right\}$ and $\deg\(w_{k,x,d}\)= d-2 \le  n-6 < n-4$, we can apply \cite{KMS}*{Proposition~4.1} with $q= w_{k,x,d}$ and $\ell=d-2$, which gives the existence of a polynomial $\psi_{k,x,d}\in F\(w_{k,x,d}\)$ such that the function $v_{k,x,d}:\R^n\to\R$ defined by
$$v_{k,x,d}\(y\):=  \big(1+\left|y\right|^2\big)^{-\frac{n}{2}}\psi_{k,x,d}\(y\)\quad\forall y\in\R^n.$$ 
solves \eqref{Lem4Eq4}. By definition of $F\(w_{k,w,d}\)$, we then obtain that there exists a family of real numbers $\(\gamma_{k,x,d,l,m}\)_{0\le m\le \[\(d-2\)/2\], 0\le l\le m+2}$ such that 
$$\psi_{k,x,d}\(y\) = \sum_{m=0}^{\[\(d-2\)/2\]} \sum_{l=0}^{m+2}\gamma_{k,x,d,l,m}\left|y\right|^{2l}\Delta_\xi^m w_{k,x,d}\(y\)\quad\forall y\in\R^n.$$ 
Moreover, by using the fact that $w_{k,x,d}$ is orthogonal to $\Span\left\{1, y_1, \dotsc, y_n\right\}$ together with successive integration by parts, we obtain
$$\Delta_{\xi}^{[\frac{d-2}{2}]} w_{k,x,d} = 0\quad\text{in }\R^n,$$
which, together with the fact that $w_{k,x,d}$ is a homogeneous polynomial of degree $d-2$, implies \eqref{Lem4Eq5}. 

\smallskip
It remains to prove \eqref{Lem4Eq6}. We provide here a PDE proof of this fact. We let $G_0$ be the Green's function of the operator $\Delta_{\xi} - \(2^*-1\) U_0^{2^*-2}$ in $\R^n$. By the stereographic projection, $G_0$ is in one-to-one correspondence with the Green's function of the operator $\Delta_{g_0} - n$ in $\mathbb{S}^n$. As a consequence, we obtain that there exists a constant $C$ depending only on $n$ such that 
$$\left|G_0\(y,z\)\right| \le C \left|y-z\right|^{2-n}\quad\forall y \neq z\in\R^n$$
(see for example the arguments in \cite{Pr}*{Lemma 2.2}). We let $\tilde{v}_{k,x,d}:\R^n\to\R$ be the function defined by
$$ \tilde{v}_{k,x,d}\(y\) := \int_{\R^n} G_0\(y,z\) U_0\(z\) w_{k,w,d}\(z\) \dz\quad\forall y\in\R^n.$$
The integral in the definition of $\tilde{v}_{k,x,d}$ is finite since 
$$\left| w_{k,x,d}\(y\) \right| \le C \(\sum_{\left|\alpha\right|=d}\left|h_{k,x,\alpha}\right|\) \left|y\right|^{d-2}\quad\forall y\in\R^n$$
for some constant $C$ depending only on $n$ and $d$, and since $d-2 \le n-6 < n-4$. Giraud-type arguments give
\begin{equation} \label{Lem4Eq12}
\left|\tilde{v}_{k,x,d}\(y\)\right| \le C\(1+\left|y\right|\)^{d+2-n}\sum_{\left|\alpha\right|=d}\left|h_{k,x,\alpha}\right| \quad \forall y \in \R^n
\end{equation}
for some constant $C$ depending only on $n$ and $d$. On the other hand, it is easy to see that $G_0\(0,z\)$ is radially symmetric in $z$. As a consequence of the orthogonality of $w_{k,x,d}$ to $\Span\left\{1, y_1, \dotsc, y_n\right\}$, we then obtain 
\begin{equation} \label{Lem4Eq13}
\tilde{v}_{k,x,d}\(0\)=\left|\nabla \tilde{v}_{k,x,d}\(0\)\right|=0.
\end{equation}
Moreover, by definition, $\tilde{v}_{k,x,d}$ solves \eqref{Lem4Eq4}. As a consequence, we obtain 
\begin{equation} \label{Lem4Eq14}
\Delta_\xi\(v_{k,x,d}- \tilde{v}_{k,x,d}\) - \(2^*-1\) U_0^{2^*-2}\(v_{k,x,d}- \tilde{v}_{k,x,d}\) = 0\quad\text{in }\R^n.
\end{equation}
By observing that $\deg\(\psi_{k,x,d}\)\le d+2$, it follows from the definition of $v_{k,x,d}$ that
$$\sup_{y \in \R^n} \(\(1+\left|y\right|\)^{n-d-2}\left|v_{k,x,d}\(y\)\right|\) < \infty.$$ 
Therefore, since $n-d-2 \ge 2$, we obtain
\begin{equation} \label{Lem4Eq15}
\left|\(v_{k,x,d}- \tilde{v}_{k,x,d}\)\(y\)\right| \to 0\quad\text{as }\left|y\right| \to \infty.
\end{equation}
By applying \cite{CL}*{ Lemma 2.4}, it follows from \eqref{Lem4Eq14} and \eqref{Lem4Eq15} that there exist $c_0, c_1, \dotsc, c_n\in\R$ such that 
$$v_{k,x,d}- \tilde{v}_{k,x,d} =  c_0 \( \frac{n-2}{2} U_0+ \< y, \nabla U_0 \> \) + \sum_{i=1}^n c_i \partial_{y_i} U_0 \quad \text{in } \R^n. $$ 
By using the explicit expression of $U_0$ together with \eqref{Lem4Eq5} and \eqref{Lem4Eq13}, we then obtain that $c_i = 0$ for all $i\in\left\{1,\dotsc,n\right\}$, i.e. $v_{k,x,d} = \tilde{v}_{k,x,d}$ in $\R^n$. The estimate \eqref{Lem4Eq12} then gives \eqref{Lem4Eq6} for $j=0$. Finally, by using \eqref{Lem4Eq4} together with standard elliptic theory, we obtain \eqref{Lem4Eq6} for $j \ge 1$. This ends the proof of Lemma~\ref{Lem4}.
\endproof

{It is useful for what follows to note that for each $k\in\N$, $x\in\S^n$ and $d\in\left\{4,\dotsc,n-4\right\}$, the function $v_{k,x,d}$ satisfies a certain orthogonality property with respect to $U_0$. Indeed, \eqref{Lem4Eq8} and \eqref{Lem4Eq9} imply that $\int_{\mathbb{S}^{n-1}} w_{k,x,d} \dsigma = 0$. Moreover, since $w_{k,x,d}$ is a polynomial of degree at most $n-6$, we have $U_0^2 w_{k,x,d} \in L^1\(\R^n\)$, and since $U_0$ is radially symmetric, a simple change of variables in polar coordinates gives $\int_{\R^n} U_0^2 w_{k,x,d} \dy = 0$. Integrating by parts \eqref{Lem4Eq4} against $U_0$ and using the equation $\Delta_\xi U_0 = U_0^{2^*-1}$ then yields
\begin{equation} \label{Lem4Eq15:bis}
\int_{\R^n} U_0^{2^*-1} v_{k,x,d} \dy = 0.
\end{equation}} 

\smallskip
From now on, we let $\(g_k\)_k$, $\(u_k\)_k$, $\(x_{1,k}\)_k$, $\(x_{2,k}\)_k$, $\(\mu_{1,k}\)_k$, $\(\mu_{2,k}\)_k$ and $\(d_k\)_k$ be as in Lemma~\ref{Lem1}, $\(\varphi_k\)_k$, $\(\hat{g}_{k,x}\)_{k,x}$, $\(\exp_{k,x}\)_{k,x}$, $\varepsilon_0$ be as in Lemma~\ref{Lem2}, $\(h_{k,x,\alpha}\)_{k,x,\alpha}$ and $d_n$ be as in Lemma~\ref{Lem3} and $\(v_{k,x,d}\)_{k,x,d}$ be as in Lemma~\ref{Lem4}. We now introduce some additional notations of radii and rescaled functions. For each $k\in\N$, we define
 \begin{equation} \label{Lem5Eq1}
 \varrho_{1,k}:=\frac{d_k}{\mu_{1,k}}\quad\text{and}\quad\varrho_{2,k}:=\sqrt{\frac{\mu_{1,k}}{\mu_{2,k}}+\frac{d_k^2}{\mu_{1,k}\mu_{2,k}}}.
 \end{equation}
For each $k\in\N$ and $i \in \left\{1,2\right\}$, we define 
\begin{equation} \label{Lem5Eq2}
\exp_{i,k}:=\exp_{k,x_{i,k}},\quad H_{i,k}:=H_{k,x_{i,k}}\quad\text{and}\quad h_{i,k,\alpha}:=h_{k,x_{i,k},\alpha}
\end{equation}
as well as the metric
\begin{equation} \label{Lem5Eq3}
\hat{g}_{i,k}:=\hat{g}_{k,x_{i,k}}\(\mu_{i,k}\cdot\)
\end{equation}
and the functions
\begin{align}
\hat{u}_{i,k}&:=\mu_{i,k}^{\frac{n-2}{2}}\big(\varphi_k\(x_{i,k},\cdot\)^{-1}u_k\big)\big(\exp_{i,k}\(\mu_{i,k}\cdot\)\big),\label{Lem5Eq4}\allowdisplaybreaks\\
\hat{w}_{i,k}&:=c_n\sum_{d=4}^{n-4}\mu_{i,k}^{d}w_{k,x_{i,k},d},\label{Lem5Eq5}\allowdisplaybreaks\\
\hat{v}_{i,k}&:=c_n\sum_{d=4}^{n-4}\mu_{i,k}^{d}v_{k,x_{i,k},d}\label{Lem5Eq6}
\end{align}
and
\begin{align*}
\check{w}_{i,k}&:=c_n\mu_{i,k}^2\sum_{a,b,c=1}^n\bigg(\partial_{y_b}\(\(H_{i,k}\)_{ab}\partial_{y_c}\(H_{i,k}\)_{ac}\)\\
&\quad-\frac{1}{2}\partial_{y_b}\(H_{i,k}\)_{ab}\partial_{y_c}\(H_{i,k}\)_{ac}+\frac{1}{4}\(\partial_{y_c}\(H_{i,k}\)_{ab}\)^2\bigg)\(\mu_{i,k}\cdot\).
\end{align*}
In the case where $3 \le n \le 7$, $\hat{w}_{i,k}=\hat{v}_{i,k}=0$ (see the discussion around \eqref{Lem4Eq3}).

\smallskip
By using \eqref{Lem3Eq6}, \eqref{Lem3Eq7} and \eqref{Lem4Eq6}, we obtain
\begin{align}
\left|c_n\mu_{i,k}^2\Scal_{\hat{g}_{i,k}}-\hat{w}_{i,k}+\check{w}_{i,k}\right|&=\bigO\(\sum_{\left|\alpha\right|=2}^{d_n}\left|h_{i,k,\alpha}\right|^2\mu_{i,k}^{2\left|\alpha\right|+2}\left|y\right|^{2\left|\alpha\right|}\)\nonumber\\
&\quad+\smallo\big(\mu_{i,k}^{n+1}\left|y\right|^{n-1}\big)\quad\text{as }k\to\infty,
\label{Lem5Eq7}\allowdisplaybreaks\\
\left|c_n\mu_{i,k}^2\Scal_{\hat{g}_{i,k}}-\hat{w}_{i,k}\right|&=\bigO\(\sum_{\left|\alpha\right|=2}^{d_n}\left|h_{i,k,\alpha}\right|^2\mu_{i,k}^{2\left|\alpha\right|}\left|y\right|^{2\left|\alpha\right|-2}\)\nonumber\\
&\quad+\smallo\big(\mu_{i,k}^2\(\mu_{i,k}\left|y\right|\)^{\max\(n-3,2\)}\big)\quad\text{as }k\to\infty\label{Lem5Eq8}
\end{align}
and
\begin{equation}\label{Lem5Eq9}
\left|\nabla^j \hat{v}_{i,k}\(y\)\right|=\bigO\(\sum_{\left|\alpha\right|=4}^{n-4}\frac{\left|h_{i,k,\alpha}\right|\mu_{i,k}^{\left|\alpha\right|}}{\(1+\left|y\right|\)^{n-\left|\alpha\right|-2+j}}\)\quad\forall j\in\N
\end{equation}
uniformly with respect to $y \in \R^n$ and $k\in\N$.

\smallskip
To prove Theorems~\ref{Th1} and \ref{Th2}, we need refined asymptotics on the functions $\(u_k\)_k$. In the following lemma, we improve the a priori estimates of Lemma~\ref{Lem1} and obtain a sharp description of $u_k + \(-1\)^i B_{i,k}$ near $x_{i,k}$, which depends on the local geometry of $g_k$ near this point. After scaling, this amounts to obtaining refined pointwise estimates on $\hat{u}_{i,k} + \(-1\)^i U_0$ in Euclidean balls of radii of order $\varrho_{i,k}$. The analysis of \cite{KMS} does not directly apply here since the functions $\(u_k\)_k$ change sign and, as a consequence, the blow-up points $\(x_{1,k}\)_k$ and $\(x_{2,k}\)_k$ are not isolated and simple (in particular, $d_k = \dist_{g_0}\(x_{1,k},x_{2,k}\)$ may tend to $0$ as $k \to \infty$). Our refined estimates are as follows:

\begin{lemma}\label{Lem5}
Let $\(g_k\)_k$, $\(u_k\)_k$, $\(x_{1,k}\)_k$, $\(x_{2,k}\)_k$, $\(\mu_{1,k}\)_k$, $\(\mu_{2,k}\)_k$ and $\(d_k\)_k$ be as in Lemma~\ref{Lem1}, $\(\varphi_k\)_k$, $\(\hat{g}_{k,x}\)_{k,x}$, $\(\exp_{k,x}\)_{k,x}$, $\varepsilon_0$ be as in Lemma~\ref{Lem2}, $\(h_{k,x,\alpha}\)_{k,x,\alpha}$ and $d_n$ be as in Lemma~\ref{Lem3} and $\(v_{k,x,d}\)_{k,x,d}$ be as in Lemma~\ref{Lem4}. Let $i\in\left\{1,2\right\}$ and $\(\varrho_{k}\)_{i,k}$, $\(\hat{u}_{i,k}\)_{k}$, $\(\hat{w}_{i,k}\)_{k}$ and $\(\hat{v}_{i,k}\)_{k}$ be as in \eqref{Lem5Eq1}, \eqref{Lem5Eq4}, \eqref{Lem5Eq5} and \eqref{Lem5Eq6}, respectively. In the case where $i=1$, assume that 
\begin{equation}\label{Lem5Eq10}
\mu_{1,k}=\smallo\(d_k\)\quad\text{(i.e. }\varrho_{1,k}\to\infty\text{)}\quad\text{as }k\to\infty
\end{equation}
and \eqref{Lem1Eq4} holds true (observe that \eqref{Lem5Eq10} implies \eqref{Lem1Eq3} since $\mu_{2,k}\le\mu_{1,k}$ for all $k\in\N$). In the case where $i=2$, we do not make any additional assumptions. Then there exist $\delta_0\in\(0,\varepsilon_0/\pi\)$ and $k_0\in\N$ such that
\begin{align}\label{Lem5Eq11}
&\sum_{j=0}^2\(1+\left|y\right|\)^j\big|\nabla^j\big(\hat{u}_{i,k}+\(-1\)^i\(U_0-\hat{v}_{i,k}\)\big)\(y\)\big|\nonumber\\
&\quad=\bigO\(\left\{\sum_{\left|\alpha\right|=2}^{d_n-1}\frac{\left|h_{i,k,\alpha}\right|^2\mu_{i,k}^{2\left|\alpha\right|}}{\(1+\left|y\right|\)^{n-2\left|\alpha\right|-2}}+\frac{\mu_{i,k}^{n-3}}{1+\left|y\right|}\quad\text{if }n\ge6\right\}+\varrho_{i,k}^{2-n}\)
\end{align}
uniformly with respect to $y\in \B_\xi\(0,\delta_0\varrho_{i,k}\)$ and $k>k_0$.
\end{lemma}

In the case where $n\in\left\{6,7\right\}$, the sum in the right-hand side of \eqref{Lem5Eq11} is empty.

\proof[Proof of Lemma~\ref{Lem5}]
We adapt the arguments from the proof of \cite{KMS}*{Proposition~5.1}, taking into account a general configuration for $\(x_{1,k}\)_k$ and $\(x_{2,k}\)_k$. In particular, we assume neither that $d_k \not \to 0$ as $k\to \infty$ nor that the blow-up points $\(x_{1,k}\)_k$ and $\(x_{2,k}\)_k$ are isolated and simple. Since $g_k\to g_0$ and $\varphi_k\(x_{i,k},x_{i,k}\)=1$, we obtain
\begin{align}
\dist_{g_0}\big(\exp_{i,k}\(y\),x_{i,k}\big)&=\(1+\smallo\(1\)\)\dist_{g_k}\big(\exp_{i,k}\(y\),x_{i,k}\big)\nonumber\\
&=\left|y\right|+\bigO\big(\left|y\right|^2\big)+\smallo\(\left|y\right|\)\quad\text{as }k\to\infty\label{Lem5Eq12}
\end{align}
uniformly with respect to $y$ in compact subsets of $\R^n$. Moreover, since $d_k\le\pi$ and $\mu_{2,k}\le\mu_{1,k}\to0$, we obtain
\begin{equation}\label{Lem5Eq13}
\mu_{1,k}\varrho_{1,k}=d_k\le\pi
\end{equation}
and
\begin{equation}\label{Lem5Eq14}
\mu_{2,k}\varrho_{2,k}=\sqrt{\mu_{1,k}\mu_{2,k}+\frac{\mu_{2,k}}{\mu_{1,k}}d_k^2}\le d_k+\bigO\(\mu_{1,k}\)\le\pi+\smallo\(1\)\quad\text{as }k\to\infty.
\end{equation}
Since $g_{k,x_{i,k}}\to g_{0,x_i}$ in $C^m\(\S^n\)$ as $k\to\infty$ for all $m\in\N$, it follows from \eqref{Lem5Eq13} and \eqref{Lem5Eq14} that there exist $\delta_0\in\(0, \varepsilon_0/\pi\)$ and $k_0\in\N$ such that, for each $k>k_0$, $\delta_0\mu_{i,k}\varrho_{i,k}$ is smaller than the injectivity radius at $x_{i,k}$ of the metric $g_{k,x_{i,k}}$ or, equivalently, $\delta_0\varrho_{i,k}$ is smaller than the injectivity radius at 0 of the metric $\hat{g}_{i,k}$ defined in \eqref{Lem5Eq3}.
By letting $k_0$ be larger if necessary, \eqref{Lem5Eq13} and \eqref{Lem5Eq14} also give
\begin{equation} \label{Lem5Eq15}
\mu_{i,k}\left|y\right| \le \varepsilon_0\quad\forall y \in \B_{\xi}\(0, \delta_0 \varrho_{i,k}\).
\end{equation}

\smallskip
We restrict ourselves to giving the proof of \eqref{Lem5Eq11} for $j=0$ as the estimates on the derivatives then follow by standard elliptic theory. Since $\varphi_k\to\varphi_0$ in $C^m\(\S^n \times \S^n\)$ for all $m\in\N$, $x_{i,k}\to x_i$ as $k \to \infty$ and $\delta_0<\varepsilon_0/\pi$, it follows from \eqref{Lem1Eq5}, \eqref{Lem2Eq1}, \eqref{Lem5Eq13} and \eqref{Lem5Eq14} that
\begin{align}\label{Lem5Eq16}
&\big|\big(\hat{u}_{i,k}+\(-1\)^iU_0\big)\(y\)\big|\nonumber\\
&\quad=\bigO\Big(\mu_{i,k}^{\frac{n-2}{2}}B_{3-i,k}\big(\exp_{i,k}\(\mu_{i,k}y\)\big)\Big)+\smallo\(U_0\(y\)\)\quad\text{as }k\to\infty
\end{align}
uniformly with respect to $y\in \B_\xi\(0,\delta_0\varrho_{i,k}\)$. Moreover, by using \eqref{Lem5Eq12}, \eqref{Lem5Eq13} and \eqref{Lem5Eq14}, we obtain
\begin{align}
&\dist_{g_0}\big(\exp_{i,k}\(\mu_{i,k}y\),x_{3-i,k}\big)\nonumber\\
&\quad\ge d_k-\dist_{g_0}\big(\exp_{i,k}\(\mu_{i,k}y\),x_{i,k}\big)\nonumber\\
&\quad\ge  \(1-\delta_0+\smallo\(1\)\)d_k+\bigO\big(\(\delta_0d_k\)^2+\delta_0\mu_{3-i,k} \big)\quad\text{as }k\to\infty \label{Lem5Eq17}
\end{align}
uniformly with respect to $y\in \B_\xi\(0,\delta_0\varrho_{i,k}\)$. By letting $\delta_0$ be smaller if necessary, it follows from \eqref{Lem5Eq17} that
\begin{equation}\label{Lem5Eq18}
\mu_{i,k}^{\frac{n-2}{2}}B_{3-i,k}\big(\exp_{i,k}\(\mu_{i,k}y\)\big)=\bigO\big(\varrho_{i,k}^{2-n}\big)
\end{equation}
uniformly with respect to $y\in \B_\xi\(0,\delta_0\varrho_{i,k}\)$ and $k>k_0$. By combining \eqref{Lem5Eq16} and \eqref{Lem5Eq18}, we obtain
\begin{equation}\label{Lem5Eq19}
\big|\big(\hat{u}_{i,k}+\(-1\)^iU_0\big)\(y\)\big|=\bigO\big(\varrho_{i,k}^{2-n}\big)+\smallo\(U_0\(y\)\)\quad\text{as }k\to\infty
\end{equation}
uniformly with respect to $y\in \B_\xi\(0,\delta_0\varrho_{i,k}\)$. Moreover, \eqref{Lem5Eq9} and \eqref{Lem5Eq15} give
\begin{equation}\label{Lem5Eq20}
\left|\hat{v}_{i,k}\(y\)\right|=\smallo\(U_0\(y\)\)\quad\text{as }k\to\infty
\end{equation}
uniformly with respect to $y\in \B_\xi\(0,\delta_0\varrho_{i,k}\)$. For each $k>k_0$, we define
$$m_{i,k}:=\max_{\overline{\B_\xi\(0,\delta_0\varrho_{i,k}\)}}\big|\big(\hat{u}_{i,k}+\(-1\)^i\(U_0-\hat{v}_{i,k}\)\big)\big|$$
and
$$\psi_{i,k}:=m_{i,k}^{-1}\big(\hat{u}_{i,k}+\(-1\)^i\(U_0-\hat{v}_{i,k}\)\big).$$
By using \eqref{Lem5Eq10} for $i=1$ and Lemma~\ref{Lem1}~(iii) for $i=2$, we obtain $\varrho_{i,k} \to \infty$ as $k \to \infty$ in both cases. By using \eqref{Lem5Eq19} and \eqref{Lem5Eq20}, we then obtain $m_{i,k}\to0$ as $k\to\infty$. By using the conformal invariance of the conformal Laplacian, we can rewrite \eqref{Sec3Eq3} as 
\begin{equation}\label{Lem5Eq21}
\Delta_{\hat{g}_{i,k}}\hat{u}_{i,k}+c_n\mu_{i,k}^2\Scal_{\hat{g}_{i,k}}\hat{u}_{i,k}=\left|\hat{u}_{i,k}\right|^{2^*-2}\hat{u}_{i,k}\quad\text{in }\B_\xi\(0,\delta_0\varrho_{i,k}\).
\end{equation}
Moreover, by using \eqref{Lem4Eq4}, we obtain
\begin{equation}\label{Lem5Eq22}
\Delta_\xi\hat{v}_{i,k}=\(2^*-1\)U_0^{2^*-2}\hat{v}_{i,k}+U_0\hat{w}_{i,k}\quad\text{in }\R^n.
\end{equation}
By using \eqref{Lem5Eq21} and \eqref{Lem5Eq22} together with the equation $\Delta_\xi U_0=U_0^{2^*-1}$, we obtain
\begin{equation}\label{Lem5Eq23}
\Delta_{\hat{g}_{i,k}} \psi_{i,k}+c_n\mu_{i,k}^2\Scal_{\hat{g}_{i,k}} \psi_{i,k} =\(2^*-1\)U_0^{2^*-2}\psi_{i,k}+m_{i,k}^{-1}f_{i,k}
\end{equation}
in $\B_\xi\(0,\delta_0\varrho_{i,k}\)$, where
\begin{align}\label{Lem5Eq24}
f_{i,k}&:=\(-1\)^i\(\Delta_{\hat{g}_{i,k}}-\Delta_\xi\)\(U_0-\hat{v}_{i,k}\)-\(-1\)^ic_n\mu_{i,k}^2\Scal_{\hat{g}_{i,k}}\hat{v}_{i,k}\nonumber\\
&\quad+\(-1\)^i U_0 \(c_n\mu_{i,k}^2\Scal_{\hat{g}_{i,k}}-\hat{w}_{i,k}\)+\left|\hat{u}_{i,k}\right|^{2^*-2}\hat{u}_{i,k}+\(-1\)^iU_0^{2^*-1}\nonumber\\
&\quad-\(2^*-1\) U_0^{2^*-2}\big(\(-1\)^i\hat{v}_{i,k}+m_{i,k}\psi_{i,k}\big) .
\end{align}
We now estimate the terms in the right-hand side of \eqref{Lem5Eq24}. Since $U_0$ is radially symmetric around 0, it follows from \eqref{Lem2Eq3} that 
\begin{equation}\label{Lem5Eq25}
\(\Delta_{\hat{g}_{i,k}}-\Delta_\xi\)U_0\(y\)=\bigO\big(\mu_{i,k}^N\left|y\right|^{N-1}\left|\nabla U_0\(y\)\right|\big)=\bigO\(\frac{\(\mu_{i,k}\left|y\right|\)^N}{\(1+\left|y\right|\)^n}\)
\end{equation}
uniformly with respect to $y\in \B_\xi\(0,\delta_0\varrho_{i,k}\)$ and $k>k_0$. We recall that, in the case where $3\le n\le 7$, $\hat{v}_{i,k}=0$ for all $k\in\N$. When $n\ge8$, by using Lemmas~\ref{Lem3} and \ref{Lem4} together with straightforward estimates and the fact that $n-4 \ge d_n$, we obtain
\begin{align}\label{Lem5Eq26}
&\(\Delta_{\hat{g}_{i,k}}-\Delta_\xi\)\hat{v}_{i,k}\(y\)-\(-1\)^ic_n\mu_{i,k}^2\Scal_{\hat{g}_{i,k}}\(y\)\hat{v}_{i,k}\(y\)\nonumber\\
&=\bigO\(\left|\nabla\hat{g}_{i,k}\(y\)\right|\left|\nabla\hat{v}_{i,k}\(y\)\right|+\left|\(\hat{g}_{i,k}-\xi\)\(y\)\right|\left|\nabla^2\hat{v}_{i,k}\(y\)\right|+\mu_{i,k}^2\left|\Scal_{\hat{g}_{i,k}}\(y\)\hat{v}_{i,k}\(y\)\right|\)\nonumber\allowdisplaybreaks\\
&=\bigO\left(\(\(1+\left|y\right|\)\left|\nabla\hat{g}_{i,k}\(y\)\right|+\left|\(\hat{g}_{i,k}-\xi\)\(y\)\right|+\mu_{i,k}^2\(1+\left|y\right|\)^2\left|\Scal_{\hat{g}_{i,k}}\(y\)\right|\)\right.\nonumber\\
&\quad\left.\times\(\(1+\left|y\right|\)^{-1}\left|\nabla\hat{v}_{i,k}\(y\)\right|+\left|\nabla^2\hat{v}_{i,k}\(y\)\right|+\(1+\left|y\right|\)^{-2}\left|\hat{v}_{i,k}\(y\)\right|\)\right)\nonumber\allowdisplaybreaks\\
&=\bigO\left(\(\sum_{\left|\alpha\right|=2}^{n-4}\left|h_{i,k,\alpha}\right|\mu_{i,k}^{\left|\alpha\right|}\left|y\right|^{\left|\alpha\right|}+\smallo\(\mu_{i,k}^{n-3}\left|y\right|^{n-3}\)\)\right.\nonumber\\
&\quad\times\left.\(\sum_{\left|\alpha\right|=4}^{n-4}\frac{\left|h_{i,k,\alpha}\right|\mu_{i,k}^{\left|\alpha\right|}}{\(1+\left|y\right|\)^{n-\left|\alpha\right|}}\)\right)\quad\text{as }k\to\infty
\end{align}
uniformly with respect to $y\in \B_\xi\(0,\delta_0\varrho_{i,k}\)$. Since $|y| \le 1+ |y|$, $\mu_{i,k}(1+|y|) \le 2 \varepsilon_0$ for $k$ large enough by \eqref{Lem5Eq15} and since $h_{i,k,\alpha}\to0$ and $2d_n\ge n-3\ge0$, by using Cauchy-Schwarz's inequality, it follows from \eqref{Lem5Eq26} that
\begin{align}\label{Lem5Eq27}
&\(\Delta_{\hat{g}_{i,k}}-\Delta_\xi\)\hat{v}_{i,k}\(y\)-\(-1\)^ic_n\mu_{i,k}^2\Scal_{\hat{g}_{i,k}}\(y\)\hat{v}_{i,k}\(y\)\nonumber\\
&\quad=\bigO\(\(1+\left|y\right|\)^{-n}\(\sum_{\left|\alpha\right|=2}^{n-4}\left|h_{i,k,\alpha}\right|\mu_{i,k}^{\left|\alpha\right|}\(1+\left|y\right|\)^{\left|\alpha\right|}\)^2 \)+\smallo\(\frac{\mu_{i,k}^{n-3}}{\(1+\left|y\right|\)^{3}}\)\nonumber\allowdisplaybreaks\\
&\quad=\bigO\(\sum_{\left|\alpha\right|=2}^{n-4}\frac{\left|h_{i,k,\alpha}\right|^2\mu_{i,k}^{2\left|\alpha\right|}}{\(1+\left|y\right|\)^{n-2\left|\alpha\right|}}\)+\smallo\(\frac{\mu_{i,k}^{n-3}}{\(1+\left|y\right|\)^{3}}\)\nonumber\allowdisplaybreaks\\
&\quad=\bigO\(\sum_{\left|\alpha\right|=2}^{d_n-1}\frac{\left|h_{i,k,\alpha}\right|^2\mu_{i,k}^{2\left|\alpha\right|}}{\(1+\left|y\right|\)^{n-2\left|\alpha\right|}}\) + \smallo\(\frac{\mu_{i,k}^{n-3}}{\(1+\left|y\right|\)^{3}}\)\quad\text{as }k\to\infty  
\end{align}
uniformly with respect to $y\in \B_\xi\(0,\delta_0\varrho_{i,k}\)$. Similarly, by using \eqref{Lem5Eq8}, \eqref{Lem5Eq15} and \eqref{Lem5Eq19}, we obtain 
\begin{align}\label{Lem5Eq28}
&U_0\(y\)\(c_n\mu_{i,k}^2\Scal_{\hat{g}_{i,k}} - \hat{w}_{i,k}\)\(y\)\nonumber\\
&\quad=\left\{\bigO\(\sum_{\left|\alpha\right|=2}^{d_n-1}\frac{\left|h_{i,k,\alpha}\right|^2\mu_{i,k}^{2\left|\alpha\right|}}{\(1+\left|y\right|\)^{n-2\left|\alpha\right|}}\)+\smallo\(\frac{\mu_{i,k}^{n-3}}{\(1+\left|y\right|\)^{3}}\)\quad\text{if }n\ge6\right\}\nonumber\\
&\qquad+\frac{\mu_{i,k}^2\(\mu_{i,k}\left|y\right|\)^{\max\(n-3.2\)}}{\(1+\left|y\right|\)^{n-2}}\quad\text{as }k\to\infty
\end{align}
and
\begin{align}\label{Lem5Eq29}
&\(\left|\hat{u}_{i,k}\right|^{2^*-2}\hat{u}_{i,k}+\(-1\)^iU_0^{2^*-1}-\(2^*-1\)U_0^{2^*-2}\big(\(-1\)^i\hat{v}_{i,k}+m_{i,k}\psi_{i,k}\big)\)\(y\)\nonumber\\
&\quad=\bigO\(U_0\(y\)^{2^*-3}\big(\hat{v}_{i,k}\(y\)^2+\(m_{i,k}\psi_{i,k}\(y\)\)^2\big)\)\nonumber\allowdisplaybreaks\\
&\quad=\left\{\bigO\(\sum_{\left|\alpha\right|=4}^{d_n-1}\frac{\left|h_{i,k,\alpha}\right|^2\mu_{i,k}^{2\left|\alpha\right|}}{\(1+\left|y\right|\)^{n-2\left|\alpha\right|}}\)+\smallo\(\frac{\mu_{i,k}^{n-3}}{\(1+\left|y\right|\)^{3}}\)\quad\text{if }n\ge6\right\}\nonumber\\
&\qquad+\bigO\(\frac{\varrho_{i,k}^{4-2n}}{\(1+\left|y\right|\)^{6-n}}\)+\smallo\(m_{i,k}U_0\(y\)^{2^*-2}\left|\psi_{i,k}\(y\)\right|\)\quad\text{as }k\to\infty
\end{align}
uniformly with respect to $y\in \B_\xi\(0,\delta_0\varrho_{i,k}\)$. We now let $\widehat{G}_{i,k}$ be the Green's function of $\Delta_{\hat{g}_{i,k}} +c_n\mu_{i,k}^2\Scal_{\hat{g}_{i,k}}$ in $\B_\xi\(0,\delta_0\varrho_{i,k}\)$ with zero Dirichlet boundary condition on $\partial \B_\xi\(0,\delta_0\varrho_{i,k}\)$. By definition of $\hat{g}_{i,k}$, it is easy to see that
$$ \widehat{G}_{i,k}\(x,y\) = \mu_{i,k}^{n-2} \widetilde{G}_{i,k}\( \mu_{i,k}x, \mu_{i,k}y \)\quad\forall x,y\in \B_{\xi}\(\delta_0 \varrho_{i,k}\),\ x\ne y,$$
where $\widetilde{G}_{i,k}$ is the Green's function of  $\Delta_{\hat{g}_{k,x_{i,k}}} +c_n\Scal_{\hat{g}_{k,x_{i,k}}} $ with Dirichlet boundary condition in $\B_{\xi}\(0, \delta_0 \mu_{i,k} \varrho_{i,k}\)$. Since $g_k\to g_0$ in $C^m\(\S^n\)$ as $k\to\infty$ for all $m\in\N$, it follows from \eqref{Lem2Eq4} and \eqref{Lem5Eq15} that $\hat{g}_{k,x_{i,k}}\to\xi$ as $k\to\infty$ in $C^m\(\B_{\xi}\(0,\varepsilon_0\)\)$ for all $m\in\N$. By using standard estimates for $\widetilde{G}_{i,k}$, which can be found for instance in \cite{R}, it is not difficult to see that $\widehat{G}_{i,k}$ satisfies
\begin{equation}\label{Lem5Eq30}
\widehat{G}_{i,k}\(y,z\)\le C\left|y-z\right|^{2-n}\quad\forall y,z\in \B_\xi\(0,\delta_0\varrho_{i,k}\)
\end{equation}
and
\begin{equation}\label{Lem5Eq31}
\big|\partial_\nu \widehat{G}_{i,k}\(y,z\)\big|\le C\left|y-z\right|^{1-n}\quad\forall y\in \B_\xi\(0,\delta_0\varrho_{i,k}\),\,z\in\partial \B_\xi\(0,\delta_0\varrho_{i,k}\)
\end{equation}
for some constant $C$ independent of $k$. A representation formula for \eqref{Lem5Eq23} now gives
\begin{align}\label{Lem5Eq32}
\psi_{i,k}\(y\)&=\int_{\B_\xi\(0,\delta_0\varrho_{i,k}\)}\widehat{G}_{i,k}\(y,\cdot\)\big(\(2^*-1\)U_0^{2^*-2}\psi_{i,k}+m_{i,k}^{-1}f_{i,k}\big)\dv_{\hat{g}_{i,k}}\nonumber\\
&\quad-\int_{\partial \B_\xi\(0,\delta_0\varrho_{i,k}\)}\partial_\nu \widehat{G}_{i,k}\(y,\cdot\)\psi_{i,k}\dsigma_{\hat{g}_{i,k}}
\end{align}
for all $y\in \B_\xi\(0,\delta_0\varrho_{i,k}\)$, where $\nu$ and $\dsigma_{\hat{g}_{i,k}}$ are the outward unit normal vector and volume element, respectively, induced by $\hat{g}_{i,k}$ on $\partial \B_\xi\(0,\delta_0\varrho_{i,k}\)$. We observe that \eqref{Lem5Eq19} and \eqref{Lem5Eq20} give
\begin{equation}\label{Lem5Eq33}
\max_{\B_\xi\(0,\delta_0\varrho_{i,k}\) \setminus \B_\xi\(0,\delta_0\varrho_{i,k}/2\)}\left|\psi_{i,k}\right|=\bigO\big(m_{i,k}^{-1}\varrho_{i,k}^{2-n}\big)
\end{equation}
uniformly with respect to $k>k_0$. First considering the case where $\left|y\right| \le\delta_0\varrho_{i,k}/2$, by putting together \eqref{Lem5Eq25}, \eqref{Lem5Eq27}, \eqref{Lem5Eq28}, \eqref{Lem5Eq29},   \eqref{Lem5Eq30}, \eqref{Lem5Eq31}, \eqref{Lem5Eq32} and \eqref{Lem5Eq33} and using straightforward integral estimates, we obtain
\begin{align}
\psi_{i,k}\(y\)&=\bigO\(\int_{\B_\xi\(0,\delta_0\varrho_{i,k}\)}\frac{\left|\psi_{i,k}\(z\)\right|\dz}{\left|y-z\right|^{n-2}\(1+\left|z\right|\)^4}+m_{i,k}^{-1}\(\left\{\sum_{\left|\alpha\right|=2}^{d_n-1}\frac{\left|h_{i,k,\alpha}\right|^2\mu_{i,k}^{2\left|\alpha\right|}}{\(1+\left|y\right|\)^{n-2\left|\alpha\right|-2}}\right.\right.\right.\nonumber\\
&\quad\left.\left.\left.+\frac{\mu_{i,k}^{n-3}}{1+\left|y\right|}\quad\text{if }n\ge6\right\}+\varrho_{i,k}^{2-n}\(1+\(\mu_{i,k}\varrho_{i,k}\)^N\)\)\)\nonumber\allowdisplaybreaks\\
&=\bigO\(\int_{\B_\xi\(0,\delta_0\varrho_{i,k}\)}\frac{\left|\psi_{i,k}\(z\)\right|\dz}{\left|y-z\right|^{n-2}\(1+\left|z\right|\)^4}+m_{i,k}^{-1}\(\left\{\sum_{\left|\alpha\right|=2}^{d_n-1}\frac{\left|h_{i,k,\alpha}\right|^2\mu_{i,k}^{2\left|\alpha\right|}}{\(1+\left|y\right|\)^{n-2\left|\alpha\right|-2}}\right.\right.\right.\nonumber\\
&\quad\left.\left.\left.+\frac{\mu_{i,k}^{n-3}}{1+\left|y\right|}\quad\text{if }n\ge6\right\}+\varrho_{i,k}^{2-n}\)\)\label{Lem5Eq34}
\end{align}
uniformly with respect to $y\in \B_\xi\(0,\delta_0\varrho_{i,k}/2\)$ and $k>k_0$. It follows from \eqref{Lem5Eq33} that \eqref{Lem5Eq34} actually remains true for $y\in \B_\xi\(0,\delta_0\varrho_{i,k}\)$. We now claim that
\begin{equation}\label{Lem5Eq35}
m_{i,k}=\bigO\(\left\{\sum_{\left|\alpha\right|=2}^{d_n-1}\left|h_{i,k,\alpha}\right|^2\mu_{i,k}^{2\left|\alpha\right|}+\mu_{i,k}^{n-3}\quad\text{if }n\ge6\right\}+\varrho_{i,k}^{2-n}\)
\end{equation}
uniformly with respect to $k>k_0$. We assume by contradiction that \eqref{Lem5Eq35} does not hold true, i.e. there exists a subsequence $\(k_j\)_{j\in\N}$ such that $k_j\to\infty$ and
\begin{equation}\label{Lem5Eq36}
\left\{\sum_{\left|\alpha\right|=2}^{d_n-1}\left|h_{i,k_j,\alpha}\right|^2\mu_{i,k_j}^{2\left|\alpha\right|}+\mu_{i,k_j}^{n-3}\quad\text{if }n\ge6\right\}+\varrho_{i,k_j}^{2-n}=\smallo\(m_{i,k_j}\)\quad\text{as }j\to\infty.
\end{equation}
Since $\left|\psi_{i,k_j}\right|\le1$ in $\B_\xi\big(0,\delta_0\varrho_{i,k_j}\big)$ and $\hat{g}_{i,k_j}\to\xi$ in $C^m_{\loc}\(\R^n\)$ as $j\to\infty$, it follows from \eqref{Lem5Eq23}, \eqref{Lem5Eq25}, \eqref{Lem5Eq27}, \eqref{Lem5Eq28} and \eqref{Lem5Eq29} and standard elliptic estimates that, up to a subsequence, $\(\psi_{i,k_j}\)_j$ converges in $C^1_{\loc}\(\R^n\)$ as $j\to\infty$ to a solution $\psi_0\in C^\infty\(\R^n\)$ of the equation
\begin{equation}\label{Lem5Eq37}
\Delta_\xi\psi_0=\(2^*-1\)U_0^{2^*-2}\psi_0\quad\text{in }\R^n.
\end{equation}
On the other hand, by using \eqref{Lem5Eq34} and \eqref{Lem5Eq36}, we obtain
\begin{equation}\label{Lem5Eq38}
\psi_{i,k_j}\(y\)=\bigO\big(\(1+\left|y\right|\)^{-2}\big)+\smallo\(1\)\quad\text{as }j\to\infty
\end{equation}
uniformly with respect to $y\in \B_\xi\big(0,\delta_0\varrho_{i,k_j}\big)$. By passing to the limit as $j\to\infty$ into \eqref{Lem5Eq38}, we then obtain
\begin{equation}\label{Lem5Eq39}
\psi_0\(y\)=\bigO\big(\(1+\left|y\right|\)^{-2}\big)
\end{equation}
uniformly with respect to $y\in\R^n$. By applying Lemma 2.4 in \cite{CL}, it follows from \eqref{Lem5Eq37} and \eqref{Lem5Eq39} that
\begin{equation}\label{Lem5Eq40}
\psi_0\(y\)=\lambda_0\frac{1-\frac{\left|y\right|^2}{n\(n-2\)}}{\(1+\frac{\left|y\right|^2}{n\(n-2\)}\)^{\frac{n}{2}}}+\sum_{i=1}^n\lambda_i\frac{y_i}{\(1+\frac{\left|y\right|^2}{n\(n-2\)}\)^{\frac{n}{2}}}\quad\forall y\in\R^n
\end{equation}
for some $\lambda_0,\dotsc,\lambda_n\in\R$. On the other hand, by using \eqref{Lem1Eq2}, \eqref{Lem1Eq4}, \eqref{Lem2Eq2}, and \eqref{Lem4Eq5}, we obtain $\psi_{i,k}\(0\)=\left|\nabla\psi_{i,k}\(0\)\right|=0$ for all $k\in\N$, which gives $\psi_0\(0\)=\left|\nabla\psi_0\(0\)\right|=0$. It then follows from \eqref{Lem5Eq39} and \eqref{Lem5Eq40} that $\lambda_0=\dotsb=\lambda_n=0$, and so $\psi_0=0$. Independently, for each $k>k_0$, by definition of $m_{i,k}$, there exists a point $y_{i,k}\in\overline{\B_\xi\(0,\delta_0\varrho_{i,k}\)}$ such that $\left|\psi_{i,k}\(y_{i,k}\)\right|=1$. It follows from \eqref{Lem5Eq38} that $\big(y_{i,k_j}\big)_j$ is bounded. This is in contradiction with the fact that $\psi_{i,k_j}\to\psi_0=0$ as $j\to\infty$ uniformly in compact subsets of $\R^n$. This proves that \eqref{Lem5Eq35} holds true. Finally, \eqref{Lem5Eq11} follows from \eqref{Lem5Eq35} together with successive iterations of \eqref{Lem5Eq34}. This ends the proof of Lemma~\ref{Lem5}.
\endproof

By using \eqref{Lem5Eq9}, \eqref{Lem5Eq11} and \eqref{Lem5Eq15} and observing that $h_{k,x_{i,k}} \to 0$ in $C^m_{\loc}\(\R^n\)$ as $k\to\infty$ for all $m\in\N$ and $d_n\le n-4$ when $n \ge 6$, we obtain
\begin{align} \label{Lem5Eq41}
& \sum_{j=0}^2\(1+\left|y\right|\)^j\big| \nabla^j \big( \hat{u}_{i,k} + \(-1\)^i  U_0\big)\(y\) \big|\nonumber\\
&\quad = \bigO\(\left\{\sum_{\left|\alpha\right|=2}^{n-4}\frac{\left|h_{i,k,\alpha}\right|\mu_{i,k}^{\left|\alpha\right|}}{\(1+\left|y\right|\)^{n-2-\left|\alpha\right|}}+\frac{\mu_{i,k}^{n-3}}{1+\left|y\right|}\quad{\text{if }n\ge8}\right\}+\varrho_{i,k}^{2-n}\)
\end{align}
uniformly with respect to $y \in \B_{\xi}\(0, \delta_0 \varrho_{i,k}\)$ and $k>k_0$. 
Here again, in the case where $n\in\left\{6,7\right\}$, the sum in the right-hand side of \eqref{Lem5Eq41} is empty. We frequently use this estimate in what follows.

\smallskip
We now write a suitable Pohozaev-type identity which we repeatedly use in the rest of our analysis:

\begin{lemma}\label{Lem6}
Let $\(\varrho_{2,k}\)_k$, $\(\hat{u}_{2,k}\)_k$ and $\(\hat{g}_{2,k}\)_k$ be as in \eqref{Lem5Eq1}, \eqref{Lem5Eq4} and \eqref{Lem5Eq3}, respectively, and $k_0$ and $\delta_0$ be as in Lemma~\ref{Lem5}. Then, for each $\delta\in\(0,\delta_0\)$ and $k>k_0$,
\begin{align}\label{Lem6Eq}
&\int_{\B_\xi\(0,\delta\varrho_{2,k}\)}\hspace{-2pt}\(\<\nabla\hat{u}_{2,k},\cdot\>_\xi+\frac{n-2}{2}\hat{u}_{2,k}\)\(\(\Delta_{\hat{g}_{2,k}}-\Delta_\xi\)\hat{u}_{2,k}+c_n\mu_{2,k}^2\Scal_{\hat{g}_{2,k}}\hat{u}_{2,k}\) \dy\nonumber\\
&\quad=\int_{\partial \B_\xi\(0,\delta\varrho_{2,k}\)}\bigg(\frac{n-2}{2}\hat{u}_{2,k}\partial_\nu\hat{u}_{2,k}+\delta\varrho_{2,k}\(\partial_\nu\hat{u}_{2,k}\)^2-\frac{\delta\varrho_{2,k}}{2}\left|\nabla\hat{u}_{2,k}\right|_\xi^2\nonumber\\
&\qquad+\frac{\delta\varrho_{2,k}}{2^*}\left|\hat{u}_{2,k}\right|^{2^*}\bigg)\dsigma,
\end{align}
where $\nu$ and $\dsigma$ are the outward unit normal vector and volume element, respectively, of the metric induced by $\xi$ on $\partial \B_\xi\(0,\delta\varrho_{2,k}\)$.
\end{lemma}

\proof[Proof of Lemma~\ref{Lem6}]
Pohozaev identities of this type have long been used in the asymptotic analysis of solutions to Yamabe-type equations. The proof consists of multiplying \eqref{Lem5Eq21} for $i=2$ by the function $\<\nabla\hat{u}_{2,k},\cdot\>_\xi$, integrating over the ball $\B_\xi\(0,\delta\varrho_{2,k}\)$ and performing clever integration by parts to obtain the boundary terms in the right-hand side of \eqref{Lem6Eq}. The details can be found, for instance, in \cite{Mar}*{Lemma~2.1}.
\endproof

We first estimate the boundary term of \eqref{Lem6Eq}. We obtain the following:

\begin{lemma}\label{Lem7}
Let $\(\varrho_{2,k}\)_k$, $\(\hat{u}_{2,k}\)_k$ and $\(\hat{g}_{2,k}\)_k$ be as in \eqref{Lem5Eq1}, \eqref{Lem5Eq4} and \eqref{Lem5Eq3}, respectively. Then
\begin{multline}\label{Lem7Eq1}
\liminf_{\delta\to0}\liminf_{k\to\infty}\Bigg(\varrho_{2,k}^{n-2}\int_{\partial \B_\xi\(0,\delta\varrho_{2,k}\)}\bigg(\frac{n-2}{2}\hat{u}_{2,k}\partial_\nu\hat{u}_{2,k}+\delta\varrho_{2,k}\(\partial_\nu\hat{u}_{2,k}\)^2\\
-\frac{\delta\varrho_{2,k}}{2}\left|\nabla\hat{u}_{2,k}\right|_\xi^2+\frac{\delta\varrho_{2,k}}{2^*}\left|\hat{u}_{2,k}\right|^{2^*}\bigg)\dsigma \Bigg)>0.
\end{multline}
\end{lemma}

\proof[Proof of Lemma~\ref{Lem7}]
Let $0< \delta \le \delta_0$ be fixed. In the rest of the proof, by passing if necessary to a subsequence in $k$, we assume
\begin{align}\label{Lem7Eq2}
&\lim_{k \to \infty}\Bigg(\varrho_{2,k}^{n-2}\int_{\partial \B_\xi\(0,\delta\varrho_{2,k}\)}\bigg(\frac{n-2}{2}\hat{u}_{2,k}\partial_\nu\hat{u}_{2,k}+\delta\varrho_{2,k}\(\partial_\nu\hat{u}_{2,k}\)^2-\frac{\delta\varrho_{2,k}}{2}\left|\nabla\hat{u}_{2,k}\right|_\xi^2\nonumber\\
&\ \ \quad +\frac{\delta\varrho_{2,k}}{2^*}\left|\hat{u}_{2,k}\right|^{2^*}\bigg)\dsigma\Bigg)\nonumber\\
&\ \ = \liminf_{k \to \infty}\Bigg(\varrho_{2,k}^{n-2}\int_{\partial \B_\xi\(0,\delta\varrho_{2,k}\)}\bigg(\frac{n-2}{2}\hat{u}_{2,k}\partial_\nu\hat{u}_{2,k}+\delta\varrho_{2,k}\(\partial_\nu\hat{u}_{2,k}\)^2-\frac{\delta\varrho_{2,k}}{2}\left|\nabla\hat{u}_{2,k}\right|_\xi^2\nonumber\\
&\ \ \quad +\frac{\delta\varrho_{2,k}}{2^*}\left|\hat{u}_{2,k}\right|^{2^*}\bigg)\dsigma\Bigg) 
\end{align}
By letting
$$\check{u}_{2,k}\(y\):=\varrho_{2,k}^{n-2}\hat{u}_{2,k}\(\varrho_{2,k}y\)\quad\text{and}\quad\check{g}_{2,k}\(y\):=\hat{g}_{2,k}\(\varrho_{2,k}y\)\quad\forall y\in\R^n,$$
we obtain
\begin{align}\label{Lem7Eq3}
&\int_{\partial \B_\xi\(0,\delta\varrho_{2,k}\)}\bigg(\frac{n-2}{2}\hat{u}_{2,k}\partial_\nu\hat{u}_{2,k}+\delta\varrho_{2,k}\(\partial_\nu\hat{u}_{2,k}\)^2-\frac{\delta\varrho_{2,k}}{2}\left|\nabla\hat{u}_{2,k}\right|_\xi^2\nonumber\\
&\qquad+\frac{\delta\varrho_{2,k}}{2^*}\left|\hat{u}_{2,k}\right|^{2^*}\bigg)\dsigma\nonumber\\
&\quad=\varrho_{2,k}^{2-n}\int_{\partial \B_\xi\(0,\delta\)}\bigg(\frac{n-2}{2}\check{u}_{2,k}\partial_\nu\check{u}_{2,k}+\delta\(\partial_\nu\check{u}_{2,k}\)^2-\frac{\delta}{2}\left|\nabla\check{u}_{2,k}\right|_\xi^2\nonumber\\
&\qquad+\frac{\delta\varrho_{2,k}^{-2}}{2^*}\left|\check{u}_{2,k}\right|^{2^*}\bigg)\dsigma.
\end{align}
By recalling \eqref{Lem1Eq5}, \eqref{Lem2Eq1}, \eqref{Lem5Eq12} and \eqref{Lem5Eq14} and since $\delta_0<\varepsilon_0/\pi$, $\varphi_k\to\varphi_0$ in $C^m\(\S^n \times \S^n\)$ for all $m\in\N$, $x_{2,k}\to x_2$, $\mu_{2,k}\le\mu_{1,k}\to0$ and $\varrho_{2,k}\to\infty$ as $k \to \infty$, we obtain
\begin{align}\label{Lem7Eq4}
&\mu_{2,k}^{\frac{n-2}{2}}\varrho_{2,k}^{n-2}\big(\varphi_k\(x_{2,k},\cdot\)^{-1}B_{2,k}\big)\big(\exp_{2,k}\(\mu_{2,k}\varrho_{2,k}y\)\big)\nonumber\\
&\quad=\(\frac{\sqrt{n\(n-2\)}\(\mu_{2,k}\varrho_{2,k}\)^2}{\mu_{2,k}^2+\(\mu_{2,k}\varrho_{2,k}\left|y\right|\)^2}\)^{\frac{n-2}{2}}\(1+\smallo\(1\)\)\nonumber\\
&\quad=\(\frac{\sqrt{n\(n-2\)}}{\left|y\right|^2}\)^{\frac{n-2}{2}}\(1+\smallo\(1\)\)\quad\text{as }k\to\infty
\end{align}
and
\begin{align}\label{Lem7Eq5}
&\mu_{2,k}^{\frac{n-2}{2}}\varrho_{2,k}^{n-2}\big(\varphi_k\(x_{2,k},\cdot\)^{-1}B_{1,k}\big)\big(\exp_{2,k}\(\mu_{2,k}\varrho_{2,k}y\)\big)\nonumber\\
&\quad=\(\frac{\(2+\smallo\(1\)\)\sqrt{n\(n-2\)}\mu_{1,k}\mu_{2,k}\varrho_{2,k}^2}{2\mu_{1,k}^2+\big(4-\mu_{1,k}^2\big)\(1-\cos\(\dist_{g_0}\(x_{1,k},x_{2,k}\)+\bigO\(\mu_{2,k}\varrho_{2,k}\left|y\right|\)\)\)}\)^{\frac{n-2}{2}}\nonumber\\
&\quad=\ell_0+\bigO\(\left|y\right|\)+\smallo\(1\)\quad\text{as }k\to\infty
\end{align}
uniformly with respect to $y\in \partial \B_\xi\(0,\delta\)$, where
$$\ell_0:=\left\{\begin{aligned}&\(\frac{\sqrt{n\(n-2\)}\dist_{g_0}\(x_1,x_2\)^2}{2\(1-\cos\(\dist_{g_0}\(x_1,x_2\)\)\)}\)^{\frac{n-2}{2}}&&\text{if }\dist_{g_0}\(x_1,x_2\)>0\\&\(n\(n-2\)\)^{\frac{n-2}{4}}&&\text{if }\dist_{g_0}\(x_1,x_2\)=0.\end{aligned}\right.$$
It follows from \eqref{Lem1Eq5}, \eqref{Lem7Eq4} and \eqref{Lem7Eq5} that
\begin{equation}\label{Lem7Eq6}
\check{u}_{2,k}\(y\)=\ell_0-\(\frac{\sqrt{n\(n-2\)}}{\left|y\right|^2}\)^{\frac{n-2}{2}}\(1+\smallo\(1\)\)+\bigO\(\left|y\right|\)\quad\text{as }k\to\infty
\end{equation}
uniformly with respect to $y\in\partial \B_\xi\(0,\delta\)$. Notice in particular that the $\bigO\(\left|y\right|\)$ in \eqref{Lem7Eq5} and \eqref{Lem7Eq6} are independent of $\delta$. Independently, by using \eqref{Lem5Eq11} together with standard elliptic theory, we obtain that up to a subsequence,
$$\check{u}_{2,k} \to \check{u}_{2,0}\quad\text{in }C^1_{\loc}\big(\B_{\xi}\(0, \delta_0\)\setminus \left\{0\right\}\big)\quad\text{as }k\to\infty$$
for some function $\check{u}_{2,0}\in C^1\(\B_{\xi}\(0, \delta_0\)  \setminus \left\{0\right\}\)$. Since $0 < \delta \le \delta_0$ is fixed, we then obtain
\begin{multline}\label{Lem7Eq7}
\lim_{k \to \infty}\int_{\partial \B_\xi\(0,\delta\)}\bigg(\frac{n-2}{2}\check{u}_{2,k}\partial_\nu\check{u}_{2,k}+\delta\(\partial_\nu\check{u}_{2,k}\)^2-\frac{\delta}{2}\left|\nabla\check{u}_{2,k}\right|_\xi^2+\frac{\delta\varrho_{2,k}^{-2}}{2^*}\left|\check{u}_{2,k}\right|^{2^*}\bigg)\dsigma\\
= \int_{\partial \B_\xi\(0,\delta\)}\bigg(\frac{n-2}{2}\check{u}_{2,0}\partial_\nu\check{u}_{2,0}+\delta\(\partial_\nu\check{u}_{2,0}\)^2-\frac{\delta}{2}\left|\nabla\check{u}_{2,0}\right|_\xi^2\bigg)\dsigma
\end{multline}
along this subsequence. On the other hand, passing to the limit into \eqref{Lem7Eq6} along this subsequence gives
\begin{equation} \label{Lem7Eq8}
\check{u}_{2,0}\(y\)= \ell_0 -\(\frac{\sqrt{n\(n-2\)}}{\left|y\right|^2}\)^{\frac{n-2}{2}}+\bigO\(\left|y\right|\)
\end{equation}
uniformly with respect to $y\in \partial \B_\xi\(0,\delta\)$. We now use \eqref{Lem7Eq8} to estimate the limit boundary integral in \eqref{Lem7Eq7}: straightforward computations show that there exists a constant $C$ independent of $\delta$ such that 
\begin{multline} \label{Lem7Eq9}
\Bigg|\int_{\partial \B_\xi\(0,\delta\)}\bigg(\frac{n-2}{2}\check{u}_{2,0}\partial_\nu\check{u}_{2,0}+\delta\(\partial_\nu\check{u}_{2,0}\)^2-\frac{\delta}{2}\left|\nabla\check{u}_{2,0}\right|_\xi^2\bigg)\dsigma\\
-\frac{1}{2}n^{\frac{n-2}{4}}\(n-2\)^{\frac{n+6}{4}}\omega_{n-1}\ell_0\Bigg| \le C \delta,
\end{multline}
where $\omega_{n-1}$ is the volume of the round $\(n-1\)$-sphere. Finally, by combining \eqref{Lem7Eq2}, \eqref{Lem7Eq3}, \eqref{Lem7Eq7} and \eqref{Lem7Eq9}, we obtain
\begin{align*}
& \liminf_{k \to \infty}\Bigg(\varrho_{2,k}^{n-2}\int_{\partial \B_\xi\(0,\delta\varrho_{2,k}\)}\bigg(\frac{n-2}{2}\hat{u}_{2,k}\partial_\nu\hat{u}_{2,k}+\delta\varrho_{2,k}\(\partial_\nu\hat{u}_{2,k}\)^2 -\frac{\delta\varrho_{2,k}}{2}\left|\nabla\hat{u}_{2,k}\right|_\xi^2\\
& \qquad +\frac{\delta\varrho_{2,k}}{2^*}\left|\hat{u}_{2,k}\right|^{2^*}\bigg)\dsigma\Bigg) \ge  \frac{1}{2}n^{\frac{n-2}{4}}\(n-2\)^{\frac{n+6}{4}}\omega_{n-1}\ell_0- C \delta,
\end{align*}
which, by considering the limit inferior as $\delta \to 0$, gives \eqref{Lem7Eq1}.
\endproof

Now considering the interior term of \eqref{Lem6Eq}, we obtain the following:

\begin{lemma}\label{Lem8}
Let $\(\varrho_{2,k}\)_k$, $\(\hat{u}_{2,k}\)_k$ and $\(\hat{g}_{2,k}\)_k$ be as in \eqref{Lem5Eq1}, \eqref{Lem5Eq4} and \eqref{Lem5Eq3}, respectively, and $k_0$ and $\delta_0$ be as in Lemma~\ref{Lem5}. Then
\begin{align}\label{Lem8Eq1}
&\int_{\B_\xi\(0,\delta\varrho_{2,k}\)}\hspace{-2pt}\(\<\nabla\hat{u}_{2,k},\cdot\>_\xi+\frac{n-2}{2}\hat{u}_{2,k}\)\(\(\Delta_{\hat{g}_{2,k}}-\Delta_\xi\)\hat{u}_{2,k}+c_n\mu_{2,k}^2\Scal_{\hat{g}_{2,k}}\hat{u}_{2,k}\)\dy\nonumber\\
&\quad=\bigO\(\sum_{\left|\alpha\right|=2}^{d_n}\left|h_{2,k,\alpha}\right|^2\mu_{2,k}^{2\left|\alpha\right|} \left|\ln\mu_{2,k}\right|^{\vartheta\(2\left|\alpha\right|,n-2\)}
+\delta\varrho_{2,k}^{2-n}\)
\end{align}
uniformly with respect to $\delta\in\(0,\delta_0\)$ and $k>k_0$, where
$$\vartheta\(s,t\):=\left\{\begin{aligned}&1&&\text{if }s=t\\&0&&\text{if }s\ne t\end{aligned}\right.\quad\forall s,t\in\R.$$
\end{lemma}

Remark that by definition of $\vartheta$, the term involving $\left|\ln\mu_{2,k}\right|$ only appears when $n$ is even and $\left|\alpha\right|= \frac{n-2}{2}$.

\proof[Proof of Lemma~\ref{Lem8}]
By writing $\hat{u}_{2,k}=\hat{u}_{2,k}+U_0-U_0$ and using \eqref{Lem2Eq3}, we obtain
\begin{align}\label{Lem8Eq2}
\left|\(\Delta_{\hat{g}_{2,k}}-\Delta_\xi\)\hat{u}_{2,k}\(y\)\right|&\le\left|\(\Delta_{\hat{g}_{2,k}}-\Delta_\xi\)\(\hat{u}_{2,k}+U_0\)\(y\)\right|+\left|\(\Delta_{\hat{g}_{2,k}}-\Delta_\xi\)U_0\(y\)\right|\nonumber\\
&=\bigO\(\left|\nabla\hat{g}_{2,k}\(y\)\right|\left|\nabla\(\hat{u}_{2,k}+U_0\)\(y\)\right|+\left|\(\hat{g}_{2,k}-\xi\)\(y\)\right|\right.\nonumber\\
&\quad\left.\times\left|\nabla^2\(\hat{u}_{2,k}+U_0\)\(y\)\right|+\mu_{2,k}^N\left|y\right|^{N-1}\left|\nabla U_0\(y\)\right|\)
\end{align}
uniformly with respect to $y\in \B_\xi\(0,\delta_0\varrho_{2,k}\)$ and $k>k_0$. On the other hand, by writing $\hat{u}_{2,k}=-U_0+\hat{u}_{2,k}+U_0$ and $c_n\mu_{2,k}^2\Scal_{\hat{g}_{2,k}}=\hat{w}_{2,k}+c_n\mu_{2,k}^2\Scal_{\hat{g}_{2,k}}-\hat{w}_{2,k}$, we obtain
\begin{align}\label{Lem8Eq3}
\<\nabla\hat{u}_{2,k}\(y\),y\>_\xi+\frac{n-2}{2}\hat{u}_{2,k}\(y\)&=-\<\nabla U_0\(y\),y\>_\xi-\frac{n-2}{2}U_0\(y\)\nonumber\\
&\quad+\bigO\(\left|\(\hat{u}_{2,k}+U_0\)\(y\)\right|+\left|y\right|\left|\nabla\(\hat{u}_{2,k}+U_0\)\(y\)\right|\)
\end{align}
and
\begin{align}\label{Lem8Eq4}
c_n\mu_{2,k}^2\Scal_{\hat{g}_{2,k}}\(y\)\hat{u}_{2,k}\(y\)&=-\hat{w}_{2,k}\(y\)U_0\(y\)+\bigO\(\left|\(c_n\mu_{2,k}^2\Scal_{\hat{g}_{2,k}}-\hat{w}_{2,k}\)\(y\)\right|\right.\nonumber\\
&\quad\left.\times\left|\hat{u}_{2,k}\(y\)\right|+\left|\hat{w}_{2,k}\(y\)\right|\left|\(\hat{u}_{2,k}+U_0\)\(y\)\right|\)
\end{align}
uniformly with respect to $y\in \B_\xi\(0,\delta_0\varrho_{2,k}\)$ and $k>k_0$. It follows from \eqref{Lem8Eq2}, \eqref{Lem8Eq3} and \eqref{Lem8Eq4} that 
\begin{align}\label{Lem8Eq5}
&\(\<\nabla\hat{u}_{2,k}\(y\),y\>_\xi+\frac{n-2}{2}\hat{u}_{2,k}\(y\)\)\(\(\Delta_{\hat{g}_{2,k}}-\Delta_\xi\)\hat{u}_{2,k}+c_n\mu_{2,k}^2\Scal_{\hat{g}_{2,k}}\hat{u}_{2,k}\)\(y\)\nonumber\\
&=\frac{n^{\frac{n-2}{2}}\(n-2\)^{\frac{n}{2}}\big(1-\left|y\right|^2\big)}{2\big(1+\left|y\right|^2\big)^{n-1}}\hat{w}_{2,k}\(y\)\nonumber\allowdisplaybreaks\\
&\quad+\bigO \Big(\left|\hat{w}_{2,k}\(y\)\right|U_0\(y\) \(\left|\(\hat{u}_{2,k}+U_0\)\(y\)\right|+\left|y\right|\left|\nabla\(\hat{u}_{2,k}+U_0\)\(y\)\right|\)\nonumber\allowdisplaybreaks\\
& \quad +\Big( U_0\(y\)+\left|y\right|\left|\nabla U_0\(y\)\right| + \left|\(\hat{u}_{2,k}+U_0\)\(y\)\right| +\left|y\right|\left|\nabla\(\hat{u}_{2,k}+U_0\)\(y\)\right|\Big) \nonumber\allowdisplaybreaks\\
&\quad \times \Big( \left|\nabla\hat{g}_{2,k}\(y\)\right|\left|\nabla\(\hat{u}_{2,k}+U_0\)\(y\)\right|+ \mu_{2,k}^N\left|y\right|^{N-1}\left|\nabla U_0\(y\)\right|\nonumber \allowdisplaybreaks\\
&\quad + \left|\(\hat{g}_{2,k}-\xi\)\(y\)\right|\left|\nabla^2\(\hat{u}_{2,k}+U_0\)\(y\)\right|+\left|\(c_n\mu_{2,k}^2\Scal_{\hat{g}_{2,k}}-\hat{w}_{2,k}\)\(y\)\right| \left|\hat{u}_{2,k}\(y\)\right|\nonumber\allowdisplaybreaks\\
&\quad+\left|\hat{w}_{2,k}\(y\)\right|\left|\(\hat{u}_{2,k}+U_0\)\(y\)\right|\Big)\Big)
\end{align}
uniformly with respect to $y\in \B_\xi\(0,\delta_0\varrho_{2,k}\)$ and $k>k_0$. Moreover, \eqref{Lem4Eq1} and \eqref{Lem5Eq41} give
\begin{align}\label{Lem8Eq6} 
&\left|\hat{w}_{2,k}\(y\)\right|\( \left|\(\hat{u}_{2,k} + U_0\)\(y\)\right|  + \left|y\right| \left| \nabla \( \hat{u}_{2,k} + U_0\)\(y\) \right|\)\nonumber\\
&\quad= \bigO \(\left\{\sum_{\left|\alpha\right|=2}^{n-4}\frac{\left|h_{2,k,\alpha}\right|^2\mu_{2,k}^{2\left|\alpha\right|}}{\(1+\left|y\right|\)^{n-2\left|\alpha\right|}}+\frac{\mu_{2,k}^{n-1}}{1+\left|y\right|}\quad\text{if }n\ge6\right\}+\mu_{2,k}^2 \varrho_{2,k}^{2-n}\)
\end{align} 
and
\begin{equation}\label{Lem8Eq7} 
\left|y\right|\left|\nabla U_0\(y\)\right|+\left|\(\hat{u}_{2,k} + U_0\)\(y\)\right|  + \left|y\right| \left| \nabla \( \hat{u}_{2,k} + U_0\)\(y\) \right| = \bigO \( U_0\(y\) \)
\end{equation} 
uniformly with respect to $y\in \B_\xi\(0,\delta_0\varrho_{2,k}\)$ and $k>k_0$. By using \eqref{Lem5Eq8} together with the fact that $d_n \le n-4$ when $n \ge 6$, straightforward estimates  give
\begin{align}\label{Lem8Eq8} 
&\left|\(c_n\mu_{2,k}^2\Scal_{\hat{g}_{2,k}}-\hat{w}_{2,k}\)\(y\)\right|\left|\hat{u}_{2,k}\(y\)\right|\nonumber\\
&\quad= \bigO \( \sum_{\left|\alpha\right|=2}^{n-4}\frac{\left|h_{2,k,\alpha}\right|^2\mu_{2,k}^{2\left|\alpha\right|}}{\(1+\left|y\right|\)^{n-2\left|\alpha\right|}} +  \frac{\mu_{2,k}^2\( \mu_{2,k} \left|y\right|\)^{\max\(n-3,2\)}}{\(1+\left|y\right|\)^{n-2}}\)
\end{align}
uniformly with respect to $y\in \B_\xi\(0,\delta_0\varrho_{2,k}\)$ and $k>k_0$. By using \eqref{Lem3Eq4} and \eqref{Lem5Eq41}, we obtain
\begin{align}\label{Lem8Eq9} 
& \left|\nabla\hat{g}_{2,k}\(y\)\right|\left|\nabla\(\hat{u}_{2,k}+U_0\)\(y\)\right|+\left|\(\hat{g}_{2,k}-\xi\)\(y\)\right|\left|\nabla^2\(\hat{u}_{2,k}+U_0\)\(y\)\right|\nonumber\\
&\qquad+\mu_{2,k}^N\left|y\right|^{N-1}\left|\nabla U_0\(y\)\right|\nonumber\\
&\quad= \bigO \(\sum_{\left|\alpha\right|=2}^{n-4}\frac{\left|h_{2,k,\alpha}\right|^2\mu_{2,k}^{2\left|\alpha\right|}}{\(1+\left|y\right|\)^{n-2\left|\alpha\right|}} +  \frac{\mu_{2,k}^2\( \mu_{2,k}\left|y\right|\)^{\max\(n-3,2\)}}{\(1+\left|y\right|\)^{n-2}}+\frac{\(\mu_{2,k}\left|y\right|\)^N}{\(1+\left|y\right|\)^n} + \mu_{2,k}^2 \varrho_{2,k}^{2-n}\)  
\end{align}
uniformly with respect to $y\in \B_\xi\(0,\delta_0\varrho_{2,k}\)$ and $k>k_0$. By plugging \eqref{Lem8Eq6}, \eqref{Lem8Eq7}, \eqref{Lem8Eq8} and \eqref{Lem8Eq9} into \eqref{Lem8Eq5}, we obtain
\begin{align} \label{Lem8Eq10} 
&\(\<\nabla\hat{u}_{2,k}\(y\),y\>_\xi+\frac{n-2}{2}\hat{u}_{2,k}\(y\)\)\(\(\Delta_{\hat{g}_{2,k}}-\Delta_\xi\)\hat{u}_{2,k}+c_n\mu_{2,k}^2\Scal_{\hat{g}_{2,k}}\hat{u}_{2,k}\)\(y\)\nonumber\\
&\quad=\frac{n^{\frac{n-2}{2}}\(n-2\)^{\frac{n}{2}}\big(1-\left|y\right|^2\big)}{2\big(1+\left|y\right|^2\big)^{n-1}}\hat{w}_{2,k}\(y\)+\bigO\Bigg(\sum_{\left|\alpha\right|=2}^{n-4}\frac{\left|h_{2,k,\alpha}\right|^2\mu_{2,k}^{2\left|\alpha\right|}}{\(1+\left|y\right|\)^{2n-2\left|\alpha\right|-2}}\nonumber\\
&\qquad+\frac{\mu_{2,k}^2\(\mu_{2,k}\left|y\right|\)^{\max\(n-3,2\)}}{\(1+\left|y\right|\)^{2n-4}}+\frac{\mu_{2,k}^2\varrho_{2,k}^{2-n}}{\(1+\left|y\right|\)^{n-2}}+\frac{\(\mu_{2,k}\left|y\right|\)^N}{\(1+\left|y\right|\)^{2n-2}}\Bigg)
\end{align}
uniformly with respect to $y\in \B_\xi\(0,\delta\varrho_{2,k}\)$, $\delta\in\(0,\delta_0\)$ and $k>k_0$. We now integrate \eqref{Lem8Eq2} in $\B_\xi\(0,\delta\varrho_{2,k}\)$. By using \eqref{Lem4Eq1} and \eqref{Lem4Eq9}, we obtain
\begin{equation} \label{Lem8Eq11}
\int_{\B_\xi\(0,\delta\varrho_{2,k}\)}\frac{n^{\frac{n-2}{2}}\(n-2\)^{\frac{n}{2}}\big(1-\left|y\right|^2\big)}{2\big(1+\left|y\right|^2\big)^{n-1}}\hat{w}_{2,k}\(y\) \dy = 0.
\end{equation}
On the other hand, when $2 \le \left|\alpha\right| \le n-4$, straightforward estimates give 
\begin{equation} \label{Lem8Eq12}
\int_{\B_\xi\(0,\delta\varrho_{2,k}\)}\frac{\mu_{2,k}^{2\left|\alpha\right|}\dy}{\(1+\left|y\right|\)^{2n-2\left|\alpha\right|-2}}=\bigO\(\left\{\begin{aligned}&\mu_{2,k}^{2\left|\alpha\right|}&&\text{if }\left|\alpha\right| < d_n \\&  \mu_{2,k}^{n-2} \left| \ln  \mu_{2,k} \right| &&\text{if }\left|\alpha\right|=d_n\\& \delta^{n-2} \varrho_{i,k}^{2-n}&&\text{if }\left|\alpha\right|>d_n\end{aligned}\right\}\) 
\end{equation}
uniformly with respect to $\delta\in\(0,\delta_0\)$ and $k>k_0$. It follows from \eqref{Lem8Eq10}, \eqref{Lem8Eq11} and \eqref{Lem8Eq12} that
\begin{align}\label{Lem8Eq13}
&\int_{\B_\xi\(0,\delta\varrho_{2,k}\)}\hspace{-2pt}\(\<\nabla\hat{u}_{2,k},\cdot\>_\xi+\frac{n-2}{2}\hat{u}_{2,k}\)\(\(\Delta_{\hat{g}_{2,k}}-\Delta_\xi\)\hat{u}_{2,k}+c_n\mu_{2,k}^2\Scal_{\hat{g}_{2,k}}\hat{u}_{2,k}\)\dy\nonumber\\
&\quad=\bigO\Bigg(\sum_{\left|\alpha\right|=2}^{d_n}\left|h_{2,k,\alpha}\right|^2 \mu_{2,k}^{2\left|\alpha\right|} \left|\ln\mu_{2,k}\right|^{\vartheta\(2\left|\alpha\right|,n-2\)} +\delta\varrho_{2,k}^{2-n}\(\delta^{1-n}\(\delta\mu_{2,k}\varrho_{2,k}\)^{\max\(n-1,4\)}\right.\nonumber\\
&\qquad\left.+\delta\(\mu_{2,k}\varrho_{2,k}\)^2+\delta^{N-n+1}\(\mu_{2,k}\varrho_{2,k}\)^N \)\Bigg)
\end{align}
uniformly with respect to $\delta\in\(0,\delta_0\)$ and $k>k_0$. Finally, \eqref{Lem8Eq1} follows from \eqref{Lem5Eq14} and \eqref{Lem8Eq13}.
\endproof

We point out that in the proofs of Lemmas~\ref{Lem6}, \ref{Lem7} and \ref{Lem8}, we only use \eqref{Lem5Eq11} with $i=2$, namely for the most concentrated bubble. We recall that while this estimate holds true without any additional assumptions in the case where $i=2$, we need to show that \eqref{Lem5Eq10} holds true in order to use it in the case where $i=1$. This is done in the next section.

\section{Proofs of Theorems~\ref{Th1} and \ref{Th2}}\label{Sec4}

In this section, we apply the analysis of Section~\ref{Sec3} to prove Theorems~\ref{Th1} and \ref{Th2}. By using Lemmas~\ref{Lem6},~\ref{Lem7} and ~\ref{Lem8}, we can first complete the proof of Theorem~\ref{Th2}:

\proof[End of proof of Theorem~\ref{Th2}]
When $n\in\left\{3,4,5\right\}$, Lemma~\ref{Lem8} gives
\begin{align*}
&\int_{\B_\xi\(0,\delta\varrho_{2,k}\)}\hspace{-2pt}\(\<\nabla\hat{u}_{2,k},\cdot\>_\xi+\frac{n-2}{2}\hat{u}_{2,k}\)\(\(\Delta_{\hat{g}_{2,k}}-\Delta_\xi\)\hat{u}_{2,k}+c_n\mu_{2,k}^2\Scal_{\hat{g}_{2,k}}\hat{u}_{2,k}\)\dy\nonumber\\
&\quad=\bigO\big(\delta\varrho_{2,k}^{2-n}\big)
\end{align*} 
uniformly with respect to $\delta\in\(0,\delta_0\)$ and $k>k_0$, which, together with Lemma~\ref{Lem6}, yields an obvious contradiction with Lemma~\ref{Lem7} as $\delta \to 0$.
\endproof

In larger dimensions $n\in\left\{6,\dotsc,10\right\}$, the Pohozaev identity of Lemma~\ref{Lem6} alone is not enough to conclude and we need to perform a more refined analysis. As a first result, we obtain a priori estimates on $\(\varrho_{1,k}\)_k$ and $\(\varrho_{2,k}\)_k$ as well as a sharp asymptotic expansion of $\Lambda_1\(\S^n,\[g_k\]\)$ as $k \to \infty$. For each $k\in\N$, we define
$$\zeta_k\(x,\mu\):=\sum_{\left|\alpha\right|=2}^{d_n}\left|h_{k,x,\alpha}\right|^2\mu_{}^{2\left|\alpha\right|}\left|\ln\mu\right|^{\vartheta\(2\left|\alpha\right|,n-2\)}\quad\forall x\in\S^n,\mu>0,$$
where $\vartheta$ is as in Lemma~\ref{Lem8}. By using the asymptotic analysis performed in Lemmas~\ref{Lem1} to \ref{Lem8}, we obtain the following:

\begin{lemma}\label{Lem9}
Let $\(\mu_{1,k}\)_k$ and $\(\mu_{2,k}\)_k$ be as in Lemma~\ref{Lem1}, $d_n$ be as in Lemma~\ref{Lem3}, $\(\varrho_{2,k}\)_k$ and $\(h_{2,k,\alpha}\)_{k,\alpha}$ be as in \eqref{Lem5Eq1} and \eqref{Lem5Eq2}, respectively, and $k_0$ be as in Lemma~\ref{Lem5}. If $6\le n\le10$, then
\begin{equation}\label{Lem9Eq1}
\varrho_{1,k}^{2-n}+\varrho_{2,k}^{2-n}=\bigO\(\max_{\S^n}\(\zeta_k\(\cdot,\mu_{1,k}\)\)\)
\end{equation}
and
\begin{equation}\label{Lem9Eq2}
\Lambda_1\(\S^n,\[g_k\]\)=\Lambda_1\(\S^n,\[g_0\]\)+\bigO\(\max_{\S^n}\(\zeta_k\(\cdot,\mu_{1,k}\)\)\)
\end{equation}
uniformly with respect to $k>k_0$.
\end{lemma}

\proof[Proof of Lemma~\ref{Lem9}]
We begin with proving that \eqref{Lem9Eq1} holds true. It follows from Lemmas~\ref{Lem6}, \ref{Lem7} and \ref{Lem8} that
\begin{equation}\label{Lem9Eq3}
\varrho_{2,k}^{2-n}=\bigO\(\zeta_k\(x_{2,k},\mu_{2,k}\)\) = \smallo \(\left \{ \begin{aligned} 
& \mu_{2,k}^4\left|\ln\mu_{2,k}\right| && \text{if } n = 6 \\
& \mu_{2,k}^4  && \text{if } 7 \le n \le 10  \\
\end{aligned} \right\}\)\quad\text{as }k\to\infty.
\end{equation}
We assume by contradiction that, up to a subsequence,
\begin{equation}\label{Lem9Eq4}
d_k=\bigO\(\mu_{1,k}\)
\end{equation}
uniformly with respect to $k>k_0$. It follows from \eqref{Lem5Eq1} and \eqref{Lem9Eq4} that
$$\frac{\mu_{2,k}}{\mu_{1,k}} = \bigO \( \varrho_{2,k}^{-2}\),$$ 
which, together with \eqref{Lem9Eq3}, gives 
\begin{equation}\label{Lem9Eq5}
\frac{\mu_{2,k}}{\mu_{1,k}} =\smallo\( \left \{ \begin{aligned} 
& \mu_{2,k}^2\left|\ln\mu_{2,k}\right|^{\frac12} && \text{if } n = 6 \\
& \mu_{2,k}^{\frac{8}{n-2}}  && \text{if } 7 \le n \le 10  \\
\end{aligned} \right\}\) = \smallo\(\mu_{2,k}\)\quad\text{as }k\to\infty.
\end{equation}
Clearly, \eqref{Lem9Eq5} contradicts the fact that $\mu_{1,k}\to0$ as $k\to\infty$. This proves that \eqref{Lem9Eq1} holds true, and thus \eqref{Lem1Eq3} and \eqref{Lem1Eq4} hold true. Moreover, by using \eqref{Lem9Eq1} and \eqref{Lem9Eq3} together with \eqref{Lem5Eq1} and the facts that $\mu_{2,k}\le\mu_{1,k}$, the function $\mu \mapsto \mu^{\frac{n-2}{2}} \left|\ln \mu\right|$ is decreasing and $\frac{n-2}{2}\le4$ when $n\le10$, we obtain
\begin{align}\label{Lem9Eq6}
\varrho_{1,k}^{2-n}&=d_k^{2-n}\mu_{1,k}^{n-2}\nonumber\\
&\sim\varrho_{2,k}^{2-n}\mu_{1,k}^{\frac{n-2}{2}}\mu_{2,k}^{-\frac{n-2}{2}}\nonumber\\
&=\bigO\(\mu_{1,k}^{\frac{n-2}{2}}\sum_{\left|\alpha\right|=2}^{d_n}\left|h_{2,k,\alpha}\right|^2\mu_{2,k}^{2\left|\alpha\right|-\frac{n-2}{2}}\left|\ln\mu_{2,k}\right|^{\vartheta\(2\left|\alpha\right|,n-2\)}\)\nonumber\\
&=\bigO\(\zeta_k\(x_{2,k},\mu_{1,k}\)\)
\end{align}
uniformly with respect to $k>k_0$. Finally, \eqref{Lem9Eq1} follows from \eqref{Lem9Eq3} and \eqref{Lem9Eq6}. 

\smallskip
We now prove \eqref{Lem9Eq2} by using \eqref{IntroEq1} and \eqref{Sec3Eq5} and estimating the energy of $u_k$. We claim that 
\begin{equation}\label{Lem9Eq7}
\int_{\S^n}\left|u_k\right|^{2^*}\dv_{g_k}=2\Lambda_1\(\S^n,\[g_0\]\)^{\frac{n}{2}}+\bigO\(\max_{\S^n}\(\zeta_k\(\cdot,\mu_{1,k}\)\)\)
\end{equation}
uniformly with respect to $k>k_0$, which, together with \eqref{IntroEq1} and \eqref{Sec3Eq5}, implies \eqref{Lem9Eq2}. We prove this claim. For each $i\in\left\{1,2\right\}$ and $\delta\in\(0,\delta_0\)$, we let $\widetilde{g}_{i,k}:=g_{k,x_{i,k}}$ be given by Lemma~\ref{Lem2}. By using \eqref{Sec3Eq3} together with a rescaling argument and the conformal covariance of the conformal Laplacian, we obtain
\begin{align}\label{Lem9Eq8}
&\int_{\B_{\widetilde{g}_{i,k}}\(x_{i,k},\delta\mu_{i,k}\varrho_{i,k}\)}\left|u_k\right|^{2^*}\dv_{\widetilde{g}_{i,k}}\nonumber\\
&\quad=\int_{\B_\xi\(0,\delta\varrho_{i,k}\)}\(\Delta_{\hat{g}_{i,k}}\hat{u}_{i,k}+c_n\mu_{i,k}^2\Scal_{\hat{g}_{i,k}}\hat{u}_{i,k}\)\hat{u}_{i,k}\dv_{\hat{g}_{i,k}}.
\end{align}
By integrating by parts, we obtain
\begin{align}\label{Lem9Eq9}
&\int_{\B_\xi\(0,\delta\varrho_{i,k}\)}\hat{u}_{i,k}\Delta_{\hat{g}_{i,k}}\hat{u}_{i,k}\dv_{\hat{g}_{i,k}}\nonumber\\
&\quad=\int_{\B_\xi\(0,\delta\varrho_{i,k}\)}\Big(U_0\Delta_{\hat{g}_{i,k}}U_0+2\big(\hat{u}_{i,k}+\(-1\)^iU_0\big)\Delta_{\hat{g}_{i,k}}\hat{u}_{i,k}\nonumber\allowdisplaybreaks\\
&\qquad-\big|\nabla\big(\hat{u}_{i,k}+\(-1\)^iU_0\big)\big|_{\hat{g}_{i,k}}^2\Big)\dv_{\hat{g}_{i,k}}+\int_{\partial \B_\xi\(0,\delta\varrho_{i,k}\)}\big(\(-1\)^iU_0\partial_\nu\big(\hat{u}_{i,k}+\(-1\)^iU_0\big)\nonumber\\
&\qquad+\big(\hat{u}_{i,k}+\(-1\)^iU_0\big)\partial_\nu\hat{u}_{i,k}\big)\dsigma_{\hat{g}_{i,k}}.
\end{align}
We recall that, by \eqref{Lem9Eq6}, Lemma~\ref{Lem5} now applies to both $i=1$ and $i=2$. {By using Lemma~\ref{Lem5},} \eqref{Lem2Eq3}, \eqref{Lem5Eq11}, {\eqref{Lem5Eq41},} \eqref{Lem8Eq12} and \eqref{Lem9Eq9}, we obtain
\begin{align}\label{Lem9Eq10}
&\int_{\B_\xi\(0,\delta\varrho_{i,k}\)}\hat{u}_{i,k}\Delta_{\hat{g}_{i,k}}\hat{u}_{i,k}\dv_{\hat{g}_{i,k}}\nonumber\\
&\ \ =\int_{\B_\xi\(0,\delta\varrho_{i,k}\)}U_0\Delta_\xi U_0\dy+2\int_{\B_\xi\(0,\delta\varrho_{i,k}\)}\big(\hat{u}_{i,k}+\(-1\)^iU_0\big)\Delta_{\hat{g}_{i,k}}\hat{u}_{i,k}\big)\dv_{\hat{g}_{i,k}}\nonumber\allowdisplaybreaks\\
&\quad\ \ +\bigO\left(\int_{\B_\xi\(0,\delta\varrho_{i,k}\)}\left(\left\{\sum_{\left|\alpha\right|=2}^{n-4}\frac{\left|h_{i,k,\alpha}\right|^2\mu_{i,k}^{2\left|\alpha\right|}}{\(1+\left|y\right|\)^{2n-2\left|\alpha\right|-2}}+\frac{\mu_{i,k}^{2n-6}}{\(1+\left|y\right|\)^4}\quad\text{if }n\ge6\right\}\right.\right.
\nonumber\allowdisplaybreaks\\
&\quad\ \ +\frac{\varrho_{i,k}^{4-2n}}{\(1+\left|y\right|\)^2}+\frac{\(\mu_{i,k}\left|y\right|\)^N}{\(1+\left|y\right|\)^{2n-2}}\Bigg)\dy\Bigg) \nonumber\allowdisplaybreaks\\
&\quad\ \  +\bigO \(\int_{\partial \B_\xi\(0,\delta\varrho_{i,k}\)}\(\sum_{\left|\alpha\right|=2}^{n-4}\frac{\left|h_{i,k,\alpha}\right|\mu_{i,k}^{\left|\alpha\right|}}{\(1+\left|y\right|\)^{2n-\left|\alpha\right|-3}} +\frac{\varrho_{i,k}^{2-n}}{\(1+\left|y\right|\)^{n-1}}\)\dsigma\)\nonumber\allowdisplaybreaks\\
&\ \ =\int_{\R^n}U_0\Delta_\xi U_0\dy+2\int_{\B_\xi\(0,\delta\varrho_{i,k}\)}\big(\hat{u}_{i,k}+\(-1\)^iU_0\big)\Delta_{\hat{g}_{i,k}}\hat{u}_{i,k}\dv_{\hat{g}_{i,k}}\nonumber\\
&\quad\ \ +\bigO\(\zeta_k\(x_{i,k},\mu_{i,k}\)+\varrho_{i,k}^{2-n}\)
\end{align}
uniformly with respect to $k>k_0$. We now estimate the last two terms in the right-hand side of \eqref{Lem9Eq8}. By using \eqref{Lem2Eq3}, \eqref{Lem4Eq9}, \eqref{Lem5Eq8}, \eqref{Lem5Eq41} and \eqref{Lem8Eq12}, we obtain
\begin{align}\label{Lem9Eq11}
&c_n\mu_{i,k}^2\int_{\B_\xi\(0,\delta\varrho_{i,k}\)}\Scal_{\hat{g}_{i,k}}\hat{u}_{i,k}^2\dv_{\hat{g}_{i,k}}\nonumber\\
&\quad=c_n\mu_{i,k}^2\int_{\B_\xi\(0,\delta\varrho_{i,k}\)}\Scal_{\hat{g}_{i,k}}\Big(U_0^2+2\hat{u}_{i,k}\big(\hat{u}_{i,k}+\(-1\)^iU_0\big)\nonumber\allowdisplaybreaks\\
&\qquad-\big(\hat{u}_{i,k}+\(-1\)^iU_0\big)^2\Big)\dv_{\hat{g}_{i,k}}\nonumber\allowdisplaybreaks\\
&\quad=\int_{\B_\xi\(0,\delta\varrho_{i,k}\)}\bigg(\hat{w}_{i,k}U_0^2+2c_n\mu_{i,k}^2\Scal_{\hat{g}_{i,k}}\hat{u}_{i,k}\big(\hat{u}_{i,k}+\(-1\)^iU_0\big)\nonumber\allowdisplaybreaks\\
&\qquad+\bigO\Big(\left|c_n\mu_{i,k}^2\Scal_{\hat{g}_{i,k}}-\hat{w}_{i,k}\right|U_0^2+\mu_{i,k}^2\left|\Scal_{\hat{g}_{i,k}}\right|\big| \hat{u}_{i,k} + \(-1\)^i U_0 \big|^2
\Big)\bigg)\dv_{\hat{g}_{i,k}}\nonumber\allowdisplaybreaks\\
&\quad=2c_n\mu_{i,k}^2\int_{\B_\xi\(0,\delta\varrho_{i,k}\)}\Scal_{\hat{g}_{i,k}}\hat{u}_{i,k}\big(\hat{u}_{i,k}+\(-1\)^iU_0\big)\dv_{\hat{g}_{i,k}}\nonumber\allowdisplaybreaks\\
&\qquad+\bigO\(\int_{\B_\xi\(0,\delta\varrho_{i,k}\)}\(\sum_{\left|\alpha\right|=2}^{n-4}\frac{\left|h_{i,k,\alpha}\right|^2\mu_{i,k}^{2\left|\alpha\right|}}{\(1+\left|y\right|\)^{2n-2\left|\alpha\right|-2}}+\frac{\mu_{i,k}^{n-1}}{\(1+\left|y\right|\)^{n-1}}
+\mu_{i,k}^4\varrho_{i,k}^{4-2n}\left|y\right|^2 \right.\right.\nonumber\allowdisplaybreaks\\
&\qquad\left.\left.+\frac{\mu_{i,k}^{N+4}\left|y\right|^{N+2}}{\(1+\left|y\right|\)^{2n-4}}\)\dy\)\nonumber\allowdisplaybreaks\\
&\quad=2c_n\mu_{i,k}^2\int_{\B_\xi\(0,\delta\varrho_{i,k}\)}\Scal_{\hat{g}_{i,k}}\hat{u}_{i,k}\big(\hat{u}_{i,k}+\(-1\)^iU_0\big)\dv_{\hat{g}_{i,k}} \nonumber\\
&\qquad +\bigO\big(\zeta_k\(x_{i,k},\mu_{i,k}\)+\varrho_{i,k}^{2-n}\big)
\end{align}
uniformly with respect to $k>k_0$. By putting together \eqref{Lem9Eq8}, \eqref{Lem9Eq10} and \eqref{Lem9Eq11} and using the equations \eqref{Lem5Eq21} and $\Delta_\xi U_0=U_0^{2^*-1}$, we obtain
{\begin{align}\label{Lem9Eq13}
&\int_{\B_{\widetilde{g}_{i,k}}\(x_{i,k},\delta\mu_{i,k}\varrho_{i,k}\)}\left|u_k\right|^{2^*}\dv_{\widetilde{g}_{i,k}}\nonumber\\
&\quad=\int_{\R^n}U_0^{2^*}\dy+2 \int_{\B_\xi\(0,\delta\varrho_{i,k}\)}\left|\hat{u}_{i,k}\right|^{2^*-2}\hat{u}_{i,k} \big( \hat{u}_{i,k} + (-1)^i U_0 \big) \dv_{\hat{g}_{i,k}} 
\nonumber\\
&\qquad+\bigO\big(\zeta_k\(x_{i,k},\mu_{i,k}\) +\varrho_{i,k}^{2-n}\big)
\end{align}
uniformly with respect to $k>k_0$. We now estimate the integral in the right-hand side of \eqref{Lem9Eq13}. Direct computations using \eqref{Lem2Eq3}, \eqref{Lem5Eq9} and \eqref{Lem5Eq11} give
\begin{align}\label{Lem9Eq13:bis}
&\int_{\B_\xi\(0,\delta\varrho_{i,k}\)}\left|\hat{u}_{i,k}\right|^{2^*-2}\hat{u}_{i,k} \big( \hat{u}_{i,k} + (-1)^i U_0 \big) \dv_{\hat{g}_{i,k}}\nonumber\\
 &\ =  \int_{\B_\xi\(0,\delta\varrho_{i,k}\)}(-1)^i U_0^{2^*-1} \hat{v}_{i,k} \dv_{\hat{g}_{i,k}} + \bigO \left(\int_{\B_\xi\(0,\delta\varrho_{i,k}\)} U_0^{2^*-2} \big| \hat{u}_{i,k} + (-1)^i U_0 \big|^2 \dv_{\hat{g}_{i,k}}  \right.\nonumber\\
&\ \quad \left.+\int_{\B_\xi\(0,\delta\varrho_{i,k}\)} U_0^{2^*-1} \big| \hat{u}_{i,k} + (-1)^i \big(U_0 - \hat{v}_{i,k} \big) \big| \dv_{\hat{g}_{i,k}}  \right)\allowdisplaybreaks\nonumber\\
 &\ =  \int_{\R^n}(-1)^i U_0^{2^*-1} \hat{v}_{i,k} \dy  +\bigO\big(\zeta_k\(x_{i,k},\mu_{i,k}\)+\varrho_{i,k}^{2-n}\big)\nonumber\\
 &\ = \bigO\big(\zeta_k\(x_{i,k},\mu_{i,k}\)+\varrho_{i,k}^{2-n}\big),
\end{align}
uniformly with respect to $k>k_0$, where the last equality follows from  \eqref{Lem4Eq15:bis} together with the definition of $\hat{v}_{i,k}$ in \eqref{Lem5Eq6}. Plugging \eqref{Lem9Eq13:bis} into \eqref{Lem9Eq13} then gives
\begin{equation}\label{Lem9Eq13:ter}
\int_{\B_{\widetilde{g}_{i,k}}\(x_{i,k},\delta\mu_{i,k}\varrho_{i,k}\)}\left|u_k\right|^{2^*}\dv_{\widetilde{g}_{i,k}}=\int_{\R^n}U_0^{2^*}\dy+\bigO\big(\zeta_k\(x_{i,k},\mu_{i,k}\)+\varrho_{i,k}^{2-n}\big)
\end{equation}
uniformly with respect to $k>k_0$.} We recall that
\begin{equation}\label{Lem9Eq14}
\int_{\R^n}U_0^{2^*}\dy=\Lambda_1\(\S^n,\[g_0\]\)^{\frac{n}{2}}.
\end{equation}
Moreover, by using \eqref{Lem5Eq10} together with the definition of $\varrho_{1,k}$ and the fact that $\mu_{2,k}\le\mu_{1,k}$, we obtain that there exists $\delta_1\in\(0,\delta_0\)$ and $k_1>k_0$ such that, for each $\delta\in\(0,\delta_1\)$ and $k>k_1$,
\begin{equation}\label{Lem9Eq15}
\B_{\widetilde{g}_{1,k}}\(x_{1,k},\delta\mu_{1,k}\varrho_{1,k}\)\cap\B_{\widetilde{g}_{2,k}}\(x_{2,k},\delta\mu_{2,k}\varrho_{2,k}\)=\emptyset.
\end{equation}
On the other hand, by using similar estimates as in the beginning of the proof of Lemma~\ref{Lem5}, we obtain
\begin{align}\label{Lem9Eq16}
&\int_{\displaystyle\S^n\setminus\bigcup_{i=1}^2\B_{\widetilde{g}_{i,k}}\(x_{i,k},\delta\mu_{i,k}\varrho_{i,k}\)}\left|u_k\right|^{2^*}\dv_{\widetilde{g}_{i,k}}\nonumber\\
&\quad=\bigO\(\sum_{i=1}^2\int_{\displaystyle\S^n\setminus \B_{\widetilde{g}_{i,k}}\(x_{i,k},\delta\mu_{i,k}\varrho_{i,k}\)}\(B_{3-i,k}\)^{2^*}\dv_{\widetilde{g}_{i,k}}\)\nonumber\\
&\quad=\bigO\big(\varrho_{1,k}^{-n}+\varrho_{2,k}^{-n}\big)
\end{align}
uniformly with respect to $k>k_1$. By using \eqref{Lem9Eq1}, {\eqref{Lem9Eq13:ter}}, \eqref{Lem9Eq14}, \eqref{Lem9Eq15} and \eqref{Lem9Eq16} together with the fact that $\mu_{2,k}\le\mu_{1,k}$, we obtain \eqref{Lem9Eq7}, which completes the proof of Lemma~\ref{Lem9}.
\endproof

We are now in position to conclude the proof of Theorem~\ref{Th1}. The idea to reach a contradiction consists in directly estimating $\Lambda_1\(\S^n,\[g_k\]\)$ with the help of suitable test-functions. These test-functions are modeled on the first-order expansion of the functions  $\(u_k\)_k$ as in \eqref{Lem5Eq11}, centered at maximum points of the functions $\(\zeta_k\)_k$ in $M$ and more concentrated than the functions $\(B_{1,k}\)_k$. We prove that these test-functions provide better competitors for $\Lambda_1\(\S^n,\[g_k\]\)$ and yield a contradiction with \eqref{Lem9Eq2}. As we mentioned in Section~\ref{Sec2}, our contradiction argument comes from the very definition of $\Lambda_2\(\S^n,\[g_k\]\)$ and its minimality. We cannot use a local sign restriction argument as in \cite{KMS} here since all local masses vanish due to the fact that $g_k \to g_0$ in $C^m\(\S^n\)$ for all $m\in\N$.

\proof[End of proof of Theorem~\ref{Th1}]
The form of the leading term of the test-functions we use is inspired from the historic work of Schoen~\cite{Sc1}. We add a lower-order term inspired from the work of Li and Zhang~\cites{LZha1,LZha2} and Khuri, Marques and Schoen~\cite{KMS} (see also the early work of Hebey and Vaugon~\cite{HV} using this idea of adding a small correction term in the test-functions). For each $k>k_0$, we let $x_k\in \S^n$ and $\mu_k\in\(0,\infty\)$ be such that 
\begin{equation} \label{Th2Eq1}
\zeta_k\(x_k,\mu_{1,k}\)=\max_{\S^n}\(\zeta_k\(\cdot,\mu_{1,k}\)\) \quad\text{and}\quad  \mu_k=\lambda\mu_{1,k},
\end{equation}
where $\lambda\in\(1,\infty\)$ is some fixed number to be chosen large later on. We let $\widetilde{g}_k:=g_{k,x_k}$ and $\exp_k:=\exp_{k,x_k}$ be as in Lemma~\ref{Lem2}, $H_k:=H_{k,x_k}$ and $h_{k,\alpha}:=h_{k,x_k,\alpha}$ be as in Lemma~\ref{Lem3}, $w_{k,d}:=w_{k,x_k,d}$ be as in \eqref{Lem4Eq1} and $v_{k,d}:=v_{k,x_k,d}$ be as in Lemma~\ref{Lem4}. Up to a subsequence, we may further assume that $x_k\to x_0\in\S^n$ as $k \to \infty$. We then define $\widetilde{g}_0:=g_{0,x_0}$. We let $\widetilde{G}_0$ be the Green's function of $L_{\widetilde{g}_0}$ in $\S^n$. In particular, for each $x \in \S^n$, $\widetilde{G}_0\(x, \cdot\) $ is a positive function in $\S^n  \setminus \left\{x\right\}$ satisfying the equation 
\begin{equation}\label{Th2Eq2}
L_{\widetilde{g}_0} \widetilde{G}_0\(x, \cdot\) = 0\quad\text{in }\S^n \setminus \left\{x\right\}.
\end{equation}
We let $\delta\in\(0,\varepsilon_0/2\)$,  where $\varepsilon_0$ is as in Lemma~\ref{Lem2}. We let $\eta$ be a smooth function on $\[0,\infty\)$ such that $\eta=1$ on $\[0,1\]$ and $\eta=0$ on $\[2, \infty\)$. We define the functions
\begin{align*}
B_k&:=\(\frac{\sqrt{n\(n-2\)}\mu_k}{\mu_k^2+\dist_{\widetilde{g}_k}\(x_k,\cdot\)^2}\)^{\frac{n-2}{2}},\allowdisplaybreaks\\
\hat{w}_k&:=c_n\sum_{d=4}^{n-4}\mu_k^dw_{k,d},\allowdisplaybreaks\\
\hat{v}_k&:=c_n\sum_{{d=4}}^{n-4}\mu_k^{{d}}v_{k,d},\allowdisplaybreaks\\
v_k&:=\mu_k^{-\frac{n-2}{2}}\hat{v}_k\circ\(\mu_k^{-1}\exp_k^{-1}\),\allowdisplaybreaks\\
\check{w}_k&:=  c_n\mu_k^2\sum_{a,b,c=1}^n\bigg(\partial_{y_b}\(\(H_k\)_{ab}\partial_{y_c}\(H_k\)_{ac}\)\\
&\quad-\frac{1}{2}\partial_{y_b}\(H_k\)_{ab}\partial_{y_c}\(H_k\)_{ac}+\frac{1}{4}\(\partial_{y_c}\(H_k\)_{ab}\)^2\bigg)\(\mu_k\cdot\),\allowdisplaybreaks\\
\Gamma_k&:=n^{\frac{n-2}{4}} \(n-2\)^{\frac{n+2}{4}} \omega_{n-1} \mu_k^{\frac{n-2}{2}}\widetilde{G}_0\(x_k, \cdot\)
\end{align*}
and
$$z_k:=\eta\(\delta^{-1}\dist_{\widetilde{g}_k}\(x_k,\cdot\)\) \(B_k - v_k\)+ \(1 - \eta\(\delta^{-1}\dist_{\widetilde{g}_k}\(x_k,\cdot\)\)\)\Gamma_k,$$
where $\omega_{n-1}$ is the volume of the round $\(n-1\)$-sphere. We also define the metric
$$\hat{g}_k:=\exp_k^*\widetilde{g}_k\(\mu_k\cdot\).$$
We recall that, by Lemma~\ref{Lem4}, $\hat{v}_k$ satisfies
\begin{equation}\label{Th2Eq3}
\Delta_{\xi} \hat{v}_k = \(2^*-1\) U_0^{2^*-2}\hat{v}_k + U_0 \hat{w}_k \quad \text{ in } \R^n.
\end{equation}
We claim that there exists $C_n >0$ depending only on $n$ such that
\begin{equation}\label{Th2Eq4}
\Lambda_1\(\S^n,\[g_0\]\)-\Lambda_1\(\S^n,\[g_k\]\)\ge C_n\zeta_k\(x_k,\mu_k\)+\smallo\(\mu_k^{n-2}\)\quad\text{as }k\to\infty.
\end{equation}
Before proving \eqref{Th2Eq4}, we first show that this estimate yields a contradiction with \eqref{Lem9Eq2} when $\lambda$ is chosen large enough, thus completing the proof of Theorem~\ref{Th1}. Since $\mu_{1,k} = \bigO \big( \varrho_{1,k}^{-1}\big)$ by \eqref{Lem5Eq13}, it follows from \eqref{Lem9Eq1} and \eqref{Th2Eq1} that 
\begin{equation}\label{Th2Eq5}
\mu_k^{n-2} = \bigO \(\zeta_k\(x_k,\mu_{1,k}\)\)
\end{equation}
uniformly with respect to $k>k_0$, where the constant in $\bigO\(\cdot\)$ depends on $\lambda$, but this is not a problem since this term is multiplied by $\smallo\(1\)$ in \eqref{Th2Eq4}. Moreover, since $\lambda>1$, by definition of $\zeta_k$, it is easy to see that 
\begin{equation}\label{Th2Eq6}
\zeta_k\(x_k,\mu_k\) \ge  \lambda^2 \zeta_k\(x_k,\mu_{1,k}\).
\end{equation}
It follows from \eqref{Th2Eq4}, \eqref{Th2Eq5} and \eqref{Th2Eq6} that
$$ \Lambda_1\(\S^n,\[g_0\]\)-\Lambda_1\(\S^n,\[g_k\]\)\ge \( \lambda^2C_n + \smallo\(1\) \)\zeta_k\(x_k,\mu_{1,k}\)\quad\text{as }k\to\infty,$$ 
which contradicts  \eqref{Lem9Eq2} when $\lambda$ is chosen large enough.

\smallskip
We now prove \eqref{Th2Eq4}. Since $\widetilde{g}_k\in\[g_k\]$ and $z_k\in C^\infty\(\S^n\)$, we obtain
\begin{equation}\label{Th2Eq7}
\Lambda_1\(\S^n,\[g_k\]\)\le\frac{\displaystyle\int_{\S^n}\(\left|\nabla z_k\right|^2_{\widetilde{g}_k}+c_n\Scal_{\widetilde{g}_k}z_k^2\)\dv_{\widetilde{g}_k}}{\displaystyle\(\int_{\S^n}\left| z_k\right|^{2^*}\dv_{\widetilde{g}_k}\)^{\frac{n-2}{n}}}.
\end{equation}
By definition of $z_k$, it is easy to see that 
\begin{equation}\label{Th2Eq8}
\int_{\S^n}\left|z_k\right|^{2^*}\dv_{\widetilde{g}_k} = \int_{\B_{\widetilde{g}_k}\(x_k,\delta\)} \left| B_k - v_k \right|^{2^*} \dv_{\widetilde{g}_k}+\smallo\(\mu_k^{n-2}\)\quad\text{as }k\to\infty.
\end{equation}
By using \eqref{Lem2Eq3}, \eqref{Lem4Eq6} and \eqref{Th2Eq8} together with similar estimates as in \eqref{Lem9Eq8}--\eqref{Lem9Eq16}, we obtain
\begin{align}\label{Th2Eq9}
\int_{\S^n}\left|z_k\right|^{2^*}\dv_{\widetilde{g}_k}&=\int_{\B_\xi\(0,\delta/\mu_k\)}\left|U_0-\hat{v}_k\right|^{2^*}\dv_{\hat{g}_k} +\smallo\( \mu_k^{n-2}\) \nonumber\\ 
&=\int_{\B_\xi\(0,\delta/\mu_k\)}\bigg(U_0^{2^*}-2^*U_0^{2^*-1}\hat{v}_k+\frac{2^*\(2^*-1\)}{2}U_0^{2^*-2}\hat{v}_k^2\nonumber\allowdisplaybreaks\\
&\quad+\bigO\(U_0^{2^*-3}\left|\hat{v}_k\right|^3\)\bigg)\dv_{\hat{g}_k}+\smallo\( \mu_k^{n-2}\)\nonumber\allowdisplaybreaks\\
&=\int_{{\B_\xi\(0,\delta/\mu_k\)}}\(U_0^{2^*}-2^*U_0^{2^*-1}\hat{v}_k+\frac{2^*\(2^*-1\)}{2}U_0^{2^*-2}\hat{v}_k^2\)\dy\nonumber\allowdisplaybreaks\\
&\quad+\bigO\(\int_{\B_\xi\(0,\delta/\mu_k\)}\sum_{\left|\alpha\right|=4}^{n-4}\frac{\left|h_{k,\alpha}\right|^3\mu_k^{3\left|\alpha\right|}}{\(1+\left|y\right|\)^{2n-3\left|\alpha\right|}}\dy\)+\smallo\( \mu_k^{n-2}\)\nonumber\allowdisplaybreaks\\
&=\Lambda_1\(\S^n,\[g_0\]\)^{\frac{n}{2}}-\frac{2^*}{2}\int_{{\B_\xi\(0,\delta/\mu_k\)}}\(2U_0^{2^*-1}\hat{v}_k-\(2^*-1\)U_0^{2^*-2}\hat{v}_k^2\)\dy\nonumber\\
&\quad
+\smallo\(\zeta_k\(x_k,\mu_k\)+\mu_k^{n-2}\)\quad\text{as }k\to\infty.
\end{align}
We now estimate the numerator in \eqref{Th2Eq7}. By using \eqref{Th2Eq2} together with the definition of $z_k$ and the facts that $x_k\to x_0$ in $\S^n$ and $\widetilde{g}_k \to \widetilde{g}_0$ in $C^m\(\S^n\)$ as $k\to\infty$ for all $m\in\N$, we obtain
\begin{align} \label{Th2Eq10}
&  \int_{\S^n \setminus \B_{\widetilde{g}_k}\(x_k, \delta\)} z_k L_{\widetilde{g}_k} z_k \dv_{\widetilde{g}_k}\nonumber\\
  &\quad=\int_{\B_{\widetilde{g}_k}\(x_k,2\delta\) \setminus \B_{\widetilde{g}_k}\(x_k, \delta\)} z_k L_{\widetilde{g}_k}\(z_k-\Gamma_k\)\dv_{\widetilde{g}_k}\nonumber\\
  &\qquad+n^{\frac{n-2}{4}}\(n-2\)^{\frac{n+2}{4}}\omega_{n-1}\mu_k^{\frac{n-2}{2}} \(  \int_{\S^n \setminus \B_{\widetilde{g}_0}\(x_0, 2 \delta\)}z_kL_{\widetilde{g}_0} \widetilde{G}_0\(x_0, \cdot\) \dv_{\widetilde{g}_k}+ \smallo\(1\) \)\nonumber\\
 &\quad =\int_{\B_{\widetilde{g}_k}\(x_k,2\delta\) \setminus \B_{\widetilde{g}_k}\(x_k, \delta\)} z_k L_{\widetilde{g}_k}\(z_k-\Gamma_k\)\dv_{\widetilde{g}_k}+ \smallo \(\mu_k^{n-2} \)\quad\text{as }k\to\infty.
\end{align}
On the other hand, since $\widetilde{g}_0$ is flat in $\B_{\widetilde{g}_0}\(x_0,\varepsilon_0\)$ by \eqref{Lem2Eq4}, we obtain (see for instance \cite{LP}) that there exists a function $R\in C^2\big(\overline{\B_{ \widetilde{g}_0}\(x_0,\varepsilon_0\)}\big)$ such that $R\(x_0\)=0$ and
\begin{equation} \label{Th2Eq11}
\widetilde{G}_0\(x_0,\cdot\)=\frac{1}{\(n-2\)\omega_{n-1}}\(\dist_{\widetilde{g}_0}\(x_0,\cdot\)^{2-n} + R\),
\end{equation}
which gives
\begin{align} \label{Th2Eq12}
B_k-\Gamma_k&=\(\sqrt{n\(n-2\)}\mu_k\)^{\frac{n-2}{2}}\(\big(\mu_k^2+\dist_{\widetilde{g}_k}\(x_k,\cdot\)^2\big)^{\frac{2-n}{2}}-\dist_{\widetilde{g}_0}\(x_0,\cdot\)^{2-n}-R\right.\nonumber\\
&\quad\left.+\(n-2\)\omega_{n-1}\(\widetilde{G}_0\(x_0,\cdot\)-\widetilde{G}_0\(x_k,\cdot\)\)\).
\end{align}
Since $R\(0\)=0$, $x_k\to0$, $\mu_k\to0$ and $\widetilde{g}_k\to\widetilde{g}_0$ in $C^m\(\S^n\)$ as $k\to\infty$ for all $m\in\N$, it follows from \eqref{Th2Eq12} that
\begin{align} \label{Th2Eq13}
\sum_{j=0}^2\delta^j\left|\nabla^j \(B_k-\Gamma_k\)\(y\)\right|&=\bigO\(\mu_k^{\frac{n-2}{2}}\sum_{j=0}^2\delta^j\left|\nabla^j R\(y\)\right|\)+\smallo\(\mu_k^{\frac{n-2}{2}}\)\nonumber\\
&=\bigO\(\delta\mu_k^{\frac{n-2}{2}}\)+\smallo\(\mu_k^{\frac{n-2}{2}}\)\quad\text{as }k\to\infty
\end{align}
uniformly with respect to $y \in \B_{ \widetilde{g}_k}\(x_k,2\delta\)\setminus\B_{ \widetilde{g}_k}\(x_k,\delta\)$, where the constant in the term $\bigO\(\delta\mu_k^{n-2}\)$ does not depend on $k$ nor on $\delta\in\(0,\varepsilon/2\)$ (the same holds true in the next estimates). Moreover, an easy consequence of \eqref{Lem5Eq9} is that
\begin{equation} \label{Th2Eq14}
\sum_{j=0}^2\delta^j\left|\nabla^j v_{k}\(y\)\right|=\smallo\(\mu_k^{\frac{n-2}{2}}\)\quad\text{as }k\to\infty
\end{equation}
uniformly with respect to $y \in \B_{ \widetilde{g}_k}\(x_k,2\delta\)\setminus\B_{ \widetilde{g}_k}\(x_k,\delta\)$. By using \eqref{Th2Eq13} and \eqref{Th2Eq14} together with the definition of $z_k$ and the facts that $R\(x_0\)=0$, $x_k\to x_0$ in $\S^n$ and $\widetilde{g}_k \to \widetilde{g}_0$ in $C^m\(\S^n\)$ as $k\to\infty$ for all $m\in\N$, we obtain
\begin{equation} \label{Th2Eq15}
\sum_{j=0}^2\delta^j\left|\nabla^j\(z_{k}-\Gamma_k\)\(y\)\right| =\bigO\(\delta\mu_k^{\frac{n-2}{2}}\)+\smallo\(\mu_k^{\frac{n-2}{2}}\)\quad\text{as }k\to\infty
\end{equation}
uniformly with respect to $y \in \B_{ \widetilde{g}_k}\(x_k,2\delta\)\setminus\B_{ \widetilde{g}_k}\(x_k,\delta\)$. By using \eqref{Th2Eq2} and \eqref{Th2Eq15}, we then obtain
\begin{equation} \label{Th2Eq16}
\int_{\B_{\widetilde{g}_k}\(x_k,2\delta\) \setminus \B_{\widetilde{g}_k}\(x_k, \delta\)} z_k L_{\widetilde{g}_k}\(z_k-\Gamma_k\)\dv_{\widetilde{g}_k}= \bigO\( \delta \mu_k^{n-2}\)+\smallo\(\mu_k^{n-2}\)\quad\text{as }k\to\infty.
\end{equation}
It follows from \eqref{Th2Eq10} and \eqref{Th2Eq16} together with an integration by parts and the definition of $z_k$ that
\begin{align} \label{Th2Eq17}
&\int_{\S^n}\(\left|\nabla z_k\right|^2_{\widetilde{g}_k}+c_n\Scal_{\widetilde{g}_k}z_k^2\)\dv_{\widetilde{g}_k}\nonumber\\
&\quad = \int_{\B_{\widetilde{g}_k}\(x_k,\delta\)}\(B_kL_{\widetilde{g}_k}B_k-2v_kL_{\widetilde{g}_k}B_k+v_kL_{\widetilde{g}_k}v_k\)\dv_{\widetilde{g}_k}+ \bigO\( \delta \mu_k^{n-2}\)\nonumber\\
  &\qquad + \smallo\(\mu_k^{n-2}\)\quad\text{as }k\to\infty.
\end{align}
By using \eqref{Lem2Eq3}, \eqref{Lem5Eq7}, \eqref{Lem5Eq8}, \eqref{Lem5Eq9} and \eqref{Th2Eq3} together with similar estimates as in \eqref{Lem9Eq9}--\eqref{Lem9Eq11}, we obtain
\begin{align}\label{Th2Eq18}
&\int_{\B_{\widetilde{g}_k}\(x_k,\delta\)}\(B_kL_{\widetilde{g}_k}B_k-2v_kL_{\widetilde{g}_k}B_k+v_kL_{\widetilde{g}_k}v_k\)\dv_{\widetilde{g}_k}\nonumber\\
&\ \ =\int_{\B_\xi\(0,\delta/\mu_k\)}\bigg(U_0\(\Delta_{\hat{g}_k}+c_n\mu_k^2\Scal_{\hat{g}_k}\)U_0-2\hat{v}_k\(\Delta_{\hat{g}_k}+c_n\mu_k^2\Scal_{\hat{g}_k}\)U_0\nonumber\allowdisplaybreaks\\
&\quad\ \ +\hat{v}_k\(\Delta_{\hat{g}_k}+c_n\mu_k^2\Scal_{\hat{g}_k}\)\hat{v}_k\bigg)\dv_{\hat{g}_k}\nonumber\allowdisplaybreaks\\
&\ \ =\int_{\B_\xi\(0,\delta/\mu_k\)}\bigg(U_0\(\Delta_\xi+\hat{w}_k-\check{w}_k\)U_0-2\hat{v}_k\(\Delta_\xi+\hat{w}_k\)U_0+\hat{v}_k\Delta_{\xi}\hat{v}_k\nonumber\allowdisplaybreaks\\
&\quad\ \ +\bigO\big(\left|c_n\mu_k^2\Scal_{\hat{g}_k}- \hat{w}_k+\check{w}_k\right|U_0^2+\left|c_n\mu_k^2\Scal_{\hat{g}_k}-\hat{w}_k\right|U_0\left|\hat{v}_k\right|+\mu_k^2\left|\Scal_{\hat{g}_k}\right|\hat{v}_k^2\nonumber\allowdisplaybreaks\\
&\quad\ \ +\(U_0+\left|\hat{v}_k\right|\)\left|\(\Delta_{\hat{g}_k}-\Delta_\xi\)U_0\right|+\left|\hat{v}_k\right|\(\left|\nabla\hat{g}_k\right|\left|\nabla\hat{v}_k\right|+\left|\hat{g}_k-\xi\right|\left|\nabla^2\hat{v}_k\right|\)\big) \bigg)\dv_{\hat{g}_k}\nonumber \allowdisplaybreaks\\
&\ \ =\Lambda_1\(\S^n,\[g_0\]\)^{\frac{n}{2}}-\int_{{\B_\xi\(0,\delta/\mu_k\)}}\(2U_0^{2^*-1}\hat{v}_k-\(2^*-1\)U_0^{2^*-2}\hat{v}_k^2\)\dy\nonumber\allowdisplaybreaks\\
&\quad\ \ -\int_{\B_\xi\(0,\delta/\mu_k\)}U_0\(U_0\check{w}_k+\hat{w}_k\hat{v}_k\)\dy+\bigO\(\sum_{\left|\alpha\right|=2}^{d_n-1}\left|h_{k,\alpha}\right|^2\mu_k^{2\left|\alpha\right|+2}\left|\ln\mu_k\right|^{\vartheta\(2\left|\alpha\right|,n-4\)}\right.\nonumber\allowdisplaybreaks\\
&\quad\ \ \left.+\sum_{\left|\alpha\right|=2}^{\[\(n-2\)/3\]}\left|h_{k,\alpha}\right|^3\mu_k^{3\left|\alpha\right|}\left|\ln\mu_k\right|^{\vartheta\(3\left|\alpha\right|,n-2\)}\)+ \smallo\(\mu_k^{n-2}\)\nonumber \allowdisplaybreaks\\
&\ \ =\Lambda_1\(\S^n,\[g_0\]\)^{\frac{n}{2}}-\int_{{\B_\xi\(0,\delta/\mu_k\)}}\(2U_0^{2^*-1}\hat{v}_k-\(2^*-1\)U_0^{2^*-2}\hat{v}_k^2\) \dy\nonumber\\
&\quad\ \ -\int_{\B_\xi\(0,\delta/\mu_k\)}U_0\(U_0\check{w}_k +\hat{w}_k\hat{v}_k\)\dy +\smallo\( \zeta_k\(x_k,\mu_k\)+\mu_k^{n-2}\)  
\text{ as }k\to\infty. 
\end{align}
It follows from \eqref{Th2Eq7}, \eqref{Th2Eq9},  \eqref{Th2Eq17} and \eqref{Th2Eq18} that
\begin{align}\label{Th2Eq19}
&\Lambda_1\(\S^n,\[g_0\]\)-\Lambda_1\(\S^n,\[g_k\]\)\nonumber\\
&=\Lambda_1\(\S^n,\[g_0\]\)^{-\frac{n-2}{2}}\int_{\B_\xi\(0,\delta/\mu_k\)}U_0\(U_0\check{w}_k+\hat{w}_k\hat{v}_k\)\dy\nonumber\allowdisplaybreaks\\
&\quad+\bigO\(\(\int_{{\B_\xi\(0,\delta/\mu_k\)}}\(2U_0^{2^*-1}\hat{v}_k-\(2^*-1\)U_0^{2^*-2}\hat{v}_k^2\) \dy\)^2+\delta \mu_k^{n-2}\)\nonumber\allowdisplaybreaks\\
&\quad+\smallo\(\zeta_k\(x_k,\mu_k\)+\mu_k^{n-2}\)\nonumber\allowdisplaybreaks\\
&=\Lambda_1\(\S^n,\[g_0\]\)^{-\frac{n-2}{2}}\int_{\B_\xi\(0,\delta/\mu_k\)}U_0\(U_0\check{w}_k+\hat{w}_k\hat{v}_k\)\dy\nonumber\\
&\quad+\bigO\(\sum_{\left|\alpha\right|=4}^{n-4}\left|h_{k,\alpha}\right|^2\mu_k^{2\left|\alpha\right|}+\delta \mu_k^{n-2}\)+\smallo\(\zeta_k\(x_k,\mu_k\)+\mu_k^{n-2}\)\ \ \text{as }k\to\infty.
\end{align}
Moreover, since $d_n\le4$ when $n\le 10$ with equality if and only if $n=10$, by definition of $\zeta_k\(x_k,\mu_k\)$, we obtain
\begin{equation}\label{Th2Eq20}
\sum_{\left|\alpha\right|=4}^{n-4}\left|h_{k,\alpha}\right|^2\mu_k^{2\left|\alpha\right|}=\smallo \(\zeta_k\(x_k,\mu_k\)+\mu_k^{n-2}\)\quad\text{as }k\to\infty.
\end{equation}
It follows from \eqref{Th2Eq19} and \eqref{Th2Eq20} that
\begin{align}\label{Th2Eq21}
\Lambda_1\(\S^n,\[g_0\]\)-\Lambda_1\(\S^n,\[g_k\]\)&=\Lambda_1\(\S^n,\[g_0\]\)^{-\frac{n-2}{2}}\int_{\B_\xi\(0,\delta/\mu_k\)}U_0\(U_0\check{w}_k+\hat{w}_k\hat{v}_k\)\dy\nonumber\\
&\quad+\bigO\( \delta \mu_k^{n-2}\)+\smallo\(\zeta_k\(x_k,\mu_k\)+\mu_k^{n-2}\)\quad\text{as }k\to\infty.
\end{align}
We now estimate the integral in the right-hand side of \eqref{Th2Eq21}, starting with the first term in this integral. We recall that 
$$H_k = \sum_{p=2}^{n-4} H_{k,p},$$ 
where $H_{k,p} = H_{k,x_k,p}$ and $H_{k,x,p}$ is as in \eqref{Lem3Eq8}. By using polar coordinates and the definition of $\check{w}_k$, we obtain
\begin{align}\label{Th2Eq22}
\int_{\B_\xi\(0,\delta/\mu_k\)}U_0^2\check{w}_k\dy
&=c_n\sum_{p,q=2}^{n-4}\mu_k^{p+q}\sum_{a,b,c=1}^n\int_{\S^n}\bigg(\frac{1}{4} \partial_{y_c}\(H_{k,p}\(y\)\)_{ab} \partial_{y_c}\(H_{k,q}\(y\)\)_{ab}\nonumber\\
&\quad-\frac{1}{2} \partial_{y_b}\(H_{k,p}\(y\)\)_{ab} \partial_{y_c}\(H_{k,q}\(y\)\)_{ac}\bigg)\dsigma\(y\)\nonumber\\
&\quad\times \int_0^{\frac{\delta}{\mu_k}}\(\frac{\sqrt{n\(n-2\)}}{1+r^2}\)^{n-2}r^{p+q+n-3}\dr.
\end{align}
To obtain \eqref{Th2Eq22}, we also use that, for each $p,q\in\left\{2,\dotsc,n-4\right\}$,
$$ \int_{\B_{\xi}\(0,\delta/\mu_k\)} \sum_{a,b,c=1}^n U_0^2 \partial_{y_b} \( \(H_{k,p}\(y\)\)_{ab} \partial_{y_c}\(H_{k,q}\(y\)\)_{ac}\) \dy = 0,$$
which follows from an integration by parts together with \eqref{Lem3Eq9} and the fact that $U_0$ is radially symmetric around 0. We observe that, since $H_k \to 0$ in $C^m\(\B_{\xi}\(0, \varepsilon_0\)\)$ as $k\to\infty$ for all $m\in\N$, for each $p,q\in\left\{2,\dotsc,n-4\right\}$ such that $q>d_n$,
\begin{align}\label{Th2Eq23}
\left|H_{k,p}\right|\left|H_{k,q}\right|\mu_k^{p+q}&= \bigO \(\left|H_{k,p}\right|^2\mu_k^{2p +2q+2-n}+\left|H_{k,q}\right|^2\mu_k^{n-2}\)  \nonumber\\
&=\smallo\(\left|H_{k,p}\right|^2\mu_k^{2p}+\mu_k^{n-2}\)\quad\text{as }k\to\infty,
\end{align}
where
$$ \left| H_{k,p} \right|:= \sum_{\left| \alpha\right| = p} \left|h_{k,\alpha}\right|. $$ 
Integrations by parts give that, for each $d\in\left\{4,\dotsc,n-3\right\}$,
\begin{align}\label{Th2Eq24}
&\int_0^{\frac{\varepsilon_0}{\mu_k}}\(\frac{\sqrt{n\(n-2\)}}{1+r^2}\)^{n-2}r^{n+d-3}\dr \nonumber\\
& \quad =\frac{n-2}{d}  \int_0^\infty\(\frac{\sqrt{n\(n-2\)}}{1+r^2}\)^{n-2}\frac{r^2-1}{1+r^2}\,r^{n+d-3}\dr +\bigO\(\mu_k^{n-d-2}\)
\end{align}
and straightforward estimates give that, for each $d\ge n-2$,
\begin{align}\label{Th2Eq25}
&\int_0^{\frac{\delta}{\mu_k}}\(\frac{\sqrt{n\(n-2\)}}{1+r^2}\)^{n-2}r^{n+d-3}\dr\nonumber\\
&\quad=\left\{\begin{aligned}&\(n\(n-2\)\)^{\frac{n-2}{2}}\left|\ln\mu_k\right|+\bigO\(1\)&&\text{if }d=n-2\\&\bigO\(\mu_k^{n-d-2}\)&&\text{if }d>n-2.\end{aligned}\right.
\end{align}
Combining \eqref{Th2Eq22}, \eqref{Th2Eq23}, \eqref{Th2Eq24} and \eqref{Th2Eq25} we obtain 
\begin{align} \label{Th2Eq26}
\int_{\B_\xi\(0,\delta/\mu_k\)}U_0^2\check{w}_k\dy & = 2c_n \( n\(n-2\)\)^{\frac{n-2}{2}}\int_0^{\mu_k}\mu^{-1}F_{k,1}\(\mu\)\dmu\nonumber\\
&\quad+\smallo\(\zeta_k\(x_k,\mu_k\)+\mu_k^{n-2}\)\quad\text{as }k\to\infty,
\end{align}
where
\begin{align*}
F_{k,1}\(\mu\)&:=\frac{n-2}{2}\sum_{p,q=2}^{d_n}\mu^{p+q} \Bigg(\sum_{a,b,c=1}^n\int_{\S^n}\bigg(\frac{1}{4} \partial_{y_c}\(H_{k,p}\(y\)\)_{ab} \partial_{y_c}\(H_{k,q}\(y\)\)_{ab}\nonumber\\
&\quad-\frac{1}{2} \partial_{y_b}\(H_{k,p}\(y\)\)_{ab} \partial_{y_c}\(H_{k,q}\(y\)\)_{ac}\bigg)\dsigma\(y\)\nonumber\\
&\quad\times\left\{\begin{aligned}&
c_{p+q} &&\text{if }p+q<n-2\\&\left|\ln\mu\right|&&\text{if } p = q =d_n\end{aligned}\right\} \Bigg)
\end{align*}
and
$$ c_{p+q} :=    \int_0^\infty\(\frac{1}{1+r^2}\)^{n-2}\frac{r^2-1}{1+r^2}\,r^{n+p+q-3}\dr.$$
As regards the second term in the integral in the right-hand side of \eqref{Th2Eq21},  by using again polar coordinates, we obtain
\begin{equation}\label{Th2Eq27}
\int_{\B_\xi\(0,\delta/\mu_k\)}U_0\hat{w}_k\hat{v}_k \dy=c_n^2 \sum_{p,q=2}^{n-4}\mu_k^{p+q}\int_{\B_\xi\(0,\delta/\mu_k\)}U_0w_{k,p}v_{k,q}\dy .
\end{equation}
By using \eqref{Lem5Eq9} and since $H_k \to 0$ in $C^m\(\B_{\xi}\(0, \varepsilon_0\)\)$ as $k\to\infty$ for all $m\in\N$, we obtain that, for each $p,q\in\left\{2,\dotsc,n-4\right\}$ such that $d:= p+q \in\left\{4,\dotsc,n-3\right\}$,
\begin{equation} \label{Th2Eq28}
\int_{\B_\xi\(0,\delta/\mu_k\)}U_0w_{k,p}v_{k,q}\dy=\int_{\R^n}U_0w_{k,p}v_{k,q} \dy+\smallo\(\mu_k^{n-d-2}\)\quad\text{as }k\to\infty.
\end{equation}
We now consider the case where $p,q\in\left\{2,\dotsc,n-4\right\}$ are such that $d:= p+q \ge n-2$. We recall that by Lemma~\ref{Lem4}, 
$$v_{k,q}\(y\)= \frac{P_{k,q}\(y\)}{\big(1+\left|y\right|^2\big)^{\frac{n}{2}}}\quad\forall y\in\R^n,$$ 
where $P_{k,q}$ is the polynomial of degree $q+2$ given by 
$$P_{k,q}\(y\):= \sum_{m=0}^{\[\(d-2\)/2\]} \sum_{l=0}^{m+2}\gamma_{k,x_k,d,l,m}\left|y\right|^{2l}\Delta_\xi^m w_{k,x_k,d}\(y\)\quad\forall y\in\R^n.$$ 
We let $P_{k,q}^{\(q+2\)}$ be the sum of terms of highest degree in $P_{k,q}$, i.e.
$$ P_{k,q}^{\(q+2\)}\(y\) =  \sum_{m=0}^{\[\(q-2\)/2\]}\gamma_{k,x_k,q,m+2,m}\left|y\right|^{2m+4} \Delta_\xi^m w_{k,q}\(y\)\quad\forall y\in\R^n.$$ 
Since $H_k \to 0$ in $C^m\(\B_{\xi}\(0, \varepsilon_0\)\)$ as $k\to\infty$ for all $m\in\N$, straightforward estimates using polar coordinates then give
\begin{align}\label{Th2Eq29}
&\int_{\B_\xi\(0,\delta/\mu_k\)}U_0w_{k,p}v_{k,q}\dy\nonumber\\
&\quad=\left\{\begin{aligned}&\(n\(n-2\)\)^{\frac{n-2}{2}}\left|\ln\mu_k\right| \int_{\S^{n-1}} w_{k,p}\(y\) P_{k,q}^{\(q+2\)}\(y\) \dsigma\(y\)+\smallo\(1\)&&\text{if }d=n-2\\&\smallo\(\mu_k^{n-d-2}\)&&\text{if }d>n-2\end{aligned}\right.
\end{align}
as $k\to\infty$. By putting together \eqref{Th2Eq23}, \eqref{Th2Eq27}, \eqref{Th2Eq28} and \eqref{Th2Eq29},
we obtain
\begin{align}\label{Th2Eq30}
\int_{\B_\xi\(0,\delta/\mu_k\)} U_0\hat{w}_k\hat{v}_k\dy &= 2c_n \( n\(n-2\)\)^{\frac{n-2}{2}}\int_0^{\mu_k}\mu^{-1}F_{k,2}\(\mu\)\dmu\nonumber\\
&\quad+\smallo\(\zeta_k\(x_k,\mu_k\)+\mu_k^{n-2}\)\quad\text{as }k\to\infty,
\end{align}
where
\begin{align*}
F_{k,2}\(\mu\)&:= c_n \sum_{p,q=2}^{d_n}\mu^{p+q}\\
&\quad\times\left\{\begin{aligned}& p\int_{\R^n}\big(1+\left|y\right|^2\big)^{2-n}w_{k, p}\(y\)v_{k,q}\(y\)\dy&&\text{if } p+q<n-2\\& d_n \left|\ln\mu\right| \int_{\S^{n-1}} w_{k,d_n}\(y\) P_{k,d_n}^{\(d_n+2\)}\(y\) \dsigma\(y\)
&&\text{if } p=q=d_n.\end{aligned}\right.
\end{align*}
The functions $F_{1,k}$ and $F_{2,k}$ have been investigated in \cite{KMS}. In particular, since $n \le 10$, Proposition~A.4 of \cite{KMS} applies\footnote{Our functions $F_{1,k}\(\mu\)$ and $F_{2,k}$ coincide with the functions $I_{1,\mu}^{\(n\)}\(H_k,H_k\)$ and $I_{2,\mu}^{\(n\)}\(H_k,H_k\)$ introduced in \cite{KMS}*{Appendix~A}. Our definition of $F_{2,k}$ seems to introduce a factor $-c_n$ with respect to $I_{2,\mu}^{\(n\)}\(H_k,H_k\)$ in \cite{KMS}, but this is a choice of convention: the function $Z\big(H_k^{\(p\)}\big)$ involved in the definition of $I_{2,\mu}^{\(n\)}\(H_k,H_k\)$ in \cite{KMS} is equal to $-c_n w_{k,p}$ where $w_{k,p}$ is as in Lemma~\ref{Lem4}.} and gives that 
there exists $C_n^{\prime}>0$ depending only on $n$ such that 
\begin{equation}\label{Th2Eq31}
F_{1,k}\(\mu\) + F_{2,k}\(\mu\)\ge C_n^{\prime}\zeta_k\(x_k,\mu\)\quad\forall\mu>0.
\end{equation}
It follows from \eqref{Th2Eq26}, \eqref{Th2Eq30} and \eqref{Th2Eq31} that
\begin{equation}\label{Th2Eq32}
\int_{\B_\xi\(0,\varepsilon_0/\mu_k\)}U_0\(U_0\check{w}_k+\hat{w}_k\hat{v}_k\)\dy\ge C_n^{\prime\prime}\zeta_k\(x_k,\mu_k\)+\smallo\(\mu_k^{n-2}\)\quad\text{as }k\to\infty
\end{equation}
for some $C_n^{\prime\prime}>0$ depending only on $n$. By plugging \eqref{Th2Eq32} into \eqref{Th2Eq21}, we finally obtain
$$
\begin{aligned}
\Lambda_1\(\S^n,\[g_0\]\)-\Lambda_1\(\S^n,\[g_k\]\)&\ge C_n\zeta_k\(x_k,\mu_k\)+\bigO\( \delta \mu_k^{n-2}\)\\
&\quad+\smallo\(\zeta_k\(x_k,\mu_k\)+\mu_k^{n-2}\)\quad\text{as }k\to\infty,
\end{aligned} $$
where $C_n:= \Lambda_1\(\S^n,\[g_0\]\)^{-\frac{n-2}{2}}C_n^{\prime\prime}>0$ and the constant in the term $\bigO\( \delta \mu_k^{n-2}\)$ is independent of $\delta \in (0, \varepsilon/2)$. A simple diagonal argument now allows us, up to a subsequence, to replace $\delta$ by a sequence $\delta_k \to 0$. This proves \eqref{Th2Eq4} and completes the proof of Theorem~\ref{Th1}.
\endproof

\end{document}